\documentclass{article}
\usepackage{graphicx}
\graphicspath{images/}
\usepackage{stackrel}
\usepackage[utf8]{inputenc}
\usepackage{latexsym, amssymb, latexsym, amsfonts, amsmath, esint, graphicx, caption, subcaption, dsfont, units, fullpage, float, setspace, cases, polynom, multirow , stmaryrd , textcomp, tikz, mathtools, amsthm, epstopdf, epsfig, hyperref, mathrsfs}
\usepackage[first=1, last=10]{lcg}
\usepackage{listings,enumerate,geometry,csquotes,url}
\usepackage{etoolbox}
\usepackage{pdfpages,xspace,xparse}
\usepackage{import,pst-plot}
\usepackage[numbers,sort]{natbib}
\usepackage[toc,page]{appendix}
\usepackage[mathscr]{euscript}
\usepackage{nicematrix}
\usepackage{enumitem}
\usepackage{comment}
\usepackage{refcount}
\DeclareMathAlphabet{\mathpzc}{OT1}{pzc}{m}{it}
\newcommand{\kell}{\mathpzc{k}}
\usepackage{stmaryrd}
\SetSymbolFont{stmry}{bold}{U}{stmry}{m}{n}
\DeclareFontFamily{U}{calligra}{}
\DeclareFontShape{U}{calligra}{m}{n}{<->callig15}{}
\newcommand{\inorm}[2]{\left\| #1 \right\|_{#2}}
\newcommand{\Eint}{\Xi(k)}

\newcommand{\D}{\mathcal{D}}
\newcommand{\enumref}[2]{\ref{#1}\textcolor{blue}{.}\ref{#2}}
\newcommand{\Omegadef}{\{k \in \C: |k|> r \text{ and } \pi/4 < \arg(k) < 3\pi/4\}}
\newcommand{\Omegaextdef}{\{k\in\C: |k|>r \text{ and } \theta_0< \arg(k) < \pi - \theta_0\}}
\newcommand{\Omegaext}{{\Omega_\mathrm{ext}}}

\newcommand{\Etwo}{\etwo{\frac12 \left\| \frac{(\beta\mathfrak n)'}{(\beta\mathfrak n)} \right\|_\D}}
\newcommand{\fracbetanu}{\frac{(\beta\mathfrak n)'}{(\beta \mathfrak n)}}
\newcommand{\fracbetanuargs}[1]{\frac{(\beta\mathfrak n)'(k,#1)}{(\beta \mathfrak n)(k,#1)}}

\newcommand{\betanuargs}[1]{(\beta\mathfrak n)(k,#1)}
\newcommand{\bigoh}{{O}}
\newcommand{\littleoh}{{o}}
\newcommand{\beq}{\begin{equation}}
\newcommand{\eeq}{\end{equation}}
\newcommand{\ba}{\begin{array}}
\newcommand{\ea}{\end{array}}
\newcommand{\bea}{\begin{eqnarray*}}
\newcommand{\eea}{\end{eqnarray*}}
\newcommand{\bc}{\begin{center}}
\newcommand{\ec}{\end{center}}
\newcommand{\bt}{\begin{table}}
\newcommand{\et}{\end{table}}

\newcommand{\sgn}{{\rm{sgn}}}
\newcommand{\la}[1]{\label{#1}}

\newcommand{\no}{\noindent}

\newcommand{\rf}[1]{(\ref{#1})}
\newcommand{\beqno}{\begin{displaymath}}
\newcommand{\eeqno}{\end{displaymath}}

\newcommand{\ie}{{\em i.e.,~}}
\newcommand{\eg}{{\em e.g.,~}}
\newcommand{\been}{\begin{enumerate}}
\newcommand{\een}{\end{enumerate}}

\newcommand{\erf}{{\rm erf}}

\newcommand{\etc}{{\em etc.}}

\newcommand{\C}{\mathbb{C}}
\newcommand{\R}{\mathbb{R}}

\newcommand{\Res}{\mathrm{Res}}
\renewcommand{\Re}{\mathrm{Re}}

\newcommand{\eone}[1]{\ensuremath{e^{#1}}}
\newcommand{\etwo}[1]{\ensuremath{\exp\left(#1\right)}}
\newcommand{\then}{\qquad \Rightarrow \qquad}
\renewcommand{\And}{\qquad \text{ and } \qquad}
\renewcommand{\bar}{\overline}

\newcommand{\figcl}[4]{\begin{figure}[H]\begin{center} \includegraphics[scale=#1]{#2} \end{center} \caption{#3} \label{fig:#4} \end{figure}}

\newcommand{\case}[2]{\left\{ \hspace*{-0.15in} \arraycolsep=0.15in\def\arraystretch{#1}\begin{array}{ll} #2 \end{array}\right.}

\newlength{\myheight}
\newlength{\mylength}

\newcounter{saveeqn}

\newtheorem{theorem}{Theorem}
\newtheorem{lemma}[theorem]{Lemma}

\newtheorem{definition}[theorem]{Definition}

\newtheorem{rem}[theorem]{Remark}

\newtheorem{assumption}[theorem]{Assumption}

\def\XXint#1#2#3{{\setbox0=\hbox{$#1{#2#3}{\int}$}
     \vcenter{\hbox{$#2#3$}}\kern-.5\wd0}}

\hypersetup{
    colorlinks=true,
    linkcolor=blue,
    citecolor=red
}
\numberwithin{equation}{section}
\usepackage{subfiles}
\title{The explicit solution of linear, dissipative, second-order initial-boundary value problems with variable coefficients}
\author{
Matthew Farkas, Bernard Deconinck\\
Department of Applied Mathematics\\
University of Washington\\
Seattle, WA 98195-2420
}
\date{\today}

\begin{document}

\maketitle

\begin{abstract}
We derive explicit solution representations for linear, dissipative, second-order Initial-Boundary Value Problems (IBVPs) with coefficients that are spatially varying, with linear, constant-coefficient, two-point boundary conditions. We accomplish this by considering the variable-coefficient problem as the limit of a constant-coefficient interface problem, previously solved using the Unified Transform Method of Fokas. Our method produces an explicit representation of the solution, allowing us to determine properties of the solution directly. As explicit examples, we demonstrate the solution procedure for different IBVPs of variations of the heat equation, and the linearized complex Ginzburg-Landau (CGL) equation (periodic boundary conditions). We can use this to find the eigenvalues of dissipative second-order linear operators (including non-self-adjoint ones) as roots of a transcendental function, and we can write their eigenfunctions explicitly in terms of the eigenvalues.
\end{abstract}
\setcounter{tocdepth}{4}
\setcounter{secnumdepth}{4}
\tableofcontents
%
%
\section{Introduction}
The Unified Transform Method (UTM), or Method of Fokas, is used to solve Initial Value Problems (IVPs) and Initial-Boundary Value Problems (IBVPs) for integrable equations. Its application to linear, constant-coefficient partial differential equations (PDEs) is particularly convenient and straightforward. The UTM leads to many new insights on PDEs and IBVPs, see for instance \cite{world_scientific,JC_fokas_collab, JC_fokas_book, treharne, JC_bernard_fokas,JC_flyer,JC_fokas_paper,DAS_diagonalization}, and references therein. Especially relevant for us, the method has been used to explicitly solve interface problems and problems with piecewise-constant coefficients, see \cite{Mantzavinos,interface_heat, interface_dispersive,interface_heat_ring,interface_kdv,interface_maps,interface_schrodinger}. The purpose of this paper is to generalize the UTM to solve variable-coefficient IBVPs. In \cite{treharne}, Fokas and Treharne use a Lax Pair approach to analyze specific variable-coefficient IBVPs. Their approach reduces the problem from solving a {\em Partial} Differential Equation to solving an {\em Ordinary} Differential Equation (ODE) by writing the solution of the PDE as an integral over the solutions to a non-autonomous ODE, but it does not provide an explicit representation of the solution. This approach, like separation of variables, is useful if the associated ODE is a second-order, self-adjoint problem on a finite domain, for which we have regular Sturm-Liouville theory \cite{ODEs}, but does not generalize well to problems that are not self adjoint, of higher order, or posed on an unbounded domain.

In our approach to variable-coefficient IBVPs, we divide the domain into $N$ parts and approximate the equation by a constant-coefficient equation on each part. We solve the resulting interface problem using the UTM as shown in \cite{Mantzavinos,interface_heat,interface_heat_ring,interface_kdv,interface_maps,interface_dispersive,interface_schrodinger}. Using Cramer's rule, the solution in each part is found as a ratio of determinants. Through the nontrivial steps of obtaining an explicit expression for the determinants and taking the limit as $N$ goes to infinity, a complicated but explicit solution expression is obtained. As in previous applications of the UTM (\eg \cite{JC_fokas_book,DAS_diagonalization} for constant-coefficient problems, \cite{interface_kdv,smith_2pointBCs} for interface problems), one of the benefits of our approach is characterizing which boundary conditions give rise to a well-posed IBVP. In particular, for the finite-interval problem, this work is consistent with Locker's work on Birkhoff regularity, \eg \cite{locker}. Since the UTM is generalizable to large classes of varying boundary conditions, IBVPs of higher order, including non-self-adjoint problems, we expect our method to generalize in these same directions as well.

In this manuscript, expanding on work presented in \cite{Farkas_Deconinck}, we construct explicit solution expressions for general, second-order IBVPs with spatially-varying coefficients and with linear, two-point boundary conditions, as integrals over known quantities. We present the solutions for the whole-line problem, the half-line problem, and the finite-interval problem. We choose to demonstrate the solution process to the three problems in all their generality starting with the simplest, so that we start with the whole-line problem in Section~\ref{sec:statement_wl}, followed by the half-line problem in Section~\ref{sec:statement_hl}, finishing with the finite-interval problem in Section~\ref{sec:statement_fi}. In Sections~\ref{sec:example_wl},~\ref{sec:example_hl},~\ref{sec:example_fi2},~and~\ref{sec:example_fi1}, we restrict to specific examples. For the finite-interval problems, our explicit representation characterizes the eigenvalues of the ODE obtained through separation of variables and gives the eigenfunctions explicitly in terms of these eigenvalues, see Section~\ref{sec:eigenvalues}. This allows for the numerical approximation of the eigenvalues, including for non-self-adjoint problems. Other numerical applications are presented in \cite{Farkas_Deconinck}. In Appendix~\ref{sec:derivations}, we present a formal derivation of the solution. In this appendix we switch the order of exposition by deriving the solution to the finite-interval problem first, as the solution of the other two problems follows from it. We finish the paper with rigorous proofs in Appendices~\ref{sec:proofs_welldefined}--\ref{sec:proofs_IC}.

Our formulae may seem complicated; however, they are similar to the solutions found in \cite{PT}, which have been used to prove a variety of properties of solutions to ODEs and eigenvalue problems. Indeed, our notation is inspired by this book. While our solutions are similar, our methods are entirely different. The reader may also find our expressions reminiscent of path integrals \cite{tak}, although those are usually used to propagate in time, unlike our spatial ``discretization'' approach. 
\section{Assumptions and definitions} \label{sec:assumptions}
\begin{figure}[tb]
    \centering
    \includegraphics[scale=0.5]{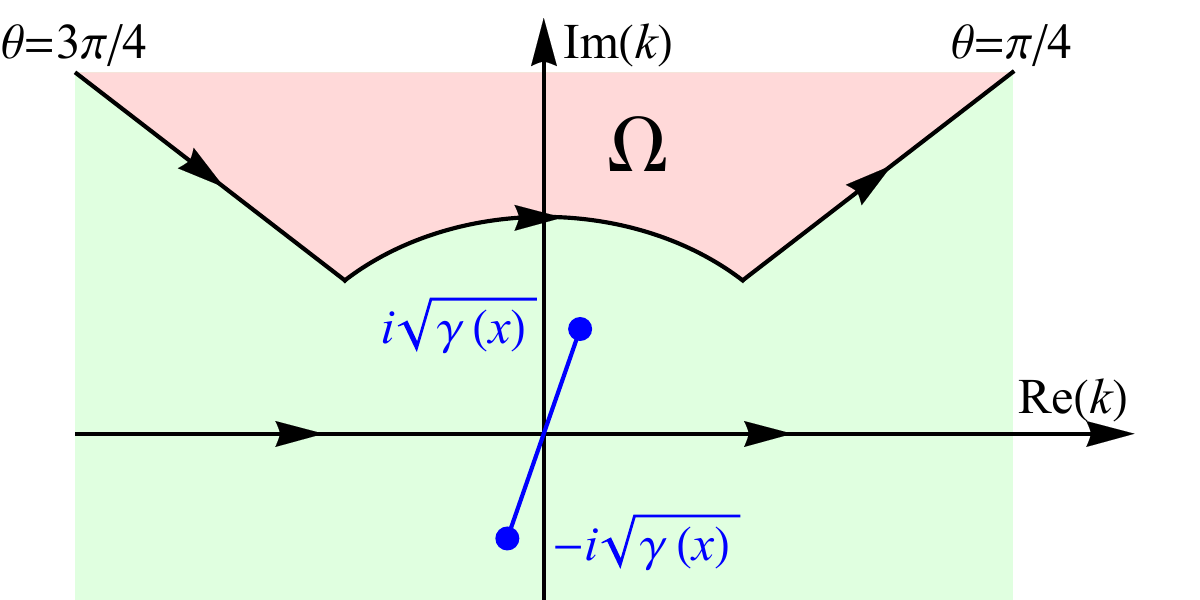}
    \caption{The region $\Omega$ with the branch cuts for $\mathfrak g(k,x)$.} 
    \label{fig:Omega}
\end{figure}

Throughout this paper, we consider the linear, second-order evolution equation with spatially variable coefficients:
\begin{subequations} \label{eqn:IVP}
    \begin{align} 
        \label{eqn:evolution_equation}
        q_t &= \alpha(x)\left(\beta(x) q_{x}\right)_x+\gamma(x)q+f(x,t), && x\in \D\subset \R, \quad t>0, \\
        q(x,0) &= q_0(x), && x \in \D, 
    \end{align}
\end{subequations}
on different domains $\D$ with (possibly) some functions $f_0(t)$, $f_1(t)$ prescribed at the boundary of $\D$. In all cases, the solution is written as 
\begin{align} \label{eqn:q_all}
    q(x,t) = \frac{1}{2\pi} \int_{\partial \Omega} \frac{\Phi(k,x,t)}{\Delta (k)} e^{-k^2t} \, dk,
\end{align}
where the functions $\Phi(k,x,t)$ and $\Delta(k)$ depend on $\D$ and the initial and boundary conditions provided. The region \mbox{$\Omega = \Omegadef$}, for some $r>\sqrt{M_\gamma}$, where $M_\gamma = \|\gamma\|_\infty$, as shown in Figure~\ref{fig:Omega}. 
In this section, we establish notation and introduce assumptions on the functions in \eqref{eqn:IVP} that we use throughout the paper.

We define \mbox{$\arg(\,\cdot\,) \in \left[-\pi/2, 3\pi/2\right)$}. We use $\D$ to denote the domain of the problem, so that $\D = \mathbb R$, $\D = (x_l,\infty)$, and $\D = (x_l,x_r)$ for the whole-line, half-line, and the finite-interval problems, respectively. The domain $\D$ is given by the open set, and we denote the closure by $\overline{\D}$. We write the $L^1$--norm over the domain $\D\subseteq \R$ as $\|\cdot\|_{\D}$. When used on a function of multiple variables, we implicitly assume a supremum norm on the other variables, \eg for a function $f(k,x)$ for $k\in \Omega\subseteq \C$ and $x\in\D$, 
$$ \| f\|_\D = \sup_{k\in \Omega} \int_\D |f(k,x)| \, dx. $$
In this way, the norms always represent fixed numbers, never functions. The notation $\mathrm{AC}(\,\cdot\,)$ represents the space of locally absolutely continuous functions on the closure of the domain. We use the `big-oh' notation $\bigoh(\,\cdot\,)$ and the `little-oh' notation $\littleoh(\,\cdot\,)$, as described in \cite{miller}.
\begin{assumption} \label{ass:alphabetagamma}
    We assume the following about the coefficient functions $\alpha,\beta,\gamma$:
    \begin{enumerate}
        \item \label{enum:arg} $\sup_{x\in\D} |\arg(\alpha(x)\beta(x)|<\pi/2$, 
        \item \label{enum:AC} $\alpha, \beta \in \mathrm{AC}({\D})$,
        \item \label{enum:lowerbound} $m_{\alpha\beta} =\inf_{x\in D}|\alpha(x)\beta(x)| >0$,
        \item \label{enum:bounded} $\alpha\beta$, ${\gamma} \in L^\infty(\D)$, and we define $M_{\alpha\beta} = \|\alpha\beta\|_\infty$ and $M_\gamma = \|\gamma\|_\infty$, 
        \item \label{enum:L1} $(\beta'/\beta - \alpha'/\alpha)$, $\gamma' \in L^1(\D)$.
    \end{enumerate}
\end{assumption}
\begin{assumption} \label{ass:ffunctions} 
    We assume the following about the inhomogeneous, initial, and boundary functions $f,q_0,f_m$:
    \begin{enumerate}
        \item \label{enum:finhomogeneous} For the inhomogeneous function $f(x,t)$, we assume 
        $f(x,\cdot) \in \mathrm{AC}((0,T))$ for each $x\in \bar \D$, and
        \begin{align*}
            \|f\|_{\D} =  \sup_{t\in[0,T]} \int_\D |f(x,t)| \, dx < \infty
            \qquad \text{ and } \qquad
            \|f_{t}\|_{\D} = \sup_{t\in[0,T]} \int_\D |f_{t}(x,t)| \, dx < \infty.
        \end{align*}
        \item \label{enum:initial} For the initial condition $q_0(x)$, we assume $q_0\in L^1(\D)$.
        \item \label{enum:fboundary} For the boundary functions $f_m(t)$, $m=0,1$, we assume $f_m \in \mathrm{AC}((0,T))$ and $f_m'\in L^\infty ((0,T))$.
    \end{enumerate}
\end{assumption}
\begin{assumption} \label{ass:alphabetagamma2}
    The finite-interval IBVP has a number of subcases. One of these, Case~\ref{enum:BC4} (see Definition~\ref{def:boundarycases}), requires the following assumption, in addition to Assumption~\ref{ass:alphabetagamma}:
    \begin{enumerate}
        \item \label{enum:AC_extra} $\beta'/\beta - \alpha'/\alpha \in \mathrm{AC}(\D)$.
        \item \label{enum:fm_extra} For the boundary functions $f_m(t)$, $m=0,1$, we require $f_m'\in \mathrm{AC}((0,T))$ and $f_m'' \in L^\infty((0,T))$.
    \end{enumerate}
\end{assumption}
\begin{rem} \label{rem:assumptions} $ $
    \begin{itemize}
        \item 
        Assumption~\enumref{ass:alphabetagamma}{enum:arg} is equivalent to restricting to problems that we call {\em fully dissipative}. This is in contrast to problems that we call {\em partially dissipative} and/or {\em partially dispersive} where $\sup_{x\in \D} |\arg(\alpha(x)\beta(x))| = \pi/2$, but $\inf_{x\in\D} |\arg(\alpha(x)\beta(x))| < \pi/2$ or {\em fully dispersive} where $|\arg(\alpha(x)\beta(x))| = \pi/2$ for all $x\in \D$.
        \item
        Assumption~\enumref{ass:alphabetagamma}{enum:AC} may seem odd considering that we derive our results from those for an interface problem in Appendix~\ref{sec:derivations}. However, in that section, we use the mean value theorem as we let the limit of the number of interfaces $N$ approach infinity, and thus we assume continuity of our function. This section can be amended to include piecewise continuous functions, but makes the solution formulas even more complicated. For simplicity, we restrict to continuous functions. Alternatively, we could employ distribution theory to extend the current results to discontinuous functions.
        \item 
        Assumptions~\enumref{ass:alphabetagamma}{enum:lowerbound}~and~\enumref{ass:alphabetagamma}{enum:bounded} are physically natural conditions to impose. 
        Assumption~\enumref{ass:alphabetagamma}{enum:L1} ensures that our solution is well defined. It may be possible to extend this to other $L^p$ spaces or other more general spaces with some more work.
        \item
        Throughout this paper, we always use Assumptions~\ref{ass:alphabetagamma}~and~\ref{ass:ffunctions}. We clearly state when Assumption~\ref{ass:alphabetagamma2} is used, which is only for the \ref{enum:BC4}th Boundary Case of the finite-interval.
    \end{itemize}
\end{rem}
\begin{definition} \label{def:mu}
    For $x\in \D$, since $\alpha(x)$ and $\beta(x)$ are continuous by Assumption~\enumref{ass:alphabetagamma}{enum:AC}, we define arguments of $\alpha(x)$ and $\beta(x)$ to be $\theta_\alpha(x)$ and $\theta_\beta(x)$, chosen to be continuous\footnote{Note that not necessarily $\theta_\alpha(x)=\arg(\alpha(x))$, given how we defined the range of $\arg(\,\cdot\,)$ above, because of the continuity requirement. For instance, if $\alpha(x)=\exp(ix)$ (and say $\beta(x)=\exp(-ix)$), we can set $\theta_\alpha(x)=x$.}, so that
    \begin{align} \label{eqn:alpha,beta_branches}
        \alpha(x) = |\alpha(x)| e^{i\theta_\alpha(x)} 
        \qquad \text{ and } \qquad
        \beta(x) = |\beta(x)| e^{i\theta_\beta(x)}.
    \end{align}
    Using this, we define, for $x\in \D$, 
    \begin{align} \label{mu}
        \mu(x) &= \frac{1}{\sqrt{\alpha(x)\beta(x)}} 
        = \frac{1}{\sqrt{|\alpha(x)\beta(x)|}} e^{-\frac{i}{2}(\theta_\alpha(x) + \theta_\beta(x))},
    \end{align}
    and, for $x\in \D$ and $k\in \C$,
    \begin{align} \label{eqn:gfrak}
        \mathfrak g(k,x) 
        &= \sqrt{1 + \frac{\gamma(x)}{k^2}}
        =\sqrt{\left| 1 + \frac{\gamma(x)}{k^2} \right|} \eone{ \frac{i}{2} \arg( 1+ \gamma(x)/k^2) } \!.
    \end{align}
    We also define $\mathfrak n(k,x) = \mu(x) \mathfrak g(k,x)$, $(\beta\mu)(x) = \beta(x) \mu(x)$, and $(\beta\mathfrak n)(k,x) = \beta(x)\mathfrak n(k,x)$, 
    \begin{align}
        \sqrt{(\beta \mu)(x)} = \sqrt[4]{\left| \frac{\beta(x)}{\alpha(x)} \right|} e^{\frac i4 ( \theta_\beta(x) - \theta_\alpha(x))},
        \qquad \text{} \qquad
        \sqrt{\mathfrak g(k,x)} = \sqrt[4]{\left| 1 + \frac{\gamma(x)}{k^2} \right|} \eone{ \frac{i}{4} \arg( 1+ \gamma(x)/k^2) }, 
    \end{align}
    and $\sqrt{(\beta \mathfrak n)(k,x)} = \sqrt{(\beta\mu)(x)} \sqrt{\mathfrak g(k,x)}$. Finally, we define
    \begin{align} \label{eqn:mathfrak_u_def}
        \mathfrak u(x) &= \frac{1}{\mu(x)} \left( \frac{\beta'(x)}{\beta(x)} - \frac{\alpha'(x)}{\alpha(x)} \right)\!,
    \end{align}
    and
    \begin{align} \label{eqn:q_alpha,f_alpha,ftilde_alpha}
        q_\alpha(x) = \frac{q_0(x)}{\alpha(x)}, 
        \qquad 
        f_\alpha(x,t) = \frac{f(x,t)}{\alpha(x)}, 
        \qquad
        \tilde f_\alpha(k^2,x,t) = \int_0^t f_\alpha(x,s) e^{k^2s} \, ds,
    \end{align}
    and $\psi_\alpha(k^2,x,t) = q_\alpha(x) + \tilde f_\alpha(k^2,x,t)$.
\end{definition}
\section{The whole-line problem} \label{sec:statement_wl}
Consider \eqref{eqn:IVP} for $x\in \mathbb{R}$ and with decay at infinity, 
\begin{subequations} \label{eqn:IBVP_wl}
    \begin{align}
        \label{eqn:PDE_wl}
        q_t &= \alpha(x)\left(\beta(x) q_{x}\right)_x+\gamma(x)q+f(x,t), && x\in \R, \quad t>0, \\
        \label{eqn:IC_wl}
        q(x,0)&=q_0(x), && x\in \R, \\ 
        \label{eqn:BC_wl}
        \lim_{|x|\to\infty} q(x,t) &= 0, && t>0.
    \end{align}
\end{subequations}
\begin{theorem}\label{theo1}
    \no Under Assumptions~\ref{ass:alphabetagamma}~and~\ref{ass:ffunctions}, the IVP \rf{eqn:IBVP_wl} has the solution
    \begin{align} \label{eqn:q_wl}
        q(x,t) = \frac{1}{2\pi} \int_{\partial \Omega} \frac{\Phi(k,x,t)}{\Delta (k)} e^{-k^2t} \, dk,
    \end{align}
    where $\Omega$ is shown in Figure~\ref{fig:Omega}. Here 
    \begin{align} \label{eqn:phidelta}
        \Phi(k,x,t) = \int_{-\infty}^{\infty} \frac{\Psi(k,x,y)\psi_\alpha (k^2,y,t)}{\sqrt{(\beta \mathfrak n)(k,x)} \sqrt{(\beta \mathfrak n)(k,y)}} \, dy
        \qquad \text{ and } \qquad 
        \Delta(k) = \sum_{n=0}^\infty \mathcal E_{2n}^{(-\infty,\infty)}(k), 
    \end{align}
    \no with, for $y<x$,
    \beq \label{eqn:Psi_wl}
        \Psi(k,x,y) = \exp\left( ik \int_y^x \mathfrak n(k,\xi) \, d\xi\right) \sum_{n=0}^\infty \sum_{\ell=0}^n (-1)^\ell  \tilde{\mathcal E}_{n-\ell}^{(-\infty,y)} (k) \mathcal E_\ell^{(x,\infty)}(k), 
    \eeq
    \no and for $y>x$, $\Psi(k,x,y) =\Psi(k,y,x)$. Here, $\mathcal E_0^{(a,b)}(k) = 1$, $\tilde{\mathcal E}_0^{(a,b)}(k) = 1$, and for $n\geq 1$,
    \begin{subequations} \label{eqn:En&Entilde}
        \begin{align}
            \label{eqn:En}
            \mathcal E_n^{(a,b)}(k) &= \frac{1}{2^{n}}\int_{a=y_0<y_1<\cdots < y_n<y_{n+1}=b}  \left(\prod_{p=1}^n \fracbetanuargs{y_p} \right) \exp\left( ik \sum_{p=0}^{n} (1-(-1)^{n-p}) \int_{y_p}^{y_{p+1}} \mathfrak n(k,\xi) d\xi\right) d\mathbf y_n , \\
            \label{eqn:Entilde}
            \tilde{\mathcal E}_n^{(a,b)}(k) &= \frac{1}{2^{n}}\int_{a=y_0<y_1<\cdots < y_n<y_{n+1}=b} \left(\prod_{p=1}^n \fracbetanuargs{y_p} \right) \exp\left( ik \sum_{p=0}^{n} (1-(-1)^{p}) \int_{y_p}^{y_{p+1}} \mathfrak n(k,\xi) d\xi\right) d\mathbf y_n,
        \end{align}
    \end{subequations}
    where $d\mathbf y_n = dy_1\cdots dy_n$ and the prime denotes the derivative with respect to the second variable.
    The function $\mathcal E_n^{(a,b)}(k)$ is defined for $b=\infty$ and if $n$ is even, for $a=-\infty$. The function $\tilde{\mathcal E}_n^{(a,b)}(k)$ is defined for $a=-\infty$ and if $n$ is even, for $b=\infty$.
    The functions $\psi_\alpha(k^2,x,t)$, $\mathfrak n(k,x)$, and $(\beta\mathfrak n)(k,x)$ are given in Definition~\ref{def:mu}.
\end{theorem}
\begin{proof}
    The formal derivation is given in Appendix~\ref{sec:derivations}. Its validity is proven in Appendices~\ref{sec:proofs_welldefined}--\ref{sec:proofs_IC}.
\end{proof}
\subsection{Example: The partially lumped heat equation} \label{sec:example_wl}
\begin{figure}[tb]
    \centering
    \includegraphics[scale=0.5]{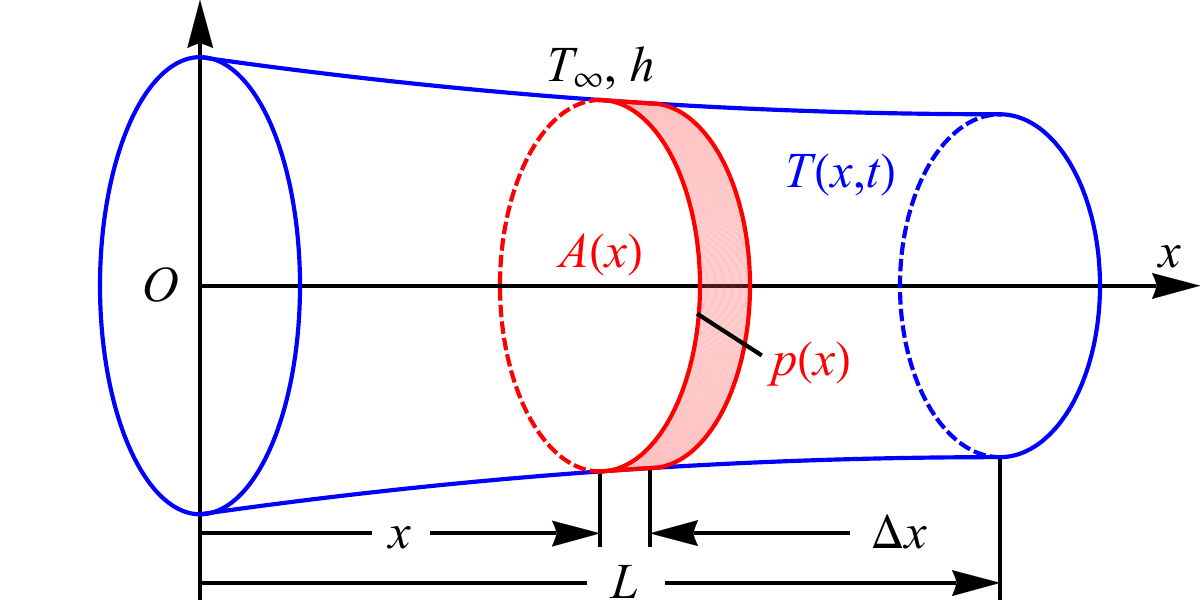}
    \caption{Terminology for the derivation of the partially lumped heat equation \cite{heat_ozosik}.} 
    \label{fig:heat_lumped}
\end{figure}
Consider the heat equation with partial lumping analysis \cite{heat_ozosik}, describing the temperature $T(x,t)$ in a body  with minimal temperature variation in the $y$ and $z$ directions with ambient temperature $T_\infty$, heat transfer coefficient $h_0$, thermal conductivity $k_0$, cross-sectional area $A(x)$, and perimeter $p(x)$, see Figure~\ref{fig:heat_lumped}. We assume the length $L$ is much greater than the width in the $y$ and $z$ directions. Ignoring temperature deviations in the $y$ and $z$-directions, this IBVP takes the form
\begin{subequations} \label{eqn:IBVP_schrodinger}
    \begin{align}
        \theta_t &= \frac{1}{A(x)} \big( A(x)\theta_x\big)_x  - C(x) \theta, && x\in \R, \quad t>0, \\ 
        \theta(x,0)&=\theta_0(x), && x\in \R, \\ 
        \lim_{|x|\to\infty} \theta(x,t) &= 0, && t>0.
    \end{align}
\end{subequations}
Here $\theta(x,t) = T(x,t) - T_\infty$ represents the difference of the temperature in the body $T(x,t)$ and the ambient temperature $T_\infty$, the function $C(x) = h_0 p(x)/(k_0 A(x)) > 0$,
and we equate the thermal diffusivity to 1 ($\alpha=1$).
Comparing this to \eqref{eqn:IBVP_wl}, we have $\alpha(x)=1/A(x)$, $\beta(x)=A(x)$, $\gamma(x) = -C(x)$, and $f(x,t)\equiv 0$. 
We require the absolute continuity of $A(x)>0$, the boundedness of $C(x)$, and the absolute integrability of $A'(x)/A(x)$ and $C'(x)$, so that Assumption~\ref{ass:alphabetagamma} is satisfied. Then we have the solution
\beq \label{eqn:theta_sol}
    \theta(x,t) = \frac{1}{2\pi} \int_{\partial \Omega}  \frac{\Phi(k,x,t)}{\Delta (k)} e^{-k^2t} \, dk,
\eeq
\no where $\Omega$ is shown in Figure~\ref{fig:Omega},
\beq
    \mathfrak n(k,x) = \sqrt{1 - \frac{C(x)}{k^2}},
\eeq
\no and $\psi_\alpha(k^2,x,t) = A(x)\theta_0(x)$. The functions $\Phi(k,x,t)$ and $\Delta(k)$ are given in \eqref{eqn:phidelta}.
\subsection{A note about the integrability conditions.}
A variable coefficient PDE in the form
\beq
    q_t = a(x) q_{xx} + b(x) q_x + c(x) q,
\eeq
can always be written in the form of \eqref{eqn:PDE_wl} as
\beq
    q_t = a(x) \exp\left(-\int_{x_0}^x \frac{b(y)}{a(y)} \, dy\right) \left[\exp\left(\int_{x_0}^x \frac{b(y)}{a(y)} \, dy\right) q_x\right]_x + c(x)q,
\eeq
which gives
\beq
    \alpha(x) = a(x)\exp\left(-\int_{x_0}^x \frac{b(y)}{a(y)} \, dy\right)\!, \qquad  \beta(x) = \exp\left(\int_{x_0}^x \frac{b(y)}{a(y)} \, dy\right)\!, \qquad \text{ and } \qquad \gamma(x) = c(x).
\eeq
From this, we have
\beq
    \fracbetanuargs{x} = \frac{1}{2}\left( \frac{\beta'(x)}{\beta(x)} -\frac{\alpha'(x)}{\alpha(x)}  + \frac{\gamma'(x)}{k^2+\gamma(x)}\right) = \frac{1}{2}\left( \frac{2b(x)}{a(x)} - \frac{a'(x)}{a(x)} + \frac{c'(x)}{k^2+c(x)}\right) \!,
\eeq
which we can see is not integrable (over an infinite or semi-infinite domain) if $a,b,c$ are constants with $ab\neq 0$. This presents a problem for our solution. However, we can make the change of variables,
\beq \label{eqn:qtransformation}
    q(x,t) = \exp\left(-\int_{x_0}^x \frac{b(y)}{2a(y)} \, dy\right)u(x,t).
\eeq
The PDE becomes
\beq
    u_t = a(x) u_{xx} + \left(\frac{a'(x)b(x) - a(x)b'(x)}{2a(x)} - \frac{b(x)^2}{4a(x)} + c(x)\right)u,
\eeq
for which, we have
\beq
    \alpha(x) = a(x), \qquad \beta(x) = 1, \qquad \gamma(x) =\frac{a'(x)b(x) - a(x)b'(x)}{2a(x)} - \frac{b(x)^2}{4a(x)} + c(x),
\eeq
and
\beq
    \fracbetanuargs{x} = \frac{1}{2}\left(-\frac{a'(x)}{a(x)} + \frac{\gamma'(x)}{k^2+\gamma(x)}\right)\!.
\eeq
In the case of constant coefficients, the integrability condition, Assumption~\enumref{ass:alphabetagamma}{enum:L1}, is satisfied and our solution is well defined. 
\subsubsection{Example: The constant-coefficient, advected heat equation}
Consider the constant-coefficient IBVP
\begin{subequations} \label{eqn:IBVP_illposed}
    \begin{align} \label{eqn:PDE_illposed}
        q_t &= q_{xx} + c q_x, && x\in \R, \quad t>0, \\ 
        q(x,0)&=q_0(x), && x\in \R, \\ 
        \lim_{|x|\to\infty} q(x,t) &= 0, && t>0.
    \end{align}
\end{subequations}
\no This problem is well posed for $c\in \R$ \cite{JC_fokas_book}. The PDE \rf{eqn:PDE_illposed} can be written in the form \rf{eqn:PDE_wl} as
\beq
    q_t = e^{-cx}\left( e^{cx} q_{x} \right)_x\!,
\eeq
\no with $\alpha(x) = e^{-cx}$, $\beta(x) = e^{cx}$, and $\gamma(x) = 0$. Since $\beta'/\beta - \alpha'/\alpha = 2c$ is not absolutely integrable over the real line, and Assumption~\enumref{ass:alphabetagamma}{enum:L1} is not satisfied. With the change of variables $q(x,t) = e^{-cx/2} u(x,t)$, the IBVP \rf{eqn:IBVP_illposed} becomes
\begin{subequations}
    \begin{align}
        u_t &= u_{xx} -\frac{c^2}{4} u, && x\in \R, \quad t>0, \\ 
        u(x,0)&=e^{cx/2} q_0(x), && x\in \R, \\ 
        \lim_{|x|\to\infty} u(x,t) &= 0, && t>0.
    \end{align}
\end{subequations}
\no Now $\alpha(x) = 1$, $\beta(x) = 1$, and $\gamma(x) = -c^2/4$, so that $\beta'/\beta - \alpha'/\alpha = 0$ and $\gamma' = 0$, and Assumption~\ref{ass:alphabetagamma} is satisfied.
This example shows that, although all evolution equations can be written in the form \rf{eqn:PDE_wl}, a transformation may be needed before the integrability conditions are met and the solution expression \rf{eqn:q_wl} applies.
\section{The half-line problem} \label{sec:statement_hl}
Consider \eqref{eqn:IVP} on the half line $x>x_l$ with a linear, constant-coefficient boundary condition and decay at infinity, 
\begin{subequations} \label{eqn:IBVP_hl}
    \begin{align}
        \label{eqn:PDE_hl}
        q_t &= \alpha(x)\left(\beta(x) q_{x}\right)_x+\gamma(x)q+f(x,t), && x>x_l, \quad t>0, \\ 
        \label{eqn:IC_hl}
        q(x,0)&=q_0(x), && x>x_l, \\ 
        \label{eqn:BC1_hl}
        f_0(t) &= a_{0}q(x_l,t)+ a_{1} q_x(x_l,t),  && t>0, \\
        \label{eqn:BC2_hl}
        \lim_{x\to\infty} q(x,t) &= 0, && t>0,
    \end{align}
\end{subequations}
with $(a_0,a_1)\neq (0,0)$.
\begin{theorem}
\no Under Assumptions~\ref{ass:alphabetagamma}~and~\ref{ass:ffunctions}, the IBVP \rf{eqn:IBVP_hl} has the solution
\beq \label{eqn:q_hl}
    q(x,t) = \frac{1}{2\pi} \int_{\partial \Omega} \frac{\Phi(k,x,t)}{\Delta(k)} e^{-k^2t} \, dk,
\eeq
\no where $\Omega$ is shown in Figure~\ref{fig:Omega}. Here
\beq \label{eqn:Delta_hl}
    \Delta(k) = 2 \sum_{n = 0}^\infty \left( \frac{(-1)^nia_0}{k\mathfrak n(k,x_l)} - a_1 \right) {\mathcal E}_{n}^{(x_l,\infty)}(k),
\eeq
and
\beq \label{eqn:Phi_hl}
    \Phi(k,x,t) = \mathcal B_0(k,x)F_0(k^2,t) + \Phi_\psi(k,x,t), 
\eeq
where
\beq \label{eqn:Phi_psi_hl}
    \Phi_\psi(k,x,t) = \int_{x_l}^{\infty} \frac{\Psi(k,x,y)\psi_\alpha(k^2,y,t)}{ \sqrt{(\beta \mathfrak n)(k,x)} \sqrt{(\beta \mathfrak n)(k,y)}} \, dy.
\eeq
The functions $\psi_\alpha(k^2,x,t)$ and $\mathfrak n(k,x)$ are defined in Definition~\ref{def:mu}. The boundary term $\mathcal B_0(k,x)$ is defined by
\beq \label{eqn:B0_hl}
    \mathcal B_0(k,x) =  \frac{4\beta(x_l) \etwo{ ik \int_{x_l}^{x} \mathfrak n(k,\xi) \, d\xi} }{\sqrt{(\beta \mathfrak n)(k,x_l)}\sqrt{(\beta \mathfrak n)(k,x)}}    \sum_{n=0}^\infty (-1)^n {\mathcal E}_{n}^{(x,\infty)}(k),
\eeq
and
\beq
\label{eqn:F_m}
F_m(k^2,t) = \int_0^t e^{k^2s} f_m(s) \, ds, \qquad m=0,1.
\eeq
Note that we will use $F_1(k^2,t)$ in the finite-interval problem in Section \ref{sec:statement_fi}.
\no For $x_l<y<x$,
\beq \label{eqn:Psi_hl}
    \Psi(k,x,y) = 4\exp\left(ik \int_{x_l}^{x} \mathfrak n(k,\xi) \, d\xi\right)\sum_{n=0}^\infty \sum_{\ell=0}^n (-1)^\ell \left(\frac{a_0}{k \mathfrak n(k,x_l)}\mathcal S_{n-\ell}^{(x_l,y)}(k) - a_1\mathcal C_{n-\ell}^{(x_l,y)}(k) \right) {\mathcal E}_\ell^{(x,\infty)}(k),
\eeq
\no and $\Psi(k,x,y) =\Psi(k,y,x)$ for $x_l<x<y$. 
${\mathcal E}_n^{(a,b)}(k)$ is defined in \rf{eqn:En}, $\mathcal C_0^{(a,b)}(k) = 1$, $\mathcal S_0^{(a,b)}(k) = 0$, and for $n\geq 1$,
\begin{subequations} \label{eqn:CnSn}
\begin{align}
\label{eqn:Cn}
\mathcal C_n^{(a,b)}(k) &= \frac{1}{2^n}\int_{a=y_0<y_1<\cdots < y_n<y_{n+1}=b}  \left(\prod_{p=1}^n \fracbetanuargs{y_p} \right) \cos\left( k \sum_{p=0}^{n} (-1)^p \int_{y_p}^{y_{p+1}}\mathfrak n(k,\xi)\,d\xi\right)  d\mathbf{y}_n, \\
\label{eqn:Sn}
\mathcal S_n^{(a,b)}(k) &= \frac{1}{2^n}\int_{a=y_0<y_1<\cdots < y_n<y_{n+1}=b}  \left(\prod_{p=1}^n \fracbetanuargs{y_p} \right) \sin\left( k \sum_{p=0}^{n} (-1)^p \int_{y_p}^{y_{p+1}}\mathfrak n(k,\xi)\,d\xi\right) d\mathbf{y}_n,
\end{align}
\end{subequations}
\no where $d\mathbf{y}_n = dy_1 \cdots dy_n$ and the prime denotes the derivative with respect to the second variable, as before.
\end{theorem}
\subsection{Example: The advected heat equation} \label{sec:example_hl}
Consider the advected heat equation on the half line with spatially variable thermal conductivity $\sigma^2(x)>0$ and velocity $c(x)$, without forcing and with homogeneous Dirichlet boundary conditions, i.e.,
\begin{subequations}
    \begin{align}
        q_t &= \left(\sigma^2(x) q_{x}\right)_x - c(x) q_x, && x>0, \quad t>0, \\ 
        q(x,0)&=q_0(x), && x>0, \\ 
        q(0,t)&= 0, && t>0,\\
        \lim_{x\to\infty} q(x,t) &= 0, && t>0.
    \end{align}
\end{subequations}
\no Here $x_l=0$, $a_0=1$ and $a_1=0$. Further,
\beq
    \alpha(x)=\exp\left( \int_0^x \frac{c(\xi)}{\sigma^2(\xi)} \, d\xi \right)
    \!, \qquad
    \beta(x) = \sigma^2(x) \exp\left( - \int_0^x \frac{c(\xi)}{\sigma^2(\xi)} \, d\xi \right)
    \!, \qquad \gamma(x) = 0,
\eeq
$f(x,t)=0$, and $f_0(t)=0$. 
We require absolute continuity of $\sigma(x)$, boundedness of $c(x)$, and since
\begin{align}
    \frac{\beta'(x)}{\beta'(x)} - \frac{\alpha'(x)}{\alpha(x)} 
    &= \frac{2\sigma'(x)}{\sigma(x)} - \frac{2c(x)}{\sigma^2(x)},
\end{align}
we require absolute integrability of $\sigma'(x)/\sigma(x)$ and $c(x)$, so that Assumption~\ref{ass:alphabetagamma} is satisfied. Note that if $\sigma(x)$ is absolutely continuous and $\sigma'(x)/\sigma(x)$ is absolutely integrable, then $\sigma(x)$ is bounded above and below.
This problem has the solution \eqref{eqn:q_hl},
\no where $\Omega$ is shown in Figure~\ref{fig:Omega}, $\mathfrak n(k,x) = 1/\sigma(x)$, 
\beq
    k\mathfrak n(k,0)\Delta(k) = 2i \sum_{n=0}^\infty (-1)^n {\mathcal E}_{n}^{(0,\infty)}(k).
\eeq
\no Since $\mathcal B_0(k,x,t) = 0$ and $\psi(k^2,y,t) = q_0(y)$, we have
\beq
    \Phi(k,x,t) = \int_{0}^{\infty} \exp\left( \frac12 \int_y^x \frac{c(\xi)}{\sigma^2(\xi)} \, d\xi\right)\frac{\Psi(k,x,y) q_0(y)}{ \sqrt{\sigma(x) \sigma(y)}} \, dy,
\eeq
\no and for $0<y<x$,
\beq
    k\mathfrak n(k,0)\Psi(k,x,y) = 4 \exp\left(ik\int_{0}^{x}\frac{d\xi}{\sigma(\xi)} \right)\sum_{n=0}^\infty \sum_{\ell=0}^n (-1)^\ell \mathcal S_{n-\ell}^{(0,y)}(k) {\mathcal E}_\ell^{(x,\infty)}(k) ,
\eeq
and $\Psi(k,x,y) =\Psi(k,y,x)$ for $0<x<y$. 
\section{The finite-interval problem} \label{sec:statement_fi}
Consider \eqref{eqn:IVP} on the finite interval $x_l<x<x_r$ with linear, constant-coefficient boundary conditions, 
\begin{subequations}
    \label{eqn:IBVP_fi}
    \begin{align} 
        \label{eqn:PDE_fi}
        q_t &= \alpha(x)\left(\beta(x) q_{x}\right)_x+\gamma(x)q+f(x,t), && x\in (x_l,x_r), \quad t>0, \\ 
        \label{eqn:IC_fi}
        q(x,0)&=q_0(x), && x\in(x_l,x_r), \\ 
        \label{eqn:BC1_fi}
        f_0(t) &= a_{11}q(x_l,t)+ a_{12} q_x(x_l,t)+ b_{11}q(x_r,t)+ b_{12} q_x(x_r,t) , && t>0, \\
        \label{eqn:BC2_fi}
        f_1(t) &= a_{21}q(x_l,t)+ a_{22} q_x(x_l,t)+ b_{21}q(x_r,t)+ b_{22} q_x(x_r,t) , && t>0. 
    \end{align}
\end{subequations}
Considering the concatenated matrix 
\begin{align} \label{eqn:ab}
    (a:b) &= \begin{pmatrix} a_{11} & a_{12} & b_{11} & b_{12} \\ a_{21} & a_{22} & b_{21} & b_{22} \end{pmatrix}\!,
\end{align}
we let $(a:b)_{i,j} = \det((a:b)_{\{1,2\},\{i,j\}})$ denote the determinant of the $2\times 2$ minor with columns at $i$ and $j$ \cite{LA_and_geom}. We require rank$(a:b)=2$ and one of the following {\em Boundary Cases}.
\begin{definition} \label{def:boundarycases}
    For $x\in \D = (x_l,x_r)$, we define the constants
    \begin{align}
        m_{\mathfrak c_0} =  \frac{(a:b)_{1,4}}{\mu(x_l)} - \frac{(a:b)_{2,3}}{\mu(x_r)}, \qquad
        m_{\mathfrak c_1} =\frac{(a:b)_{1,4}}{\mu(x_l)} + \frac{(a:b)_{2,3}}{\mu(x_r)} , \qquad
        m_{\mathfrak s} = \frac{(a:b)_{1,3}}{\mu(x_l)\mu(x_r)},
    \end{align}
    and $\mathfrak u_{\pm} =  \mathfrak u(x_r) \pm  \mathfrak u(x_l)$, where $\mu(x)$ and $\mathfrak u(x)$ are defined in Definition~\ref{def:mu}.
    We define the following {\em Boundary Cases:}
    \begin{enumerate}
        \item \label{enum:BC1} $(a:b)_{2,4}\neq 0$,
        \item \label{enum:BC2} $(a:b)_{2,4} = 0$ and $m_{\mathfrak c_0} \neq 0$,
        \item
        \label{enum:BC3} $(a:b)_{2,4} = 0$, $m_{\mathfrak c_0} = 0$, $m_{\mathfrak c_1} = 0$, and $(a:b)_{1,3} \neq 0$,
        \item
        \label{enum:BC4} $(a:b)_{2,4} = 0$, $m_{\mathfrak c_0} = 0$, $m_{\mathfrak c_1} \neq 0$, and $m_{\mathfrak c_1} \mathfrak u_{+} - 8m_{\mathfrak s} \neq 0$.
    \end{enumerate}
\end{definition}
Please refer to Remark~\ref{rem:BC} for an interpretation of these different Boundary Cases. 
\begin{theorem}
\no Under Assumptions~\ref{ass:alphabetagamma}~and~\ref{ass:ffunctions} (and for Boundary Case~\ref{enum:BC4}, Assumption~\ref{ass:alphabetagamma2}), the IBVP \rf{eqn:IBVP_fi} has the solution
\beq
    \label{eqn:q_fi}
    q(x,t) = \frac{1}{2\pi} \int_{\partial \Omega} \frac{\Phi(k,x,t)}{\Delta(k)} e^{-k^2t} \, dk,
\eeq
\no where $\Omega$ is shown in Figure~\ref{fig:Omega}. 
We define 
\begin{equation} \label{eqn:Xi_def}
    \Eint = \exp\left(ik \int_{x_l}^{x_r} \mathfrak n(k,\xi) \, d\xi\right)\!,
\end{equation}
where $\mathfrak n(k,x)$ is defined in Definition~\ref{def:mu}.
Then
\beq \label{eqn:Delta_fi}
    \Delta(k) = 2i\,\Xi(k) \left( \mathfrak a(k) + \sum_{n=0}^\infty \mathfrak c_n(k) \mathcal C_n^{(x_l,x_r)}(k)  + \sum_{n=0}^\infty \mathfrak s_n(k) \mathcal S_n^{(x_l,x_r)}(k) \right)\!,
\eeq
\no with 
\begin{subequations} \label{eqn:Delta_a,cn,sn}
    \begin{align}
        \mathfrak a(k) &= \frac{\beta(x_r)(a:b)_{1,2}+\beta(x_l)(a:b)_{3,4}}{k \sqrt{\betanuargs{x_l} }\sqrt{\betanuargs{x_r}}}, \\
        \mathfrak c_n(k) &= (-1)^{n}\frac{(a:b)_{1,4}}{k \mathfrak n(k,x_l)} - \frac{(a:b)_{2,3}}{k \mathfrak n(k,x_r)}, \\
        \mathfrak s_n(k) &= (-1)^n (a:b)_{2,4} + \frac{(a:b)_{1,3}}{k^2 \mathfrak n(k,x_l)\mathfrak n(k,x_r)} .
    \end{align}
\end{subequations}
\no The numerator of \rf{eqn:q_fi} is 
\begin{subequations}
    \beq \label{eqn:Phi_fi}
        \Phi(k,x,t) = \mathcal B_0(k,x) F_0(k^2,t) + \mathcal B_1(k,x)F_1(k^2,t) + \Phi_\psi(k,x,t),
    \eeq
    where
    \beq \label{eqn:Phi_psi_fi_def}
        \Phi_\psi(k,x,t) = \int_{x_l}^{x_r} \frac{\Psi(k,x,y)\psi_\alpha (k^2,y,t)}{ \sqrt{\betanuargs{x}} \sqrt{\betanuargs{y}}} \, dy.
    \eeq
    \no The function $\psi_\alpha(k^2,x,t)$ is defined in Definition~\ref{def:mu}, $F_m(k^2,t)$ is defined in \eqref{eqn:F_m}, and the boundary terms $\mathcal B_{0}(k,x)$ and $\mathcal B_1(k,x)$ are given by
    \begin{align} \label{eqn:Bm} \nonumber
        \mathcal B_{2-j}(k,x) &= (-1)^{j} \frac{4\Eint}{\sqrt{\betanuargs{x}}} \left\{ \frac{\beta(x_r)}{\sqrt{\betanuargs{x_r}}} \left[-\frac{a_{j1} }{k \mathfrak n(k,x_l)}  \sum_{n=0}^\infty \mathcal S_n^{(x_l,x)}(k)  + a_{j2} \sum_{n=0}^\infty \mathcal C_n^{(x_l,x)}(k) \right] \right. \\
        &  \hspace*{1.08in} \left. +\frac{\beta(x_l)}{\sqrt{\betanuargs{x_l}}} \left[ \frac{b_{j1}}{k \mathfrak n(k,x_r)}  \sum_{n=0}^\infty \mathcal S_n^{(x,x_r)}(k) + b_{j2} \sum_{n=0}^\infty (-1)^n \mathcal C_n^{(x,x_r)}(k)\right] \right\}\!, ~~ j=1,2. 
    \end{align}
\end{subequations}
\no Further, for $x_l<y<x<x_r$,
\begin{subequations} \label{eqn:Psi_fi}
    \begin{align} \nonumber \label{eqn:Psi_fi1}
        \Psi(k,x,y) &= 4 \Eint \left\{ - (a:b)_{2,4} \sum_{n=0}^\infty \sum_{\ell=0}^n (-1)^\ell \mathcal C_{n-\ell}^{(x_l,y)}(k)\mathcal C_\ell^{(x,x_r)}(k) + \frac{(a:b)_{1,3}}{k^2 \mathfrak n(k,x_l) \mathfrak n (k,x_r)}\sum_{n=0}^\infty \sum_{\ell=0}^n \mathcal S_{n-\ell}^{(x_l,y)}(k)\mathcal S_\ell^{(x,x_r)}(k) \right. \\ \nonumber
        &\hspace*{0.7in}+ \frac{(a:b)_{1,4}}{k \mathfrak n(k,x_l)} \sum_{n=0}^\infty \sum_{\ell=0}^n (-1)^\ell \mathcal S_{n-\ell}^{(x_l,y)}(k) \mathcal C_\ell^{(x,x_r)}(k) - \frac{(a:b)_{2,3}}{k \mathfrak n(k,x_r)} \sum_{n=0}^\infty \sum_{\ell=0}^n \mathcal C_{n-\ell}^{(x_l,y)}(k) \mathcal S_\ell^{(x,x_r)}(k) \\ 
        &\hspace*{0.7in}-\left. \frac{\beta(x_r)(a:b)_{1,2}}{k \sqrt{\betanuargs{x_l}} \sqrt{\betanuargs{x_r}}} \sum_{n=0}^\infty \mathcal S_n^{(y,x)}(k) \right\}\!,
    \end{align}
    \no and, for $x_l<x<y<x_r$,
    \begin{align} \nonumber \label{eqn:Psi_fi2}
        \Psi(k,x,y) &= 4\Eint \left\{-(a:b)_{2,4} \sum_{n=0}^\infty \sum_{\ell=0}^n (-1)^\ell \mathcal C_{n-\ell}^{(x_l,x)}(k)\mathcal C_\ell^{(y,x_r)}(k)  +\frac{(a:b)_{1,3}}{k^2 \mathfrak n(k,x_l) \mathfrak n(k,x_r)} \sum_{n=0}^\infty \sum_{\ell=0}^n \mathcal S_{n-\ell}^{(x_l,x)}(k) \mathcal S_\ell^{(y,x_r)}(k) \right. \\ \nonumber 
        &\hspace*{0.7in}+ \frac{(a:b)_{1,4}}{k \mathfrak n(k,x_l)} \sum_{n=0}^\infty \sum_{\ell=0}^n (-1)^\ell  \mathcal S_{n-\ell}^{(x_l,x)}(k) \mathcal C_\ell^{(y,x_r)}(k) - \frac{(a:b)_{2,3}}{k \mathfrak n(k,x_r)} \sum_{n=0}^\infty \sum_{\ell=0}^n  \mathcal C_{n-\ell}^{(x_l,x)}(k)\mathcal S_\ell^{(y,x_r)}(k) \\ 
        &\hspace*{0.7in}\left. -\frac{\beta(x_l) (a:b)_{3,4}}{k \sqrt{\betanuargs{x_l}} \sqrt{ \betanuargs{x_r}}} \sum_{n=0}^\infty \mathcal S_n^{(x,y)}(k) \right\}\!.
    \end{align}
\end{subequations}
\no Note that $\Psi(k,x,y) \neq \Psi(k,y,x)$ unless $\beta(x_r)(a:b)_{1,2}=\beta(x_l)(a:b)_{3,4}$. The functions $\mathcal C_n^{(a,b)}(k)$ and $\mathcal S_n^{(a,b)}(k)$ are defined in \rf{eqn:Cn} and \rf{eqn:Sn}, respectively. 
\end{theorem}
\newcommand{\nonzero}[1]{\blue{\underline{#1}}}
\begin{rem} \label{rem:BC} 
    We use {\nonzero{underline}} to denote non-zero terms in this remark. Further, we use row reduction and the fact that the order of equations \eqref{eqn:BC1_fi} and \eqref{eqn:BC2_fi} is irrelevant.
    \begin{enumerate}
        \item If $(a:b)_{2,4} \neq 0$, the most general form of the matrix $(a:b)$ in \eqref{eqn:ab} is
        $$ (a:b) = \begin{pmatrix} 
            a_{11} & \nonzero{a_{12}} & b_{11} & 0 \\ 
            a_{21} & 0 & b_{21} & \nonzero{b_{22}} 
        \end{pmatrix}\!. $$
        This case includes the classical Neumann and Robin boundary conditions at both boundaries. We refer to these as {\em Robin-type boundary conditions}. In the case of constant coefficients, this is Birkhoff regular \cite{locker}.
        \item \label{enum:BC2_reduction} If $(a:b)_{2,4} = 0$ and $m_{\mathfrak c_0} \neq 0$, the most general form of the matrix $(a:b)$ in \eqref{eqn:ab} are
        $$(a:b) =  \begin{pmatrix} 
            a_{11} & {\nonzero{a_{12}}} & 0 & 0 \\ 
            a_{21} & 0 & {\nonzero{b_{21}}} & 0
        \end{pmatrix}\!, ~ \begin{pmatrix} 
            0 & 0 & b_{11} & \nonzero{b_{12}} \\
            \nonzero{a_{21}} & 0 & b_{21} & 0
        \end{pmatrix}\!, ~ \begin{pmatrix} 
            \nonzero{a_{11}} & 0 & 0 & 0 \\ 
            0 & a_{22} & b_{21} & \nonzero{b_{22}}
        \end{pmatrix}\!, ~ \begin{pmatrix} 
            0 & 0 & \nonzero{b_{11}} & 0 \\
            a_{21} & \nonzero{a_{22}} & 0 & b_{22}            
        \end{pmatrix}\!; $$
        or 
        $$(a:b) = \begin{pmatrix} 
            \nonzero{a_{11}} & 0 & \nonzero{b_{11}} & 0  \\
            0 & \nonzero{a_{22}} & b_{21} & \nonzero{b_{22}}
        \end{pmatrix} \!,
        ~~~~~~ \text{where} ~~~~~~
        \frac{{\nonzero{a_{11}}} \, {\nonzero{b_{22}}}}{\mu(x_l)} + \frac{{\nonzero{a_{22}}} \, {\nonzero{b_{11}}}}{\mu(x_r)}  \neq 0.
        $$
        This case includes a Robin boundary condition on the left (or right) and a Dirichlet boundary condition on the right (or left). It also includes the classical periodic `boundary conditions'. We refer to these as {\em mixed-type or periodic-type boundary conditions}. In the case of constant coefficients, this is Birkhoff regular \cite{locker}.
        \item \label{enum:rem:BC3} If $(a:b)_{2,4} = 0$, $m_{\mathfrak c_0} = 0$, $m_{\mathfrak c_1} = 0$, and $(a:b)_{1,3} \neq 0$, the most general form of the matrix $(a:b)$ in \eqref{eqn:ab} are
        $$ (a:b) = \begin{pmatrix} 
            \nonzero{a_{11}} & 0 & 0 & b_{12} \\ 
            0 & 0 & \nonzero{b_{21}} & 0 
        \end{pmatrix} \qquad \text{ or } \qquad (a:b) = \begin{pmatrix} 
            \nonzero{a_{11}} & 0 & 0 & 0 \\
            0 & a_{22} & \nonzero{b_{21}} & 0             
        \end{pmatrix}\!.$$
        This case includes the case of the classical Dirichlet boundary conditions. We refer to these as {\em Dirichlet-type boundary conditions}. In the case of constant coefficients, this is Birkhoff regular for the case of Dirichlet boundary conditions $($\ie if $a_{22} = 0 = b_{12}$ or, equivalently, if $(a:b)_{1,2} = 0 = (a:b)_{3,4})$ and is Birkhoff irregular if $a_{22} \neq 0 $ or $b_{12} \neq 0$ {$($}or, equivalently, if $(a:b)_{1,2} \neq 0$ or $(a:b)_{3,4}\neq 0)$ \cite{locker}.
        \item If $(a:b)_{2,4} = 0$, $m_{\mathfrak c_0} = 0$, $m_{\mathfrak c_1} \neq 0$, and $m_{\mathfrak c_1} \mathfrak u_+ - 8m_{\mathfrak s} \neq 0$, the most general form of the matrix $(a:b)$ is
        $$(a:b) = \begin{pmatrix} 
            \nonzero{a_{11}} & 0 & \nonzero{b_{11}} & 0 \\
            a_{21} & \nonzero{a_{22}} & 0 & \nonzero{b_{22}}
        \end{pmatrix}\!,
        ~~~~ \text{ where } ~~~~ 
        \frac{{\nonzero{a_{11}}} \, {\nonzero{b_{22}}}}{\mu(x_l)}   + \frac{ {\nonzero{a_{22}}} \, {\nonzero{b_{11}}} }{\mu(x_r)} = 0
        ~~~~ \text{ and } ~~~~
         a_{21}  \neq -\frac{1}{4} \mu(x_l) {\nonzero{a_{22}}} \mathfrak u_+.
        $$
        This case does not include any classical boundary conditions. Instead, it is an interface problem on a circle.  In the case of constant coefficients, this is Birkhoff irregular \cite{locker}.
    \end{enumerate}
\end{rem}
\subsection{Example: The heat equation with homogeneous, Dirichlet boundary conditions} \label{sec:example_fi2}
Consider the heat equation on the finite interval with spatially varying thermal conductivity $\sigma^2(x)$ without forcing and with homogeneous Dirichlet boundary conditions, \ie
\begin{subequations}
    \begin{align}
        q_t &= \big(\sigma^2(x) q_{x}\big)_x, && x\in(0,1), \quad t>0, \\ 
        q(x,0)&=q_0(x), && x\in(0,1), \\ 
        q(0,t)&= 0, && t>0,\\
        q(1,t) &= 0, && t>0.
    \end{align}
\end{subequations}
\no We let $x_l=0$, $x_r=1$, $\alpha(x)=1$, $\beta(x) = \sigma^2(x)$, $\gamma(x) = 0$, $f(x,t)=0$, $f_m(t)=0$ ($m=0,1$), and
\beq
    (a:b) = \begin{pmatrix} 1 & 0 & 0 & 0 \\ 0 & 0 & 1 & 0 \end{pmatrix}\!.
\eeq
Since $(a:b)_{2,4} = 0$, $m_{\mathfrak c_0} = 0$, $m_{\mathfrak c_1} = 0$, and $(a:b)_{1,3} = 1 \neq 0$, this is an example of Boundary Case~\ref{enum:BC3} and it is {\em regular}.
\no We require absolute continuity of $\sigma(x)$, integrability of $q_0(x)$, and absolutely integrability of $\sigma'(x)/\sigma(x)$. This has the solution
\beq
    q(x,t) = \frac{1}{2\pi} \int_{\partial \Omega} \frac{\Phi(k,x)}{\Delta (k)}  e^{-k^2t}\, dk,
\eeq
\no where $\Omega$ is shown in Figure~\ref{fig:Omega}. Since $\mathfrak n(k,x) = 1/\sigma(x)$, $\mathfrak a(k) = 0$, $\mathfrak c_n(k) = 0$,  and $\mathfrak s_n(k) = \sigma(0)\sigma(1)/k^2$, then
\beq
    \frac{k^2\Delta(k)}{2i\sigma(0)\sigma(1)} = \etwo{ik \int_0^1 \frac{d\xi}{\sigma(\xi)}} \sum_{n=0}^\infty  \mathcal S_n^{(0,1)}(k),
\eeq
\no and since $\mathcal B_0(k,x) = \mathcal B_1(k,x) = 0$ and $\psi_\alpha(k^2,y,t) = q_0(y)$, we have
\beq
    \Phi(k,x) = \int_0^1 \frac{\Psi(k,x,y)q_0(y)}{\sqrt{\sigma(x)\sigma(y)}} \, dy,
\eeq
\no where, for $0<y<x<1$,
\beq
\frac{k^2\Psi(k,x,y)}{4\sigma(0)\sigma(1)} =  \sum_{n=0}^\infty \sum_{\ell=0}^n \mathcal S_{n-\ell}^{(0,y)}(k) \mathcal S_\ell^{(x,1)}(k),
\eeq
\no and $\Psi(k,x,y) = \Psi(k,y,x)$ for $0<x<y< 1$. This is the same solution given in \cite{Farkas_Deconinck}. It reduces to the solution given in \cite{JC_fokas_book} for constant $\sigma(x)$.
\subsection{Example: The CGL equation with periodic boundary conditions} \label{sec:example_fi1}
\label{sec:Ginzburg-Landau}
The complex Ginzburg-Landau (CGL) equation is the nonlinear PDE
\begin{align} \label{eqn:NLcGLE}
    A_t = (1+ia(x)) A_{xx} + A - (1 + i b(x))|A|^2A,
\end{align}
where $a,b$ are real functions of $x$. In the special case $a(x) = 0 = b(x)$, \eqref{eqn:NLcGLE} is the real Ginzburg-Landau equation. If $a(x),b(x)\to\infty$, \eqref{eqn:NLcGLE} becomes the Nonlinear Schr\"odinger (NLS) equation \cite{ginzburg-landau}.
Consider the linearized (about $A=0$), CGL equation with periodic boundary conditions:
\begin{subequations}
    \begin{align}
        A_t &= (1+ia(x)) A_{xx} + A, && x\in(0,1), \quad t>0, \\ 
        A(x,0)&=A_0(x), && x\in(0,1), \\ 
        A(0,t)&= A(1,t), && t>0,\\
        A_x(0,t) &= A_x(1,t), && t>0.
    \end{align}
\end{subequations}
\no Here $x_l=0$, $x_r=1$, $\alpha(x) = 1+ia(x)$, $\beta(x) = 1$, $\gamma(x) = 1$, $f(x,t) = 0$, $f_0(t) = 0=f_1(t)$, and
\beq
    (a:b)=\begin{pmatrix}
        1 & 0 & -1 & 0 \\ 0 & 1 & 0 & -1
    \end{pmatrix}\!.
\eeq
We assume $a(x)\in \R$ and $a \in \mathrm{AC}(\D)$, so that Assumption~\ref{ass:alphabetagamma} is satisfied. Here, 
\begin{align}
    \mu(x) = \frac{e^{-\frac{i}{2}\arctan(a(x))}}{\sqrt[4]{1+a(x)^2}}  
    \qquad \text{ and }\qquad
    \mathfrak g(k) = \sqrt{1 + \frac{1}{k^2}},
\end{align}
and $\mathfrak n(k,x) = \mu(x) \mathfrak g(k)$, where the square root in $\mathfrak g(k)$ is defined in \eqref{eqn:gfrak}. Since $(a:b)_{2,4}=0$ and $m_{\mathfrak c_0} \neq 0$, this is a Boundary Case~\ref{enum:BC2} example, which is {\em regular}. For simplicity, we assume that $a(x)$ is periodic, \ie $a(0)=a(1)$. This problem has the solution
\beq
    q(x,t) 
    = - \frac{1}{2\pi i} \int_{\partial \Omega} \frac{\tilde \Phi(k,x)}{ \tilde \Delta (k)} e^{-k^2t} \, dk,
\eeq
\no where we define $\tilde \Delta(k) = k\mathfrak n(k,0)\Xi(-k) \Delta(k)/(4i)$ and $\tilde \Phi(k,x) = -k \mathfrak n(k,0) \Xi(-k) \Phi(k,x)/4$, and where $\Omega$ is shown in Figure~\ref{fig:Omega}. Here, 
$\mathfrak a(k) = 2/(k \mathfrak n(k,0))$,
$\mathfrak c_n(k) = -(1 + (-1)^n)/(k\mathfrak n(k,0)) $, 
$\mathfrak s_n(k) = 0$, 
\no and since
$\mathcal B_0(k,x,t) = 0 = \mathcal B_1(k,x,t)$ and $\psi_\alpha(k^2,x,t) = A_0(x)/(1+ia(x))$,
\beq \label{eqn:Delta_GLE}
    \tilde \Delta (k) =   1 - \sum_{n=0}^\infty \mathcal C_{2n}^{(0,1)}(k) 
    \qquad \text{ and } \qquad 
    \tilde \Phi(k,x) = \int_0^1 \frac{\tilde \Psi(k,x,y) A_0(y)}{(1+ia(y))\sqrt{ \mathfrak n(k,x) } \sqrt{\mathfrak n(k,y)} }  \, dy.
\eeq
We define $\tilde \Psi(k,x,y) = -k \mathfrak n(k,0) \Xi(-k) \Psi(k,x)/4$,
\beq
    \tilde \Psi(k,x,y) = \sum_{n=0}^\infty \sum_{\ell=0}^n (-1)^\ell \mathcal S_{n-\ell}^{(0,y)}(k) \mathcal C_\ell^{(x,1)}(k) + \sum_{n=0}^\infty \sum_{\ell=0}^n \mathcal C_{n-\ell}^{(0,y)}(k) \mathcal S_\ell^{(x,1)}(k)  + \sum_{n=0}^\infty \mathcal S_n^{(y,x)}(k) ,
\eeq
\no for $0<y<x<1$, and $\tilde \Psi(k,x,y) = \tilde \Psi(k,y,x)$ for $0<x<y< 1$.
\subsection{Sturm-Liouville Problems: Eigenvalues and Eigenfunctions} \label{sec:eigenvalues}
\begin{theorem} \label{thm:efuns}
    The Sturm--Liouville problem 
    \beq \label{eqn:eigenproblem}
        \alpha(x) \left( \beta(x)  y'\right)' +\gamma(x) y = \lambda y,
    \eeq
    \no with boundary conditions
    \begin{subequations} \label{eqn:eigenproblem_BC}
        \begin{align}
            \label{eqn:eigenproblem_BC_1}
            a_{11} y(x_l)+a_{12} y'(x_l) + b_{11} y(x_r) + b_{12} y'(x_r) &= 0, \\
            \label{eqn:eigenproblem_BC_2}
            a_{21} y(x_l)+a_{22} y'(x_l) + b_{21}y(x_r) + b_{22} y'(x_r) &= 0,
        \end{align}
    \end{subequations}
    has the eigenfunctions 
    \beq \label{eqn:efuns}
        X_m(x) = \frac{C_m}{\sqrt{(\beta\mathfrak n)(\kappa_m,x)}} \sum_{n=0}^\infty \mathcal C_n^{(x_l,x)}(\kappa_m) + \frac{S_m}{\sqrt{(\beta\mathfrak n)(\kappa_m,x)}} \sum_{n=0}^\infty \mathcal S_n^{(x_l,x)}(\kappa_m), 
    \eeq
    \no corresponding to the eigenvalues $\lambda_m = -\kappa_m^2$, where $\{\kappa_m\}_{m=1}^\infty$ are the zeros of $\Delta(k)$ \rf{eqn:Delta_fi}. Here,
    \begin{subequations}
        \begin{align} 
            C_m &= - \frac{a_{12} \kappa_m \mathfrak n(\kappa_m,x_l)}{\sqrt{(\beta\mathfrak n)(\kappa_m,x_l)}} -  \frac{b_{11}} {\sqrt{(\beta\mathfrak n)(\kappa_m,x_r)}} \sum_{n=0}^\infty \mathcal S_n^{(x_l,x_r)}(\kappa_m) - \frac{b_{12} \kappa_m \mathfrak n(\kappa_m,x_r)}{\sqrt{(\beta\mathfrak n)(\kappa_m,x_r)}} \sum_{n=0}^\infty (-1)^n \mathcal C_n^{(x_l,x_r)}(\kappa_m), \\
            S_m &= \frac{a_{11}}{\sqrt{(\beta\mathfrak n)(\kappa_m,x_l)}} + \frac{b_{11}}{\sqrt{(\beta\mathfrak n) (\kappa_m,x_r)}} \sum_{n=0}^\infty \mathcal C_n^{(x_l,x_r)}(\kappa_m)  - \frac{b_{12} \kappa_m \mathfrak n(\kappa_m,x_r)}{\sqrt{(\beta\mathfrak n)(\kappa_m,x_r)}} \sum_{n=0}^\infty (-1)^n \mathcal S_n^{(x_l,x_r)}(\kappa_m). 
        \end{align}
    \end{subequations}
\end{theorem}
\begin{proof}
    Using \eqref{eqn:Enderivatives} in \eqref{eqn:efuns} gives that the eigenfunctions solve the eigenvalue equation \eqref{eqn:eigenproblem}. Inserting \eqref{eqn:efuns} into the boundary conditions \eqref{eqn:eigenproblem_BC_1}, we find
    \begin{align}
        a_{11} X_m(x_l) + a_{12} X_m'(x_l) + b_{11} X_m(x_r) + b_{12} X_m'(x_r) = C_m S_m - S_m C_m = 0.
    \end{align}
    For \eqref{eqn:eigenproblem_BC_2}, we find
    \begin{align} \label{eqn:eigenproblem_BC_2_intermediate} \nonumber
        &a_{21} X_m(x_l)+a_{22} X_m'(x_l) + b_{21} X_m(x_r) + b_{22} X_m'(x_r) \\ \nonumber
        &= C_m \left[ \frac{a_{21}}{\sqrt{(\beta\mathfrak n)(\kappa_m,x_l)}} + \frac{b_{21}}{\sqrt{(\beta\mathfrak n)(\kappa_m,x_r)}} \sum_{n=0}^\infty \mathcal C_n^{(x_l,x_r)}(\kappa_m) - \frac{b_{22} \kappa_m \mathfrak n(\kappa_m,x_r)}{\sqrt{(\beta \mathfrak n)(\kappa_m,x_r)}} \sum_{n=0}^\infty (-1)^n \mathcal S_n^{(x_l,x_r)}(\kappa_m) \right] \\
        &~~~+ S_m \left[ \frac{a_{22} \kappa_m\mathfrak n(\kappa_m,x_l)}{\sqrt{(\beta\mathfrak n)(\kappa_m,x_l)}} + \frac{b_{21}}{\sqrt{\beta\mathfrak n)(\kappa_m,x_r)}} \sum_{n=0}^\infty \mathcal S_n^{(x_l,x_r)}(\kappa_m) + \frac{b_{22}\kappa_m \mathfrak n(\kappa_m,x_r)}{\sqrt{(\beta\mathfrak n)(\kappa_m,x_r)}} \sum_{n=0}^\infty (-1)^n \mathcal C_n^{(x_l,x_r)}(\kappa_m) \right]\!.
    \end{align}
    Expanding \eqref{eqn:eigenproblem_BC_2_intermediate}, we obtain
    \begin{align} \label{eqn:eigenproblem_BC_2_intermediate_2} 
        \nonumber
        & a_{21} X_m(x_l)+a_{22} X_m'(x_l) + b_{21} X_m(x_r) + b_{22} X_m'(x_r) \\ \nonumber
        &= \frac{(a:b)_{1,2} \kappa_m}{\beta(x_l)} + \frac{(a:b)_{3,4} \kappa_m}{\beta(x_r)} \left[\sum_{n=0}^\infty (-1)^n \mathcal C_n^{(x_l,x_r)}(\kappa_m) \sum_{n=0}^\infty \mathcal C_n^{(x_l,x_r)}(\kappa_m) + \sum_{n=0}^\infty (-1)^n \mathcal S_n^{(x_l,x_r)}(\kappa_m)\sum_{n=0}^\infty \mathcal S_n^{(x_l,x_r)}(\kappa_m) \right] \\
        &~~~+\frac{\kappa_m^2\mathfrak n(\kappa_m,x_l) \mathfrak n(\kappa_m,x_r)}{\sqrt{(\beta\mathfrak n)(\kappa_m,x_l)} \sqrt{(\beta\mathfrak n)(\kappa_m,x_r)}} \left[ \sum_{n=0}^\infty \mathfrak c_n(k) \mathcal C_n^{(x_l,x_r)}(k) +  \sum_{n=0}^\infty \mathfrak s_n(k) \mathcal S_n^{(x_l,x_r)}(k) \right]\!.
    \end{align}
    Using the identity
    \beq \label{eqn:BC_identity}
        1 = \sum_{n=0}^\infty (-1)^n \mathcal C_n^{(x_l,x_r)}(k) \sum_{n=0}^\infty \mathcal C_n^{(x_l,x_r)}(k) + \sum_{n=0}^\infty (-1)^n \mathcal S_n^{(x_l,x_r)}(k)\sum_{n=0}^\infty \mathcal S_n^{(x_l,x_r)}(k),
    \eeq
    in \eqref{eqn:eigenproblem_BC_2_intermediate_2}, this becomes
    \begin{align}
        & a_{21} X_m(x_l)+a_{22} X_m'(x_l) + b_{21} X_m(x_r) + b_{22} X_m'(x_r) = \frac{\kappa_m^2\mathfrak n(\kappa_m,x_l) \mathfrak n(\kappa_m,x_r)}{\sqrt{(\beta\mathfrak n)(\kappa_m,x_l)} \sqrt{(\beta\mathfrak n)(\kappa_m,x_r)}} \Delta(\kappa_m) = 0,
    \end{align}
    and the second boundary condition \eqref{eqn:eigenproblem_BC_2} is satisfied.
    
    To prove \eqref{eqn:BC_identity}, we define the right-hand side as $\mathfrak e_1$ and rewrite it as a Cauchy product, obtaining
    \begin{align}
        \mathfrak e_1 &= \sum_{n=0}^\infty  \sum_{\ell=0}^n (-1)^\ell\left[  \mathcal C_\ell^{(x_l,x_r)}(k) \mathcal C_{n-\ell}^{(x_l,x_r)}(k) + \mathcal S_\ell^{(x_l,x_r)}(k) \mathcal S_{n-\ell}^{(x_l,x_r)}(k) \right]\!.
    \end{align}
    Letting $\ell \to n-\ell$ in the inner sum, we see that
    \begin{align}
         \sum_{\ell=0}^n (-1)^\ell\left[  \mathcal C_\ell^{(x_l,x_r)}(k) \mathcal C_{n-\ell}^{(x_l,x_r)}(k) + \mathcal S_\ell^{(x_l,x_r)}(k) \mathcal S_{n-\ell}^{(x_l,x_r)}(k) \right] = 0,
    \end{align}
    for odd $\ell$.
    The $n=0$ term is $1$. The $n \geq 2$ even terms are $0$ and thus gives $\mathfrak e_1 = 1$. For $n=2$, we show
    \begin{align}\nonumber
        0 &= \int_{x_l}^{x_r} dz_1 \int_{z_1}^{x_r} dz_2 \, \cos\left( \int_{x_l}^{x_r} - \int_{x_l}^{z_1} + \int_{z_1}^{z_2} - \int_{z_2}^{x_r} \nu(k,\xi) \, d\xi \right) \\\nonumber
        &- \int_{x_l}^{x_r} dy_1 \int_{x_l}^{x_r} dz_1 \, \cos\left( \int_{x_l}^{y_1} - \int_{y_1}^{x_r} - \int_{x_l}^{z_1} + \int_{z_1}^{x_r} \nu(k,\xi) \, d\xi \right) \\
        &+ \int_{x_l}^{x_r} dy_1 \int_{y_1}^{x_r} dy_2 \, \cos\left( \int_{x_l}^{y_1} - \int_{y_1}^{y_2} + \int_{y_2}^{x_r} - \int_{x_l}^{x_r} \nu(k,\xi) \, d\xi \right)\!.
    \end{align}
    \no Let $I_j$ denote the three integrals above, in order. Since the first and the last term are equal and equal to
    \begin{equation}
        I_1 = I_3 = \int_{x_l}^{x_r} dy_1 \int_{y_1}^{x_r} dy_2 \, \cos\left( 2\int_{y_1}^{y_2} \nu(k,\xi) \, d\xi \right)\!,
    \end{equation}
    \no and since the second term is
    \begin{align}\nonumber
    I_2 &= - \int_{x_l}^{x_r} dy_1 \int_{x_l}^{x_r} dz_1 \, \cos\left( 2 \int_{z_1}^{y_1} \nu(k,\xi) \, d\xi \right) \\\nonumber
    &= - \int_{x_l}^{x_r} dy_1 \int_{x_l}^{y_1} dz_1 \, \cos\left( 2 \int_{z_1}^{y_1} \nu(k,\xi) \, d\xi \right) - \int_{x_l}^{x_r} dy_1 \int_{y_1}^{x_r} dz_1 \, \cos\left( 2 \int_{z_1}^{y_1} \nu(k,\xi) \, d\xi \right) \\
    &= -2 \int_{x_l}^{x_r} dy_1 \int_{y_1}^{x_r} dz_1 \, \cos\left( 2 \int_{z_1}^{y_1} \nu(k,\xi) \, d\xi \right)\!,
    \end{align}
    \no and so the $n=2$ term is $0$. The other $n$ terms are similar.
\end{proof}
Comparing to P\"{o}schel and Trubowitz~\cite{PT} ($\alpha(x)=\beta(x)=1$, $\gamma(x) = -q(x)$, $\lambda = k^2$), we find that
\beq 
y_1(x, k^2, q) = \sqrt{\frac{\mathfrak n(k,0)}{\mathfrak n(k,x)}} \sum_{n=0}^\infty \mathcal C_n^{(0,x)}(k) \And y_2(x, k^2, q) = \frac{1}{k \sqrt{\mathfrak n(k,0)} \sqrt{\mathfrak n(k,x)}} \sum_{n=0}^\infty \mathcal S_n^{(0,x)}(k), 
\eeq
\no where $\mathfrak n(k,x) = \sqrt{1-q(x)/k^2}$.
\subsubsection{Example: Eigenvalues for the CGL equation with periodic boundary conditions}
We revisit the complex Ginzburg-Landau equation described in Section~\ref{sec:Ginzburg-Landau}, setting $a(x) = x\sin(2\pi x)$. The associated eigenvalue problem is of the form
\begin{align}
    (1+ia(x)) y'' + y &= \lambda y, \qquad y(0)= y(1), \qquad y'(0) = y'(1).
\end{align}
The eigenvalues $\lambda_m = -\kappa_m^2$ are related to the zeroes $\kappa_m$ ($m=0,1,2,\ldots$) of $\tilde \Delta(k)$ \eqref{eqn:Delta_GLE}. Since $\tilde \Delta(k)$ is even in $k$, if $\kappa_m$ is a root, so is $-\kappa_m$, and each gives rise to the same eigenvalue. 
Since $\mathfrak g(\pm i) = 0$, and
\begin{align}
    \mathcal C_n^{(0,1)}(0) &= \frac{1}{2^n} \int_{0<\cdots<1} \left( \prod_{p=1}^n \frac{\mu'(y_p)}{\mu(y_p)} \right) \, d\mathbf{y}_n = \frac{1}{2^n n!} \left( \int_0^1 \frac{\mu'(y)}{\mu(y)} \, dy \right)^n = \frac{1}{2^n n!} \left( \log\left( \frac{\mu(1)}{\mu(0)} \right) \right)^n = 0,
\end{align}
then $\tilde \Delta(\pm i) = 0$, and $\kappa_0 = i$ (and $-i$) is an exact double root of $\tilde \Delta(k)$, and $\lambda_0 = -\kappa_0^2 = 1$ is an exact eigenvalue of the problem, which can be confirmed directly (with the constant eigenfunction).
We define
\begin{align}
    \mathfrak m(x)  = \int_0^x \mu(\xi) \, d\xi 
    \qquad \text{ and } \qquad 
    \eta(y) = \fracbetanuargs{x} = \frac{\mu'(x)}{\mu(x)} = \frac{2\pi x\cos(2\pi x) + \sin(2\pi x)}{2i - 2x\sin(2\pi x)}.
\end{align}
We truncate \eqref{eqn:Delta_GLE} at order $n=N$ and denote as $\tilde \Delta_N(k)$. Denoting $\kell = k\mathfrak g(k)$, the zeroth-order approximations of the roots of $\tilde \Delta(k)$ are
\begin{align}
    \tilde \Delta_0(k) = 1 - \cos\left(\mathfrak m(1) \kell \right) = 0 \then \kappa_m^{(0)} = \pm \frac{\sqrt{4m^2\pi^2 -\mathfrak m(1)^2}}{\mathfrak m(1)}, \qquad m=1,2,3,\ldots.
\end{align}
As in the case $\kappa_0 = \pm i$, these approximations are double roots. However, the actual eigenvalues are simple roots that are near these points. The next-order approximations $\kappa_m^{(1)}$ are the roots of 
\begin{align}
    0 = \tilde \Delta_1(k) = 1 - \cos\left(\mathfrak m(1) \kell \right) - \mathcal C_2^{(0,1)}(k).
\end{align}
In order to compute $\mathcal C_2^{(0,1)}(k)$, we use an interpolation function for $\mathfrak m(x)$, and rewrite
\begin{align}
    k\sum_{p=0}^n (-1)^p \int_{y_p}^{y_{p+1}} \mathfrak n(k,\xi)\,d\xi 
    = k \mathfrak g(k) \sum_{p=0}^n (-1)^p (\mathfrak m(y_{p+1}) - \mathfrak m(y_p))
    = \kell \left( \mathfrak m(1) - 2\sum_{p=0}^n (-1)^p \mathfrak m(y_p) \right)\!.
\end{align}
Then we use \eqref{eqn:Cn} to compute the $\tilde \Delta_1(k)$. We use a root finding algorithm to find the roots, using that
\begin{align}
    \partial_\kell \mathcal C_n^{(0,1)}(k) &= \frac{1}{2^n}\int_{0<\cdots<1} \left( \prod_{p=1}^n \eta(y_p) \right) \cos\left( \kell \left( \mathfrak m(1) - 2\sum_{p=0}^n (-1)^p \mathfrak m(y_p) \right) \right)\left( \mathfrak m(1) - 2\sum_{p=0}^n (-1)^p \mathfrak m(y_p) \right) d\mathbf{y}_n.
\end{align}
\begin{table}[tb]
    \begin{center}
        \begin{tabular}{|c|c|c|c|c|} \hline
            Method: & $\lambda_1$ & $\lambda_2$ & $\lambda_3$ & $\lambda_4$ \\ \hline
            chebfun  & $-41.585 + 3.3357 i$ & $-41.689 + 7.7171i$ & $-170.71 + 19.919i$ & $-170.62 + 23.463i$ \\
            NDEigenvalues & $-41.585 + 3.3364 i$ & $-41.689 + 7.7167 i$ & $-170.73 + 19.929 i$ & $-170.65 + 23.464 i$ \\
            { Hill's Method} & $-41.585 + 3.3358 i$ & $-41.689 + 7.7171 i$ & $-170.71 + 19.919 i$ & $-170.62 + 23.463 i$ \\
            FindRoot: $\Delta_0(k)$ & $-42.012 + 5.3928 i$ & $-42.012 + 5.3928 i$ & $-171.05 + 21.571 i$ & $-171.05 + 21.571 i$ \\
            FindRoot: $\Delta_1(k)$ & $-41.595 + 3.3501 i$ & $-41.671 + 7.7097 i$ & $-170.73 + 19.949 i$ & $-170.60 + 23.434 i$ \\
            FindRoot: $\Delta_2(k)$ & $-41.585 + 3.3356 i$ & $-41.689 + 7.7172 i$ & $-170.70 + 19.916 i$ & $-170.63 + 23.466 i$ \\ \hline
        \end{tabular}
    \caption{\label{tb:eigenvalues} Eigenvalues of the system \rf{eqn:eigenproblem} calculated using MATLAB's chebfun package, Mathematica's NDEigenvalues, Hill's method \cite{deconinck_kutz}, and a root finding algorithm on $\tilde \Delta_N(k)$ for $N=0,1,2$.}
    \end{center}
\end{table}
The results are shown in Table~\ref{tb:eigenvalues}.
%

%
%
\begin{appendices}
\section{Derivations} \label{sec:derivations}
In this appendix, we derive the solution expressions for the finite-interval, half-line, and whole-line IBVPs, in that order. The solution for each successive problem is obtained from the preceding one in a straightforward manner.
\subsection{The finite-interval problem}

\label{sec:derivations_fi}
\figcl{0.45}{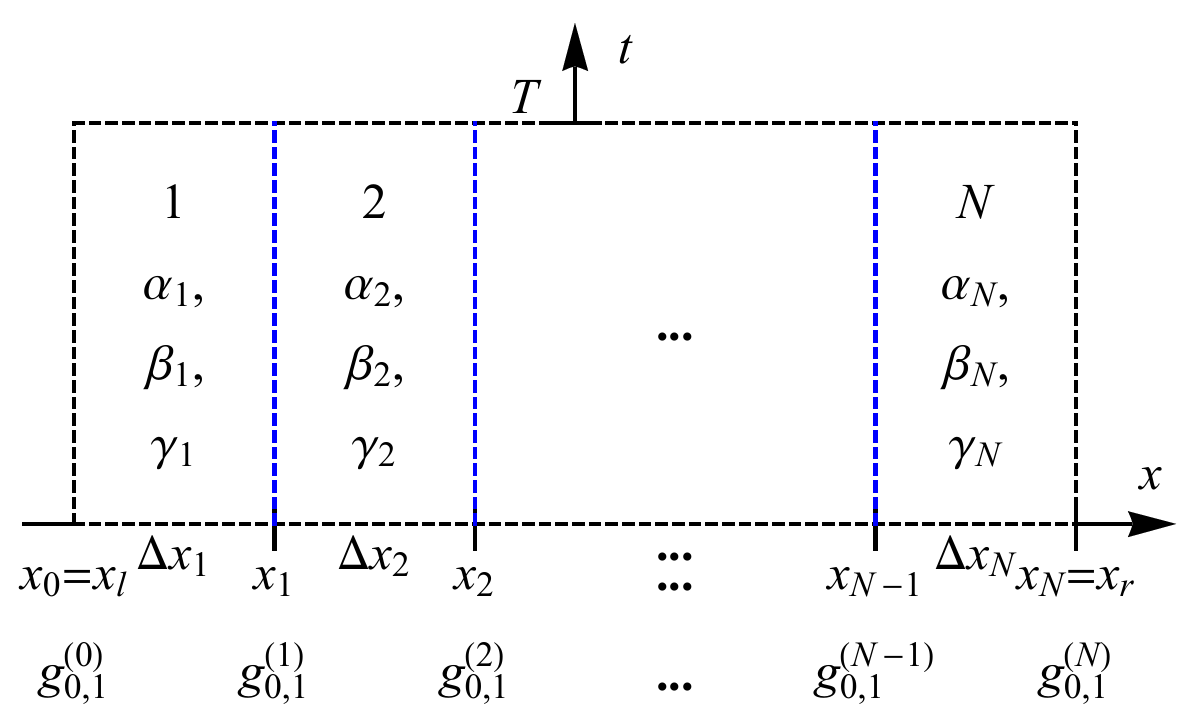}{A partition of the finite interval $[x_l,x_r]$.}{partition_fi}
To consider the finite-interval problem \eqref{eqn:IBVP_fi}, we form a partition $\{x_j, j=0, \ldots, N\}$ of the interval $[x_l,x_r]$, see Figure~\ref{fig:partition_fi}. {{For simplicity, we assume that the partition is evenly spaced, \textit{i.e.}, $\Delta x_j = \Delta x=(x_r-x_l)/N$ for $j=1,\ldots,N$, although this assumption may be relaxed easily.}} On each subinterval, we solve the evolution equation \eqref{eqn:PDE_fi} with constant-coefficient approximations $\alpha_j$, $\beta_j$, $\gamma_j$, $j=1, \ldots, N$ for $\alpha(x)$, $\beta(x),$ and $\gamma(x)$ (such that $\alpha_j\to \alpha(x_j)$, \etc, in the limit as $N\to\infty$), with the initial condition restricted to the subinterval. At each interface $x_j$, $j=1, \ldots, N-1$, we impose continuity of the solution and a jump discontinuity on the derivative, corresponding to the evolution equation, \textit{i.e.}, {{we solve the following interface problem:}}
\begin{subequations}
    \begin{align}
        \label{eqn:interface_PDE}
        q_t^{(j)} &= \alpha_j\beta_j q_{xx}^{(j)} +\gamma_j q^{(j)} + f(x,t), && x \in (x_{j-1}, x_{j}), \quad t>0,\quad j= 1, \ldots, N, \\ 
        \label{eqn:interface_IC}
        q^{(j)}(x,0) &= q_0(x), && x \in (x_{j-1}, x_{j}), \quad t>0,\quad j= 1, \ldots, N,  \\ 
        \label{eqn:interface_continuity} 
        q^{(j)}(x_{j}, t) &= q^{(j+1)}(x_{j}, t), && t>0, \quad j= 1, \ldots, N-1, \\ 
        \label{eqn:interface_jump} \beta_{j}q_x^{(j)} (x_{j},t) &= \beta_{j+1}q_x^{(j+1)}(x_{j}, t), && t>0, \quad j= 1, \ldots, N-1,
    \end{align}
\end{subequations}
\no with the boundary conditions
\begin{subequations} \label{eqn:interface_bc}
    \begin{align}
        a_{11}q^{(1)}(x_l,t)+a_{12}q_x^{(1)}(x_l,t)+b_{11}q^{(N)}(x_r,t) +b_{12}q_x^{(N)}(x_r,t) &=f_0(t), && t>0, \\ 
        a_{21}q^{(1)}(x_l,t)+a_{22}q_x^{(1)}(x_l,t) +b_{21}q^{(N)}(x_r,t) +b_{22}q_x^{(N)}(x_r,t)&=f_1(t), && t>0.
    \end{align}
\end{subequations}
\no The jump discontinuity in the derivative \eqref{eqn:interface_jump} can be derived by dividing the PDE \eqref{eqn:PDE_fi} by $\alpha(x)$ and integrating over a small interval containing $x_j$. 
Following \cite{interface_heat, interface_schrodinger, interface_heat_ring, interface_maps, interface_kdv, interface_dispersive}, we find the {\em local relations}
\beq
    \label{eqn:localrelation_fi}
    \left(e^{-i\kappa x+w_jt}q^{(j)}(x,t)\right)_t = \alpha_j\beta_j\left(e^{-i\kappa x+w_jt}\left( q_x^{(j)}(x,t) + i\kappa  q^{(j)}(x,t) \right)\right)_x+e^{-i\kappa x+w_jt}f(x,t),
\eeq
\no for $x \in (x_{j-1}, x_{j})$, $1 \leq j \leq N,$ and $w_j(\kappa) = \alpha_j\beta_j \kappa^2 - \gamma_j$. We define the ``transforms''
\begin{subequations} \label{eqn:transforms}
    \begin{align}
        \hat q_0^{(j)}(k)&= \frac{1}{\alpha_j}\int_{x_{j-1}}^{x_{j}} e^{-ik y} q_0(y) \, dy , && j= 1, \ldots, N, \\
        \hat q^{(j)}(k,t)&=\frac{1}{\alpha_j}\int_{x_{j-1}}^{x_{j}}e^{-ik y}q^{(j)}(y,t)\,dy,&& j= 1, \ldots, N, \\ 
        \tilde f_j(k,t)&= \frac{1}{\alpha_j}\int_0^t ds \int_{x_{j-1}}^{x_{j}} e^{-ik y+Ws}f(y,s)\, dy, && j= 1, \ldots, N,\\
        F_m(W,t)&= \int_{0}^t e^{Ws} f_m(s) \, ds, && m=0,1,\\
        g_m^{(j)}(W, t) &= \int_0^t e^{Ws} q_{mx}^{(j)}(x_{j}, s) \, ds, && j= 0, \ldots, N, \quad m =0,1,
    \end{align}
\end{subequations}
\no with $q_{mx}^{(0)}(x_l,t) = q_{mx}^{(1)}(x_l,t)$, for consistency at $j=0$. Using the interface conditions \eqref{eqn:interface_continuity} and \eqref{eqn:interface_jump}, we have
\begin{align}
    g_0^{(j)}(W,t) = \int_0^t e^{Ws} q^{(j+1)}(x_j,s) \, ds, ~~~~~~
    g_1^{(j)}(W,t) = \frac{\beta_{j+1}}{\beta_j} \int_0^t e^{Ws} q_x^{(j+1)}(x_j,s)\, ds, ~~~~~~ j= 0, \ldots, N-1,
\end{align}
\no where we define $\beta_0=\beta_1$, again for consistency. From the boundary conditions \eqref{eqn:interface_bc}, 
\begin{subequations} \label{eqn:BC}
    \begin{align}
        \label{eqn:BC0}
        a_{11}g_0^{(0)}(k^2,t)+a_{12}g_1^{(0)}(k^2,t)+b_{11}g_0^{(N)}(k^2,t)+b_{12} g_1^{(N)}(k^2,t)&=F_0(k^2,t),\\
        \label{eqn:BC1}
        a_{21}g_0^{(0)}(k^2,t)+a_{22}g_1^{(0)}(k^2,t)+b_{21}g_0^{(N)}(k^2,t)+b_{22} g_1^{(N)}(k^2,t)&=F_1(k^2,t).
    \end{align}
\end{subequations}
Integrating the local relations \eqref{eqn:localrelation_fi} over $D_j = (x_{j-1}, x_{j})\times (0,T)$, we find
\beq
    \alpha_j\tilde f_j(\kappa ,T) = \int_0^Tdt\int_{x_{j-1}}^{x_{j}}dx\, \left[\left(e^{-i\kappa x + w_jt} q^{(j)}\right)_t -\alpha_j\beta_j \left(e^{-i\kappa x+w_jt}\left(q_x^{(j)}+i\kappa  q^{(j)}\right)\right)_x\right]\!,
    ~~~~ j= 1, \ldots, N.
\eeq
\no Using Green's theorem, 
\begin{align} \nonumber
    \alpha_j \tilde f_j(\kappa ,T) &= \int_{x_{j-1}}^{x_{j}} e^{-i\kappa x} q_0(x) \, dx -e^{w_jT}\int_{x_{j-1}}^{x_{j}} e^{-i\kappa x} q^{(j)}(x,T) \,dx \\ \nonumber
    &~~~+ \alpha_j e^{-i\kappa x_{j}}\int_0^Te^{w_jt}\left(\beta_jq_x^{(j)}(x_{j}, t) + i\beta_j\kappa  q^{(j)}(x_{j}, t) \right)\, dt \\ 
    &~~~- \alpha_j e^{-i\kappa x_{j-1}}\int_0^T e^{w_jt}\left(\beta_jq_x^{(j)}(x_{j-1}, t) + i\beta_j\kappa  q^{(j)}(x_{j-1}, t)\right)\,dt,~~~~j= 1, \ldots, N,
\end{align}
\no which are rewritten as {\em global relations} using \eqref{eqn:transforms},
\begin{align} \label{eqn:GR} \nonumber 
    e^{w_jt}\hat q^{(j)}(\kappa ,t) &= \hat q_0^{(j)}(\kappa ) - \tilde f_j(\kappa ,t) + e^{-i\kappa x_{j}}\left(\beta_{j}g_1^{(j)}(w_j,t) + i\beta_{j}\kappa  g_0^{(j)}(w_j,t)\right) \\
    &~~~- e^{-i\kappa x_{j-1}} \left(\beta_{j-1}g_1^{(j-1)}(w_j,t) + i\beta_j\kappa   g_0^{(j-1)}(w_j,t)\right)\!, ~~~~~~ j= 1, \ldots, N.
\end{align}
As in \cite{Mantzavinos, interface_maps, interface_schrodinger}, it is convenient for the first arguments of $g_m^{(j)}(w_j,t)$ to be identical. We transform the independent variable $\kappa$ in the $j$th equation as 
\beq \label{eqn:nuj}
    \kappa  = \nu_j(k) =  \frac{k}{\sqrt{\alpha_j\beta_j}}\sqrt{1+\frac{\gamma_j}{k^2}}, ~~~~~~ j= 1, \ldots, N.
\eeq
We do not worry about the branch cuts here. The resulting branch cuts in the solution are defined in Section~\ref{sec:assumptions} and proven to be correct in Appendices~\ref{sec:proofs_welldefined}--\ref{sec:proofs_IC}.
However, since we assume that $\gamma(x)$ is bounded, there are no branch cuts for $|k|> \sqrt{M_\gamma}$ for $k \in \Omega$ where $M_\gamma>0$ is defined in Assumption~\enumref{ass:alphabetagamma}{enum:bounded}. Until we take the limit, we suppress the $k$ dependence of $\nu_j(k)$. Our global relations \eqref{eqn:GR} become 
\begin{align} \nonumber
    e^{k^2t}\hat q^{(j)}(\nu_j, t)  &= \hat q_0^{(j)}\left(\nu_j\right)- \tilde f_j(\nu_j,t) + e^{-i\nu_j x_{j}}\left( \beta_{j} g_1^{(j)}(k^2,t) +i\beta_j\nu_j g_0^{(j)}(k^2,t)\right) \\ 
    &~~~- e^{- i\nu_jx_{j-1}}\left( \beta_{j-1} g_1^{(j-1)}(k^2,t) + i\beta_j\nu_j g_0^{(j-1)}(k^2,t)\right)\!, ~~~~j= 1, \ldots, N.
\end{align}
\no These relations are valid for $k\in \C$, since the domains are bounded. Letting $k\mapsto -k$, (and $\nu_j\mapsto-\nu_j$), gives $2N$ equations, along with \eqref{eqn:BC} for $2N+2$ unknowns. We write this linear system of equations in matrix form as 
\beq
    \mathcal A_N(k) X_N(k^2,t) =  Y_N(k,t)-e^{k^2t} \mathcal Y_N(k, t),
\eeq
\no where 
\begin{subequations}
    \begin{align}
        X_N(k^2,t) &= \left(g_0^{(0)}(k^2,t), \ldots,  g_0^{(N)}(k^2,t), \beta_0g_1^{(0)}(k^2,t), \ldots, \beta_{N}g_1^{(N)} (k^2,t) \right)^\top\!\!\!\!\!,\\ \nonumber
        Y_N(k,t) &= \left(0,\hat q_0^{(1)}(\nu_1), \ldots, \hat q_0^{(N)}(\nu_N), \hat q_0^{(1)}(-\nu_1), \ldots, \hat q_0^{(N)}(-\nu_N),0\right)^\top\!\!\!\!\!\\
        &~~~- \left(-F_0(k^2,t), \, \tilde f_1(\nu_1,t), \ldots, \tilde f_{N}(\nu_N,t), \tilde f_1(-\nu_1,t), \ldots, \tilde f_{N}(-\nu_N,t), \, -F_1(k^2,t) \right)^\top\!\!\!\!\!, \\
        \mathcal Y_N(k,t) &= \left(0,\hat q^{(1)}\left(\nu_1,t\right),\ldots, \hat q^{(N)}\left(\nu_N,t\right), \hat q^{(1)}\left(-\nu_1,t\right), \ldots, \hat q^{(N)}\left(-\nu_N,t\right),0\right)^\top\!\!\!\!\!,
    \end{align}
\end{subequations}
\no and
\beq
    \mathcal A_N(k) = \left(\begin{NiceTabular}{ccccc:ccccc} 
        $a_{11}$ & 0 & $\cdots$ & 0 & $b_{11}$ & $a_{12}/\beta_0$ & 0 & $\cdots$ & 0 & $b_{12}/\beta_N$ \\ 
        \hdottedline  \\ 
        & & $A_{N,N+1}(k)$ & & & & & $B_{N,N+1}(k)$ & &  \\ \\ 
        \hdottedline \\ 
        & & $A_{N,N+1}(-k)$ & & & & & $B_{N,N+1}(-k)$ & & \\\\ 
        \hdottedline $a_{21}$ & 0 & $\cdots$ & 0 & $b_{21}$ & $a_{22}/\beta_0$ & 0 & $\cdots$ & 0 & $b_{22}/\beta_N $ 
    \end{NiceTabular}\right)\!,
\eeq
\no where 
\beq
    A_{N,N+1}(k) = \begin{pmatrix} 
        i \beta_1\nu_1 e^{-i\nu_1 x_0} & -i\beta_1\nu_1 e^{-i\nu_1x_1} & 0 & \cdots  & 0 \\ 
        \vdots & \ddots &  \ddots & \ddots & \vdots \\ 
        0 & \cdots & 0 & i\beta_N\nu_N e^{-i\nu_{N} x_{N-1}} & -i\beta_N\nu_N e^{-i\nu_{N} x_{N}} 
    \end{pmatrix}\!,
\eeq
\no and 
\beq
    B_{N,N+1}(k) = \begin{pmatrix} 
        e^{-i\nu_1 x_0} & -e^{-i\nu_1x_1} & 0 & \cdots  & 0  \\ 
        \vdots & \ddots &  \ddots & \ddots & \vdots \\ 
        0 & \cdots & 0 & e^{-i\nu_{N} x_{N-1}} & -e^{-i\nu_{N} x_{N}} 
    \end{pmatrix}\!.
\eeq
\no Here $A_{N,N+1}(k)$ and $B_{N,N+1}(k)$ are $N\times (N+1)$-dimensional matrices. 
Since the contribution involving $\mathcal Y_N$ along the contour $\partial \Omega$ (see below) is zero \cite{interface_maps}, it suffices to solve $\mathcal A_N X_N = Y_N$ for the unknown functions $g_m^{(j)}$. This is further justified in Appendices~\ref{sec:proofs_welldefined}--\ref{sec:proofs_IC}. Using Cramer's rule, 
\beq
    X_N^{(j)}(k^2,t) = g_0^{(j-1)}(k^2,t) = \frac{\det\big(\mathcal A_N^{(j)}(k)\big)}{\det\big(\mathcal A_N(k)\big)}, ~~~~~~j= 1, \ldots, N+1,
\eeq
\no where the matrix $\mathcal A_N^{(j)}(k)$ is $\mathcal A_N(k)$ with the $j$th column replaced by $Y_N$. If we multiply this equation by $k e^{-k^2t}$ and integrate over $\partial \Omega$, where $\Omega = \Omegadef$ for some $r>\sqrt{M_\gamma}$, see Figure~\ref{fig:Omega}, we can invert the time ``transform'' $g_0^{(j-1)}(k^2,t)$ \cite{JC_fokas_book}, to find 
\beq
q^{(j)}(x_{j-1},t) = \frac{1}{i\pi}\int_{\partial \Omega} \frac{\det(\mathcal A_N^{(j)}(k))}{\det(\mathcal A_N(k))} k e^{-k^2t}  \, dk, ~~~~~~ j= 1, \ldots, N+1.
\eeq
\no This gives the solution at the interface boundary points, which is all that is needed to consider the limit $q(x,t) = \lim_{N\to\infty} q^{(j)}(x_{j-1},t)$, {{where the $N$ dependence of $q^{(j)}(x_{j-1},t)$ is implicit}}.
Alternatively, we could compute the full solution of the interface problem as in \cite{interface_heat, interface_schrodinger, interface_heat_ring, interface_maps} and calculate that limit. This gives the same result. 
\par Define
\beq \label{eqn:DNdef}
    D_N(k) =  i^{ N+1}\frac{\det(\mathcal A_N(k))}{\nu_1\nu_N} \left(\prod_{p=1}^{N-1} \frac{1}{\Lambda_p^+} \right)\!,
\eeq
\no with 
\begin{align}
    \Lambda_p^{\boldsymbol{\ell}} = (\beta\nu)_{p+1} + (-1)^{\ell_p+\ell_{p+1}}(\beta\nu)_p
    ~~~~~ \text{ and } ~~~~~
    \Lambda_p^{\pm} = (\beta\nu)_{p+1} \pm (\beta\nu)_{p},
\end{align}
where $\boldsymbol{\ell} \in \{0,1\}^N$, so that $\ell_p, \ell_{p+1}\in \{0,1\}$ and $(\beta\nu)_j = \beta_j\nu_j$. Note that $\Lambda_p^{\boldsymbol{\ell}} = \Lambda_p^{+}$ when $\ell_p= \ell_{p+1}$ and $\Lambda_p^{\boldsymbol{\ell}} = \Lambda_p^-$ when $\ell_p\neq \ell_{p+1}$. For $N\leq 8$, we explicitly verify using Mathematica that
\begin{align}\nonumber \label{eqn:DN}
    D_N(k) &=  2i\left\{\frac{\beta_N(a:b)_{1,2}+\beta_1(a:b)_{3,4}}{\sqrt{(\beta\nu)_1}\sqrt{(\beta\nu)_N}}\frac{\sqrt{(\beta\nu)_N}}{\sqrt{(\beta\nu)_1}} \left(\prod_{p=1}^{N-1} 
    \frac{2(\beta\nu)_p}{\Lambda_p^+} 
    \right)  + (a:b)_{2,4} \mathfrak S_{N,1}^{(1,N)}(k)   + \frac{(a:b)_{1,3}}{\nu_1\nu_N} \mathfrak S_{N,0}^{(1,N)}(k)\right. \\
    &~~~\hspace{0.3in}\left.  + \frac{(a:b)_{1,4}}{\nu_1} \mathfrak C_{N,1}^{(1,N)}(k) - \frac{(a:b)_{2,3}}{\nu_N} \mathfrak C_{N,0}^{(1,N)}(k)  \right\}\!,
\end{align}
where $(a:b)_{i,j} = \det((a:b)_{\{1,2\},\{i,j\}})$ is the determinant of the minor of maximal size with columns at $i$ and $j$ \cite{LA_and_geom} of the concatenated matrix $(a:b)$. We define
\begin{subequations}\label{eqn:frakdef}
    \begin{align}
        \label{eqn:frakCndef}
        \mathfrak C_{N,\lambda}^{(q,s)}(k) = \sum_{\substack{\boldsymbol \ell \in \{0,1\}^{s-q+1} \\ \ell_q=0}} (-1)^{\lambda\ell_s} \left(\prod_{p=q}^{s-1} 
        \frac{\Lambda_p^{\boldsymbol{\ell}}}{ \Lambda_p^+}
        \right)  \cos\left(\sum_{p=q}^{s} (-1)^{\ell_p} \nu_p \Delta x\right)\!, \\
        \label{eqn:frakSndef}
        \mathfrak S_{N,\lambda}^{(q,s)}(k) = \sum_{\substack{\boldsymbol \ell \in \{0,1\}^{s-q+1} \\ \ell_q=0}} (-1)^{\lambda \ell_s} \left(\prod_{p=q}^{s-1} 
        \frac{\Lambda_p^{\boldsymbol{\ell}}}{\Lambda_p^+} \right) \sin\left(\sum_{p=q}^{s} (-1)^{\ell_p} \nu_p \Delta x\right)\!,
    \end{align}
\end{subequations}
\no where $\lambda = 0,1$. We do not prove this result for general $N$. Its justification follows indirectly from the proofs in Appendices~\ref{sec:proofs_welldefined}--\ref{sec:proofs_IC}.
We can show that
\beq \label{eqn:prod_asymptotics1}
    \prod_{p=1}^{N-1} \frac{2(\beta\nu)_p}{\Lambda_p^+} = \exp\left(\sum_{p=1}^{N-1}\left( \ln(2(\beta\nu)_p) - \ln(\Lambda_p^+)\right) \right) = \frac{\sqrt{(\beta\nu)_1}}{\sqrt{(\beta\nu)_N}} + \bigoh(\Delta x), 
\eeq
\no as $N\to\infty$ and $\Delta x\to 0^+$.
Similarly, 
\beq
    \label{eqn:MVT&prod_asymptotics}
    \frac{\Lambda_p^-}{\Lambda_p^+} = \frac12 \frac{(\beta\mathfrak n)'(k,x_p)}{(\beta\mathfrak n)(k,x_p)} \Delta x + \bigoh ((\Delta x)^2) 
    \qquad \text{ and } \qquad 
    \prod_{p=\ell}^{m-1} \frac{2(\beta\nu)_p}{\Lambda_p^+} = \frac{\sqrt{(\beta\mathfrak n)(k,x_\ell)}}{\sqrt{(\beta\mathfrak n)(k,x_m)}} + \bigoh(\Delta x),
\eeq
as $\ell, m, N\to\infty$, $\Delta x \to 0^+$, and where the prime denotes the derivative with respect to the spatial variable, and $\mathfrak n(k,x)$, $(\beta\mathfrak n)(k,x)$ are defined in Definition~\ref{def:mu}.
Note that to use \eqref{eqn:MVT&prod_asymptotics}, we assume $(\beta\mathfrak n)(k,x)$ is a smooth function of $x$. If this function has a countable number of discontinuities, it is possible to proceed, but we have to account for the jumps. 
 
We wish to consider \eqref{eqn:DN} as $N\to\infty$ $(${\em i.e.,} $\Delta x \to 0^+)$. To this end, we break up the sum in \eqref{eqn:frakCndef} by the number of times $n$ the entries of the vector $\boldsymbol{\ell}=(\ell_q,\ldots, \ell_s)$ switch from 0 to 1 or from 1 to 0, \textit{e.g.}, $(0,\ldots, 0, 1, \ldots, 1)$ switches once, so $n=1$. We sum over where the possible switches of each order $n$ can occur \textit{i.e.}, $q-1< y_1 < y_2 < \cdots < y_n < s$. At the location of each switch, $\Lambda_p^{\boldsymbol{\ell}}/\Lambda_p^+ = \Lambda_p^- / \Lambda_p^+$, whereas $\Lambda_p^{\boldsymbol{\ell}}/\Lambda_p^+ =1$ otherwise. 
Defining $y_0=q-1$ and $y_{n+1}=s$, this gives
\begin{align}
    \mathfrak C_{N,\lambda}^{(q,s)}(k) &= \sum_{n=0}^{s-q} \sum_{\stackbin{~}{y_0< y_1 < \cdots < y_n < y_{n+1}}} (-1)^{\lambda n} \left(\prod_{p=1}^{n} \frac{\Lambda_{y_p}^{-}}{ \Lambda_{y_p}^+}\right)  \cos\left(\sum_{p=0}^n (-1)^p \sum_{r=y_p+1}^{y_{p+1}} \nu_p \Delta x\right)\!.
\end{align}
Using \eqref{eqn:MVT&prod_asymptotics}, we arrive at a sum of $n$-dimensional Riemann sums which limit to $n$-dimensional integrals, giving 
\beq \label{eqn:Cnlimit}
    \mathfrak C_{N,\lambda}^{(q,s)}(k) = \sum_{n=0}^\infty (-1)^{\lambda n} \mathcal C_n^{(x_q,x_s)}(k) + \bigoh(\Delta x), 
\eeq
where $\mathcal C_n^{(a,b)}(k)$ is defined in \eqref{eqn:Cn}. The limit of \eqref{eqn:frakSndef} is
\beq \label{eqn:Snlimit}
    \mathfrak S_{N,\lambda}^{(q,s)}(k) = \sum_{n=0}^\infty (-1)^{\lambda n} \mathcal S_n^{(x_q,x_s)}(k)+ \bigoh(\Delta x),
\eeq
obtained the same way, with $\mathcal S_n^{(a,b)}(k)$ defined in \eqref{eqn:Sn}. No more rigor is required at this point, as we prove in Appendices~\ref{sec:proofs_welldefined}--\ref{sec:proofs_IC} that our result is a solution under less restrictive assumptions needed to justify these steps.

Using \eqref{eqn:prod_asymptotics1}, \eqref{eqn:Cnlimit}, and \eqref{eqn:Snlimit} in \eqref{eqn:DN}, we have that 
\beq
    \Delta(k) = \lim_{N\to\infty} \Xi(k) D_N(k),
\eeq
gives \eqref{eqn:Delta_fi}.
For the numerator, similar to $D_N(k)$ in \eqref{eqn:DNdef}, we define
\beq
    \label{eqn:ENdef}
    E_N(k,j,t) = i^N\frac{2\det(\mathcal A_N^{(j)}(k))}{\nu_1\nu_N}\left(\prod_{p=1}^{N-1} \frac{1}{\Lambda_p^+} \right)\! , 
\eeq
\no and use a cofactor expansion along the $j$th column of $\mathcal A_N^{(j)}$, so that
\beq
    E_N(k,j,t) = \sum_{m=1}^{N+1} Y_N^{(m)}M_N^{(m,j)}(k) + \sum_{m=1}^{N+1} Y_N^{(m+N+1)} M_N^{(m+N+1,j)}(k),
\eeq
\no where $M_N^{(m,j)}(k)$ are cofactors of the matrix $\mathcal A_N^{(j)}$, scaled by the same factor as in \eqref{eqn:ENdef}. {{With $x=x_j=j\Delta x = j/N$, which is fixed}}, we let
\beq
    \label{eqn:B0}
    \mathcal B_0(k,x) = k \Xi(k) \lim_{N\to\infty} M_N^{(1,j)}(k) \qquad \text{ and } \qquad  \mathcal B_1(k,x) = k \Xi(k) \lim_{N\to\infty} M_N^{(2N+2,j)}(k).
\eeq
\no Since, for $m=1, \ldots, 2N$,  
\beq
    Y_N^{(m+1)} = \hat q_0^{(m)}(\nu_{m}) - \tilde f_{m}(\nu_m,t) = \frac{e^{-i\nu_mx_m}}{\alpha_m} \left(q_0(x_m) - \int_0^t f(x_m,s) e^{k^2s} \, ds \right)\Delta x + \bigoh ((\Delta x)^2),
\eeq
\no so that, for $m=1, \ldots, N$, 
\beq
    Y_N^{(m+1)} = \frac{e^{-i\nu_mx_m} \psi_N^{(m)}(k^2,t)}{\alpha_m} \, \Delta x + \bigoh((\Delta x)^2)
    \quad \text{ and } \quad 
    Y_N^{(m+N+1)} = \frac{e^{i\nu_mx_m}\psi_N^{(m)}(k^2,t)}{\alpha_m} \, \Delta x + \bigoh((\Delta x)^2),
\eeq
\no which defines $\psi_N^{(m)}(k^2,t)$.
Then
\begin{align}\label{eqn:Phi_Nlimit}
    \Phi(k,x,t) &= \lim_{N\to \infty} k\Xi(k) E_N(k,j,t),
\end{align}
which gives \eqref{eqn:Phi_fi}. Here
\begin{align}  \label{eqn:Phi_psi}
    \Phi_\psi(k,x,t) &= \lim_{N\to\infty}k \Xi(k) \sum_{m=1}^{N}\frac{\psi_N^{(m)}(k^2,t)}{\alpha_m}\left( e^{-i\nu_mx_m} M_N^{(m+1,j)}(k) + e^{i\nu_mx_m} M_N^{(m+N+1,j)}(k) \right) \Delta x,
\end{align}
\no where we let $y=x_m = m\Delta x = m/N$, which is kept fixed. This gives \eqref{eqn:Phi_psi_fi_def}, after defining
\begin{subequations}
    \begin{align}
        \Psi(k,x,y) &= \lim_{N\to\infty} \Psi_N^{(j,m)}(k) = \lim_{N\to\infty} \Xi(k) \sqrt{(\beta\nu)_m} \sqrt{(\beta\nu)_j} \left( e^{-i\nu_m x_m} M_N^{(m+1,j)}(k) + e^{i \nu_mx_m }M_N^{(m+N+1,j)}(k) \right)\!, \\
        \psi_\alpha(k^2,y,t) &= \lim_{N\to\infty} \frac{\psi_N^{(m)}(k^2,t)}{\alpha_m} = \frac{q_0(y)}{\alpha(y)} - \int_0^t \frac{f(y,s)}{\alpha(y)} e^{k^2s} \, ds,
    \end{align}
\end{subequations}
\no which defines $\Psi_N^{(j,m)}$.

For the boundary term $\mathcal B_0(k,x)$, similar to \eqref{eqn:DN}, we explicitly verify using Mathematica for
$N\leq 8$ and for $1= m < j\leq N$,
\begin{align} \nonumber
    M_N^{(1,j)}(k) &= \frac{4}{\sqrt{(\beta \nu)_{j-1}}}\left\{ \frac{\beta_N}{\sqrt{(\beta \nu)_{j-1}}} \left(\prod_{p=j-1}^{N-1} \frac{2(\beta\nu)_p}{\Lambda_p^+}\right) \left[ -\frac{a_{21}}{\nu_1} \mathfrak S_{N,0}^{(1,j-1)}(k) +a_{22} \mathfrak C_{N,0}^{(1,j-1)}(k) \right] \right. \\
    &\hspace{1in} + \left. \frac{\sqrt{(\beta\nu)_{j-1}}}{\nu_1} \left(\prod_{p=1}^{j-1} \frac{2(\beta\nu)_p}{\Lambda_p^+}\right) \left[ \frac{b_{21}}{\nu_N} \mathfrak S_{N,0}^{(j,N)}(k) + b_{22} \mathfrak C_{N,1}^{(j,N)}(k)\right] \right\}\!,
\end{align}
\no and using \eqref{eqn:MVT&prod_asymptotics}, \eqref{eqn:Cnlimit}, \eqref{eqn:Snlimit}, and \eqref{eqn:B0}, we find \eqref{eqn:Bm} for $j=2$.
\no Similarly, for the other boundary term $\mathcal B_1(k,x)$, we find \eqref{eqn:Bm} for $j=1$.

For the remaining terms, for $1\leq m < j \leq N$, we derive
{ \begin{align} \nonumber
    \Psi_N^{(j,m)}(k) &= 4\Xi(k) \frac{\sqrt{(\beta\nu)_j}}{\sqrt{(\beta\nu)_m}} \left(\prod_{p=m}^{j-1} 
    \frac{2(\beta\nu)_p}{\Lambda_p^+}
    \right)  \bigg[ - (a:b)_{2,4}\mathfrak C_{N,0}^{(1,m)}(k) \mathfrak C_{N,1}^{(j,N)}(k) +  \frac{(a:b)_{1,3}}{\nu_1\nu_N}  \mathfrak S_{N,0}^{(1,m)}(k) \mathfrak S_{N,0}^{(j,N)}(k) \\ \nonumber
    &~~~\hspace{2in}+ \frac{(a:b)_{1,4}}{\nu_1} \mathfrak S_{N,0}^{(1,m)}(k) \mathfrak C_{N,1}^{(j,N)}(k) - \frac{(a:b)_{2,3}}{\nu_N} \mathfrak C_{N,0}^{(1,m)}(k) \mathfrak S_{N,0}^{(j,N)}(k) \bigg] \\
    &~~~- \frac{4(a:b)_{1,2}\beta_N}{\sqrt{(\beta\nu)_1}\sqrt{(\beta\nu)_N}} \frac{\sqrt{(\beta\nu)_j}}{\sqrt{(\beta\nu)_{j-1}}} \frac{\sqrt{(\beta\nu)_m}}{\sqrt{(\beta\nu)_1}} \left(\prod_{p=1}^{m} \frac{(\beta\nu)_p}{\Lambda_p^+} \right) \frac{\sqrt{(\beta\nu)_N}}{\sqrt{(\beta\nu)_{j-1}}}
    \left(\prod_{p=j-1}^{N-1} \frac{(\beta\nu)_p}{\Lambda_p^+} \right) \mathfrak S_{N,0}^{(m+1,j-1)}(k).
\end{align}}
\no Taking the limit using \eqref{eqn:MVT&prod_asymptotics}, \eqref{eqn:Cnlimit}, and \eqref{eqn:Snlimit}, as before letting $x_j\to x$ and $x_m\to y$, we find for $x_r \leq y<x \leq x_\ell$,
\begin{align} \nonumber
    \Psi(k,x,y) &= 4\Xi(k)\left\{ -(a:b)_{2,4} \left( \sum_{n=0}^\infty \mathcal C_n^{(x_l,y)} \right) \left( \sum_{n=0}^\infty (-1)^n \mathcal C_n^{(x,x_r)} \right) + \frac{(a:b)_{1,3}}{k^2\mathfrak n(k,x_l)\mathfrak n(k,x_r)} \left( \sum_{n=0}^\infty \mathcal S_n^{(x_l,y)} \right) \left( \sum_{n=0}^\infty \mathcal S_n^{(x,x_r)} \right) \right. \\ \nonumber
    &~~~\hspace{0.5in}+ \frac{(a:b)_{1,4}}{k\mathfrak n(k,x_l)}\left( \sum_{n=0}^\infty \mathcal S_n^{(x_l,y)} \right) \left( \sum_{n=0}^\infty (-1)^n \mathcal C_n^{(x,x_r)} \right) - \frac{(a:b)_{2,3}}{k\mathfrak n(k,x_r)} \left( \sum_{n=0}^\infty \mathcal C_n^{(x_l,y)} \right) \left( \sum_{n=0}^\infty \mathcal S_n^{(x,x_r)} \right)\\
    &~~~\hspace{0.5in}\left.-\frac{(a:b)_{1,2} \beta(x_r)}{k\sqrt{(\beta\mathfrak n)(k,x_l)}\sqrt{(\beta\mathfrak n)(k,x_r)}} \sum_{n=0}^\infty \mathcal S_n^{(y,x)} \right\}\!,
\end{align}
which may be rewritten as \eqref{eqn:Psi_fi1}.
\no Similarly, for $x_l\leq x<y\leq x_r$, we find \eqref{eqn:Psi_fi2}.
Finally, we have
\beq
    q(x,t) = \lim_{N\to \infty} \frac{1}{i\pi}\int_{\partial \Omega} \frac{\det(\mathcal A_N^{(j)}(k))}{\det(\mathcal A_N(k))} k e^{-k^2t}  \, dk= \lim_{N\to \infty} \frac{1}{2\pi}\int_{\partial \Omega} \frac{k\Xi(k) E_N(k,j,t)}{\Xi(k) D_N(k)}  e^{-k^2t}  \, dk,
\eeq
which gives \eqref{eqn:q_fi}.
\subsection{The half-line problem}

\label{sec:derivations_hl}
We obtain the solution of the half-line problem by taking the limit as $x_r\to \infty$ of the solution of the finite-interval problem \eqref{eqn:IBVP_fi} with $f_1(t)=0$ and
\beq \label{eqn:ab_hl_derivation}
    (a:b) = \begin{pmatrix} a_0 & a_1 & 0 & 0 \\ 0 & 0 & 1 & 0 \end{pmatrix}\!. 
\eeq
In this limit, \eqref{eqn:q_fi} becomes \eqref{eqn:q_hl} with the same $\Omega$, shown in Figure~\ref{fig:Omega}. This process is detailed below.

Using \eqref{eqn:ab_hl_derivation}, we find the coefficients $\mathfrak a(k) = 0$, $\mathfrak c_n(k) = -a_1/(k\mathfrak n(k,x_r)),$ and $\mathfrak s_n(k) = a_0/(k^2\mathfrak n(k,x_l) \mathfrak n(k,x_r))$. We define $\tilde \Delta(k) = k \mathfrak n(k,x_r) \Delta(k)$,  and \eqref{eqn:Delta_fi} becomes
\beq
    \tilde \Delta(k)  =2i\Xi(k) \left\{ \frac{a_0}{k\mathfrak n(k,x_l)} \sum_{n=0}^\infty \mathcal S_n^{(x_l,x_r)}(k) -a_1 \sum_{n=0}^\infty \mathcal C_n^{(x_l,x_r)}(k)  \right\}\!.
\eeq
We have $\mathcal B_1(k,x)=0$, and \eqref{eqn:Bm} for $j=2$ becomes
\beq
    \tilde{\mathcal B}_0(k,x) = k\mathfrak n(k,x_r) \mathcal B_0(k,x) = \frac{4\beta(x_l) \Xi(k)}{\sqrt{(\beta\mathfrak n)(k,x_l)}\sqrt{(\beta \mathfrak n)(k,x)}} \sum_{n=0}^\infty \mathcal S_n^{(x,x_r)}(k).
\eeq
\no For $x_l\leq y<x\leq x_r$, \eqref{eqn:Psi_fi1} becomes
\begin{align}
    \tilde{\Psi}(k,x,y) = k\mathfrak n(k,x_r)\Psi(k,x,y) = 4\Xi(k) \left\{ \frac{a_0}{k\mathfrak n(k,x_l)}\sum_{n=0}^\infty \sum_{\ell=0}^n \mathcal S_\ell^{(x,x_r)}(k) \mathcal S_{n-\ell}^{(x_l,y)}(k) - a_1 \sum_{n=0}^\infty \sum_{\ell=0}^n \mathcal S_\ell^{(x,x_r)}(k) \mathcal C_{n-\ell}^{(x_l,y)}(k)\right\}\!, 
\end{align}
\no and, for $x_l<x<y< x_r$, $\tilde \Psi(k,x,y) = \tilde \Psi(k,y,x)$. From \eqref{eqn:Phi_fi},
\beq
    \tilde\Phi(k,x,t) = k \mathfrak n(k,x_r) \Phi(k,x,t)= \tilde{\mathcal B}_0(k,x) + \tilde \Phi_\psi(k,x,t),
\eeq
where
\beq
    \tilde \Phi_\psi(k,x,t) = k \mathfrak n(k,x_r) \Phi_\psi(k,x,t) = \int_{x_l}^{x_r} \frac{\tilde \Psi(k,x,y)\psi_\alpha(k^2,y,t)}{\sqrt{(\beta \mathfrak n)(k,x)(\beta \mathfrak n)(k,y)}} \, dy,
\eeq
\no and $\psi_\alpha(k^2,y,t)$ is defined in Definition~\ref{def:mu}.

To take the limit as $x_r\to\infty$, using \eqref{eqn:Sn}, we write
\begin{align} \nonumber \label{eqn:Cnhl}
    \etwo{\int_a^{x_r} ik \mathfrak n(k,\xi) \, d\xi} \mathcal S_n^{(a,x_r)}(k) &= \frac{1}{2i\cdot 2^{n}} \int_{a<\cdots <x_r} \left(\prod_{p=1}^n \frac{(\beta\mathfrak n)'(k,y_p)}{(\beta \mathfrak n)(k,y_p)} \right) \left[ \exp\left(\sum_{p=0}^n (1+(-1)^p) \int_{y_p}^{y_{p+1}} ik \mathfrak n(k,\xi) \, d\xi \right) \right. \\
    & \hspace*{1in} \left. - \exp\left(\sum_{p=0}^n (1-(-1)^p) \int_{y_p}^{y_{p+1}} ik \mathfrak n(k,\xi) \, d\xi \right) \right] dy_1\cdots dy_n.
\end{align}
Since $\Re(ik\mathfrak n(k,x))< 0$ for all $k \in \Omega$ and all $x>x_l$, see Lemma~\ref{lem:Re_inu} in Section~\ref{sec:proofs_welldefined}, it follows that
\begin{align}
    \exp\left(\int_{y_n}^{x_r} ik \mathfrak n(k,\xi) \, d\xi \right) \to 0,
\end{align}
as $x_r \to \infty$. Thus the term in \eqref{eqn:Cnhl} which does not contain the $p=n$ term survives and considering even and odd $n$ separately, we conclude
\begin{align}
    \exp\left( \int_a^{x_r} ik \mathfrak n(k,\xi) \, d\xi\right) \mathcal S_n^{(a,x_r)}(k) \to -\frac{(-1)^n}{2i}{\mathcal E}_n^{(a,\infty)}(k),~~\mbox{as}~~ x_r\to \infty,
\end{align}
where ${\mathcal E}_n^{(a,b)}(k)$ is defined in \eqref{eqn:Entilde}.
\no Similarly, 
\begin{align}
    \exp\left( \int_a^{x_r} ik \mathfrak n(k,\xi) \, d\xi \right) \mathcal C_n^{(a,x_r)}(k) &\to \frac12 {\mathcal E}_n^{(a,\infty)}(k), ~~\mbox{as}~~ x_r\to \infty.
\end{align}
Therefore, we have as $x_r\to \infty$,
\beq
    -2i \tilde \Delta(k) \to  2 \sum_{n=0}^\infty \left(  \frac{(-1)^n ia_0}{k \mathfrak n(k,x_l)}-a_1\right){\mathcal E}_n^{(x_l,\infty)}(k),
\eeq
\no and
\beq
    -2i\tilde{\mathcal B}_0(k,x) \to  \frac{4\beta(x_l) \etwo{\int_{x_l}^{x} ik \mathfrak n(k,\xi) \, d\xi}}{\sqrt{(\beta\mathfrak n)(k,x_l)} \sqrt{(\beta \mathfrak n)(k,x)}}  \sum_{n=0}^\infty (-1)^n {\mathcal E}_n^{(x,\infty)}(k),
\eeq
\no and, for $x_l\leq y<x$,
\beq
    -2i\tilde{\Psi}(k,x,y) \to 4\etwo{\int_{x_l}^{x}i \nu(k,\xi) \, d\xi} \sum_{n=0}^\infty \sum_{n=0}^n (-1)^\ell  \left(\frac{a_0}{k \mathfrak n(k,x_l)}\mathcal S_{n-\ell}^{(x_l,y)}(k)- a_1\mathcal C_{n-\ell}^{(x_l,y)}(k) \right) {\mathcal E}_\ell^{(x,\infty)}(k),
\eeq
\no and similarly for $x_l\leq x<y $. These final results combine to give \eqref{eqn:q_hl}. 
\subsection{The whole-line problem}

\label{sec:derivations_wl}
We repeat the process from the previous section, now letting  $x_l\to -\infty$. Starting from the half-line solution \eqref{eqn:q_hl} with $f_0(t)=0$, $a_0=1$, and $a_1=0$. 
The denominator in \eqref{eqn:q_hl} is determined by
\beq
    k\mathfrak n(k,x_l)\Delta(k) = 2i\sum_{n=0}^\infty (-1)^n {\mathcal E}_{n}^{(x_l,\infty)}(k) .
\eeq
\no Since $\mathcal B_0(k,x,t) = 0$, we also have from \eqref{eqn:Phi_hl}
\beq
    \Phi(k,x,t) = \int_{-\infty}^{\infty} \frac{\Psi(k,x,y)\psi_\alpha(k^2,y,t)}{\sqrt{(\beta\mathfrak n)(k,x)(\beta\mathfrak n)(k,y)}} \, dy,
\eeq
\no where $\psi_\alpha(k^2,y,t)$ is defined in Definition~\ref{def:mu}. For $x_l<y<x$, \eqref{eqn:Psi_hl} becomes
\beq
    k \mathfrak n(k,x_l)\Psi(k,x,y) = 4\exp\left(\int_{x_l}^{x}ik \mathfrak n(k,\xi) \, d\xi\right) \sum_{n=0}^\infty \sum_{\ell=0}^n (-1)^\ell   \mathcal S_{n-\ell}^{(x_l,y)} (k){\mathcal E}_\ell^{(x,\infty)}(k),
\eeq
\no and $\Psi(k,x,y) =\Psi(k,y,x)$ for $x_l<x<y$.
Since ${\mathcal E}_n^{(x_l,\infty)}(k) \to 0$ if $n$ is odd, and
from \eqref{eqn:Cnhl}, we have
\beq
    \etwo{\int_{x_l}^{b}ik \mathfrak n(k,\xi) \, d\xi}\mathcal S_n^{(x_l,b)}(k) \to - \frac1{2i} \tilde{\mathcal E}_n^{(-\infty,b)}(k),
\eeq
\no as $x_l\to-\infty$. Therefore, 
\beq
    k \mathfrak n(k,x_l)\Delta(k) \to 2i \sum_{\substack{n=0\\n\text{ even}}}^\infty  {\mathcal E}_n^{(-\infty,\infty)}(k).
\eeq
\no For $x_l<y<x$,
\beq
    k\mathfrak n(k,x_l)\Psi(k,x,y) = 2i\exp\left( \int_{y}^x ik \mathfrak n(k,\xi)\, d\xi \right) \sum_{n=0}^\infty \sum_{\ell=0}^n (-1)^\ell \tilde {\mathcal E}_{n-\ell}^{(-\infty,y)} (k) {\mathcal E}_\ell^{(x,\infty)}(k),
\eeq
\no and $\Psi(k,x,y) =\Psi(k,y,x)$ for $x_l<x<y$. Combining these results gives \eqref{eqn:q_wl}.
\section{Proofs: the solution expressions are well defined} \label{sec:proofs_welldefined}
\newcommand{\Jn}{{\mathcal{J}}_n^{(a,b)}[\sigma_{p,n}](k)}
\newcommand{\Jnp}[1]{{\mathcal{J}}_n^{(a,b)}[#1](k)}
\newcommand{\En}{{\mathcal{E}}_n^{(a,b)}(k)}
\newcommand{\Entilde}{\tilde{\mathcal{E}}_n^{(a,b)}(k)}
\newcommand{\Enbar}{\bar{\mathcal{E}}_n^{(a,b)}(k)}
\newcommand{\argi}{\arg}
Prior to proving that the solution expression \eqref{eqn:q_all} solves the evolution equation \ref{eqn:evolution_equation}  and satisfies the initial and boundary conditions for the problem considered, we show in this appendix that this expression is well defined for all problems considered. We refer to the whole-line, half-line, and {\em regular} finite-interval problems as {\em regular problems}, and the {\em irregular} finite-interval problems as {\em irregular problems}. Throughout, we need Assumptions~\ref{ass:alphabetagamma}~and~\ref{ass:ffunctions} from Section~\ref{sec:assumptions}. For Boundary Case \ref{enum:BC4} of the finite-interval problem, we also require Assumption~\ref{ass:alphabetagamma2}. Note that Assumption~\enumref{ass:alphabetagamma2}{enum:AC_extra} is not needed for all {\em irregular problems}, only for Boundary Case \ref{enum:BC4}. Therefore, we will be explicit as to when Assumption~\ref{ass:alphabetagamma2} is required.

In this appendix, we denote the $r$ dependence of $\Omega$ explicitly as \mbox{$\Omega(r) = \{k\in \C: |k| >r $ and $\pi/4 < \arg(k) < 3\pi/4\}$.} We define $\arg(\,\cdot\,) \in [-\pi/2,3\pi/2) $ with $\theta = \arg(k)$. The following lemma characterizes some properties of the coefficient functions $\alpha,\beta$ that follow from the assumptions.
\begin{lemma} \label{lem:alphabetabounds}
    If $\alpha(x) \beta(x)$ is not identically zero, the following are equivalent:
    \begin{enumerate}
        \item[a.] $\alpha\beta \in L^\infty(\D)$, 
        \item[b.] $\alpha \in L^\infty(\D)$, 
        \item[c.] $\beta \in L^\infty(\D)$, 
    \end{enumerate}
    as are the following:
    \begin{enumerate}
        \item[i.] $m_{\alpha\beta} = \inf_{x\in \D}|\alpha(x)\beta(x)| >0$, 
        \item[ii.] $m_\alpha = \inf_{x\in \D}|\alpha(x)| >0$, 
        \item[iii.] $m_\beta = \inf_{x\in \D}|\beta(x)| >0$.
    \end{enumerate}
\end{lemma}

\begin{proof}
    Under Assumption ~\enumref{ass:alphabetagamma}{enum:lowerbound}, there exists an $x_0 \in \D$ such that $0 < |\alpha(x_0)\beta(x_0)| < \infty$. From    Assumptions~\enumref{ass:alphabetagamma}{enum:AC}~and~\enumref{ass:alphabetagamma}{enum:L1},
    \begin{subequations}
        \begin{align}
            \left|\frac{\beta(x)}{\beta(x_0)}\right| 
            &= \left|\frac{\alpha(x)}{\alpha(x_0)} \exp\left(\int_{x_0}^x \frac{\beta'(y)}{\beta(y)} - \frac{\alpha'(y)}{\alpha(y)} \, dy\right)\right| 
            \leq \left|\frac{\alpha(x)}{\alpha(x_0)} \right| \exp\left( \left\|\frac{\beta'}{\beta} - \frac{\alpha'}{\alpha}\right\|_\D \right) = E\left|\frac{\alpha(x)}{\alpha(x_0)} \right|\!,
        \end{align}
        with $E = \exp(\|\beta'/\beta - \alpha'/\alpha\|_\D)$. 
        We conclude that $b \Rightarrow c$, $b\Rightarrow a$, $a\Rightarrow c$, $iii\Rightarrow ii$, $iii\Rightarrow i$, and $i\Rightarrow ii$. Similarly,
        \begin{align}
            \left|\frac{\alpha(x)}{\alpha(x_0)}\right|  \leq  E\left|\frac{\beta(x)}{\beta(x_0)} \right|\!,
        \end{align}
    \end{subequations}
    so that $c\Rightarrow b$, $c\Rightarrow a$, $a\Rightarrow b$, $ii\Rightarrow iii$, $ii\Rightarrow i$, and $i\Rightarrow iii$.
\end{proof}
Next, Lemmas~\ref{lem:nubounds}--\ref{lem:Re_inu} present some properties of the functions $\mathfrak n(k,x)$ and $(\beta\mathfrak n)(k,x)$.
\begin{lemma} \label{lem:nubounds}
    For $|k|\geq r>\sqrt{M_\gamma}$, where $M_\gamma$ is defined in Assumption~\enumref{ass:alphabetagamma}{enum:bounded}, we have
    \begin{align}
        m_{\mathfrak n} = \frac{1}{\sqrt{M_{\alpha\beta}}} \sqrt{1 - \frac{M_\gamma}{r^2}} \leq \left|\mathfrak n(k,x)\right| \leq \frac{1}{\sqrt{m_{\alpha\beta}}} \sqrt{1 + \frac{M_\gamma}{r^2}} = M_{\mathfrak n},
    \end{align}
    \no which defines $m_{\mathfrak n},M_{\mathfrak n}>0$. From this, we also have $m_{\mathfrak n} \leq |\mu(x)| \leq M_{\mathfrak n}$.
\end{lemma}
\begin{proof} 
    The proof is trivial from the definition of $\mathfrak n(k,x)$ in Definition~\ref{def:mu} using Assumptions~\enumref{ass:alphabetagamma}{enum:lowerbound}~and~\enumref{ass:alphabetagamma}{enum:bounded}.
\end{proof}
\begin{lemma} \label{lem:beta_nu}  
    For $|k|\geq r>\sqrt{M_\gamma}$, $(\beta\mathfrak n)'/(\beta \mathfrak n) \in L^{1}(\D)$, and under Assumption~\ref{ass:alphabetagamma2}, $\mathfrak u \in \mathrm{AC}(\D)$.
\end{lemma}
\begin{proof}
    The function
    \begin{align} \label{eqn:betanu_expansion}
        \fracbetanuargs{x} = \frac12\left(\frac{\beta'(x)}{\beta(x)}-\frac{\alpha'(x)}{\alpha(x)} + \frac{\gamma'(x)}{k^2+\gamma(x)}\right) \in L^1(\D),
    \end{align}
    for $|k|\geq r>\sqrt{M_\gamma}$, by Assumption~\enumref{ass:alphabetagamma}{enum:L1}.
    On the finite interval, by Assumptions~\enumref{ass:alphabetagamma}{enum:AC}~and~\enumref{ass:alphabetagamma}{enum:lowerbound}, $\mu \in \mathrm{AC}(\D)$ and, from Assumption~\enumref{ass:alphabetagamma2}{enum:AC_extra}, $\mathfrak u \in \mathrm{AC}(\D)$ and thus $\mathfrak u \in L^\infty(\D)$ and $\mathfrak u'\in L^1(\D)$ \cite{folland}.
\end{proof}
\begin{figure}[tb]
    \centering
    \includegraphics[height=0.25\textwidth]{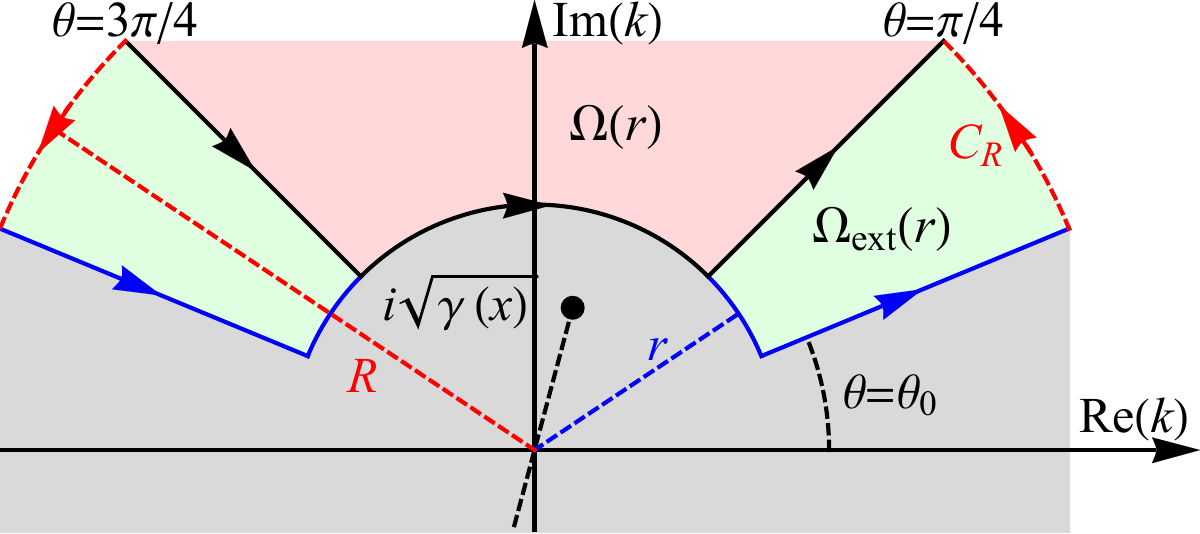}
    \caption{The region $\Omegaext(r) = \Omegaextdef$ described in Lemma~\ref{lem:Re_inu} ($\Omega(r) \,\cup$ the green regions) and the contour $C_R = \{k\in\C:|k|=R$ and $\theta_0 < \arg(k) < \pi/4 $ or $ 3\pi/4 < \theta < \pi-\theta_0 \}$.}
    \label{fig:Omegaext}
\end{figure}
\begin{lemma} \label{lem:Re_inu}
    There exists an $r>\sqrt{M_\gamma}$, $m_{i\mathfrak n}>0$, and $0 < \theta_0<\pi/4$ such that 
    \begin{align} \label{eqn:Re_inu}
        \Re(ik\mathfrak n(k,x)) \leq -m_{i\mathfrak n}|k|,
    \end{align}
    for $k\in\Omegaext(r)$, where $\Omegaext(r) = \Omegaextdef$, see Figure~\ref{fig:Omegaext}.
\end{lemma}
\begin{proof}
    With $\phi = \argi(k \mathfrak n(k,x))$, $\Theta = \sup_{x\in\D} |\argi(\alpha(x)\beta(x))|$, $\psi=\argi(1+\gamma(x)/k^2)$, (and $\theta = \arg(k)$), we have from the definition of $\mathfrak n(k,x)$ in Definition~\ref{def:mu},
    \begin{align} \label{eqn:nuarg}
        \theta - \frac12(\Theta - \psi) \leq \phi \leq  \theta + \frac12(\Theta + \psi),
    \end{align}
    see Figure~\ref{fig:args}.
    \begin{figure}[b]
        \centering
        \begin{subfigure}[tb]{0.4\textwidth}
            \centering
            \includegraphics[width=\textwidth]{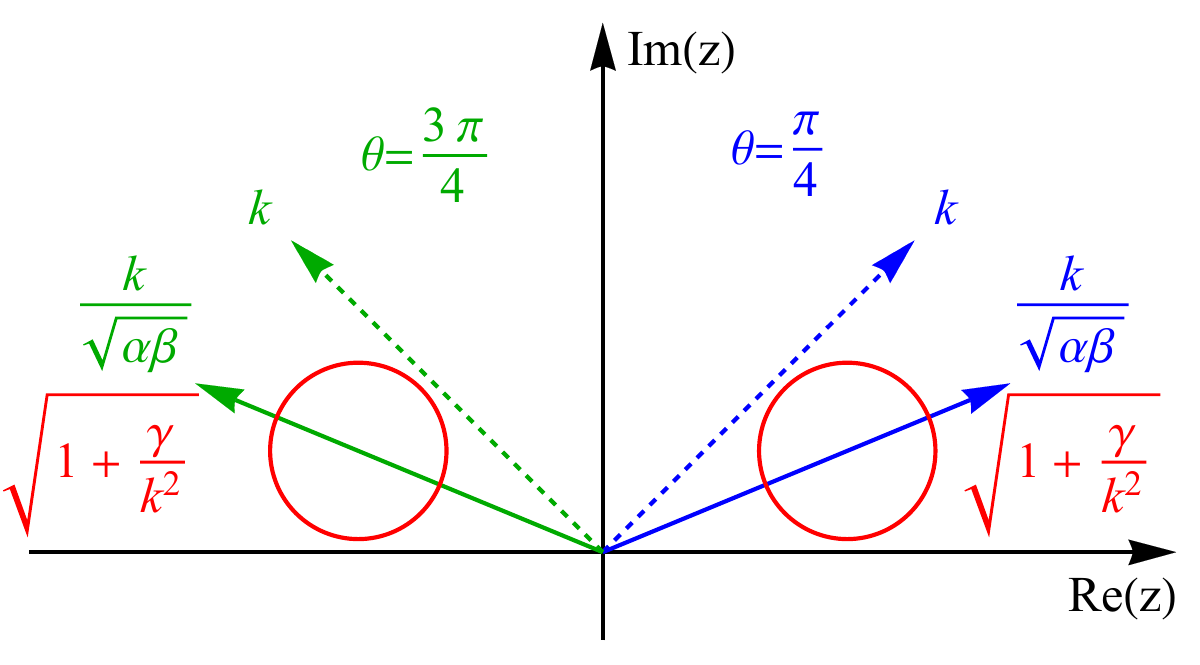}
            \caption{}
            \label{fig:args}
        \end{subfigure}
        \hspace*{0.3in}
        \begin{subfigure}[tb]{0.4\textwidth}
            \centering
            \includegraphics[width=\textwidth]{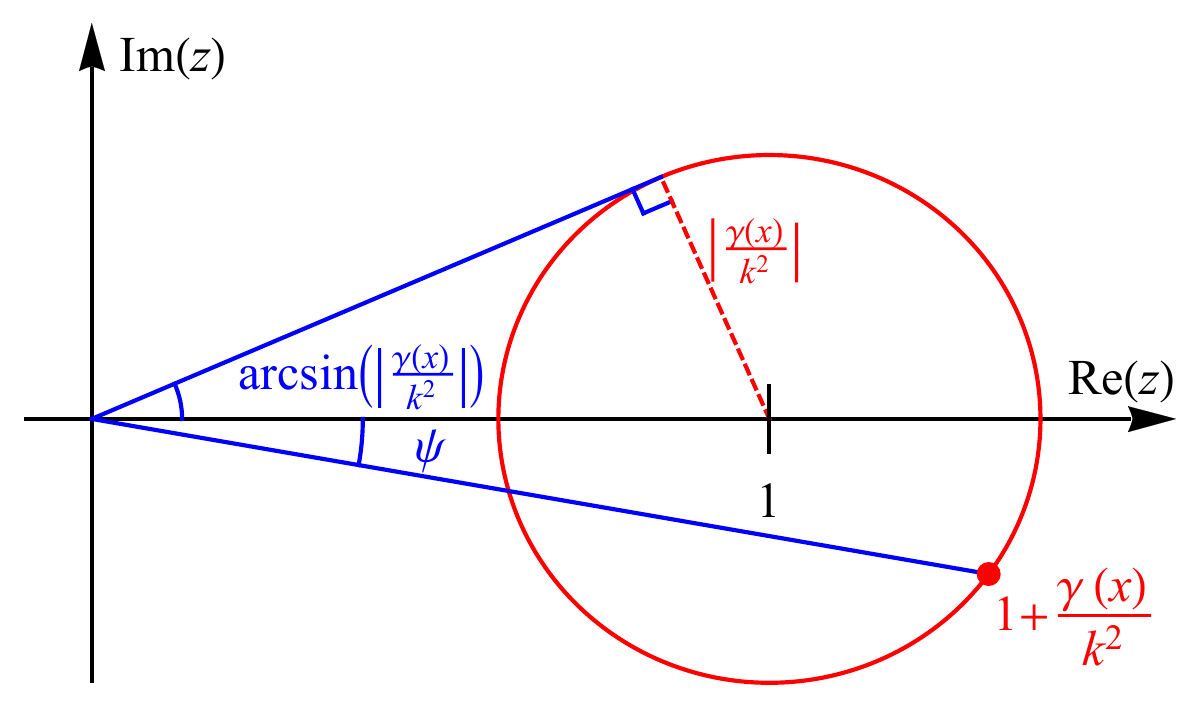}
            \caption{}
            \label{fig:args_circle}
        \end{subfigure}
        \caption{(a) The arguments of $k\mathfrak n(k,x)$ and its components, (b) $\psi = \argi\left(1+\gamma(x)/k^2\right)$.}
        \label{fig:argnu}
    \end{figure}
    Using Assumption~\enumref{ass:alphabetagamma}{enum:bounded} and Figure~\ref{fig:args_circle}, 
    \begin{align} \label{eqn:g_arg_bound}
        |\psi| \leq \arcsin\!\left(\left|\frac{\gamma(x)}{k^2}\right| \right) \leq \arcsin\!\left(\frac{M_\gamma}{r^2}\right) \leq \frac{2M_\gamma}{r^2} .
    \end{align}
    Since $0 \leq \Theta < \pi/2$, we can choose $r>\sqrt{M_\gamma}$ large enough so that 
    \begin{align}
        \left|\psi\right| \leq \frac12\left( \frac\pi2-\Theta\right)\!,
    \end{align}
    which gives, from \eqref{eqn:nuarg},
    \begin{align} \label{eqn:argnubounds}
        \theta - \theta_1  = \theta - \frac14\left( \Theta + \frac\pi2 \right) {{\leq \theta - \frac12\left(\Theta + |\psi|\right)}} \leq \phi \leq  {{\theta + \frac12\left(\Theta + |\psi|\right) \leq}} \theta + \frac14\left( \Theta + \frac\pi2 \right)  = \theta + \theta_1 ,
    \end{align}
    which defines $0 < \theta_1 < \pi/4$. 
    For $|k|>r$ and $\theta_1\leq \argi(k)\leq \pi - \theta_1$, we have that $0\leq \phi\leq \pi$, so that $\Re(ik\mathfrak n(k,x)) \leq 0$. In particular, $\Re(ik \mathfrak n(k,x)) < 0$ for $k\in \Omegaext(r)$. 
    More specifically, using that $\sin(\phi)\geq \phi(\pi - \phi)/\pi$ for $0 \leq \phi \leq \pi$, then for $\theta_1 \leq \theta \leq \pi -\theta_1$, we have
    \begin{align}
        \Re(ik \mathfrak n(k,x)) &= -|k\mathfrak n(k,x)|\sin(\phi) \leq -\frac{1}{\pi}m_{\mathfrak n}|k|\phi(\pi - \phi) \leq  -\frac{1}{\pi}m_{\mathfrak n}|k|(\theta-\theta_1)(\pi - \theta_1 - \theta).
    \end{align}
    Finally, \eqref{eqn:Re_inu} follows from choosing $\theta_0$ such that $0 \leq \theta_1 < \theta_0 < \pi/4$ and letting \mbox{$m_{i\mathfrak n} = m_{\mathfrak n} (\theta_0 - \theta_1)/4$}.
\end{proof}
Having established some properties of the coefficient functions $\alpha,\beta$ and the {\textit{dispersion functions}} $\mathfrak n$, $(\beta \mathfrak n)$, we define a generalization $\Jn$ of the {\textit{accumulation functions}} $\mathcal E_n^{(a,b)}(k)$, $\tilde{\mathcal E}_n^{(a,b)}(k)$, $\mathcal C_n^{(a,b)}(k)$, and $\mathcal S_n^{(a,b)}(k)$, and we show some relations between these functions. Further, we show these functions are bounded and well defined, and we find their large-$k$ asymptotics.
\begin{definition} \label{def:Jndefinitions} 
    \begin{subequations} \label{eqn:Jndefinition}
        With $r$ from Lemma~\ref{lem:Re_inu}, $(a,b)\subseteq \D$ and $k\in \Omegaext(r)$. For integer $n>0$, we define the function
        \begin{align} \label{eqn:Jndefinition_n>0}
            \Jn &= \frac{1}{2^n} \int_{a=y_0<y_1<\cdots<y_n<y_{n+1}=b} \left(\prod_{p=1}^n \fracbetanuargs{y_p} \right) \exp\left( \sum_{p=0}^n \sigma_{p,n} \int_{y_p}^{y_{p+1}} ik\mathfrak n(k,\xi) \, d\xi\right) d\mathbf y_n ,
        \end{align}
        where $\sigma_{p,n}$ is a non-negative integer-valued function of $n$ and $p=0,1,\ldots,n$. Here we require for any $p$ that $\sigma_{p,n}\neq \sigma_{p+1,n}$. For $n<0$, we define $\Jn = 0$, and for $n=0$, we define
        \begin{equation} \label{eqn:J0}
            \mathcal J_0^{(a,b)}[\sigma_{0,0}](k) = \exp\left(\sigma_{0,0} \int_{a}^{b} ik\mathfrak n(k,\xi) \, d\xi\right)\!.
        \end{equation}
        The function $\mathcal J_n^{(a,b)}[\sigma_{p,n}](k)$ is defined as $a\to-\infty$ if $\sigma_{0,n}=0$, and as $b\to\infty$ if $\sigma_{n,n}=0$ $($if $\D$ is unbounded$)$. 
    \end{subequations}
        
    Finally, we define
    \begin{align} \label{eqn:mathscrCnSn}
        \mathscr C_n^{(a,b)}(k) = \etwo{\int_a^bik \mathfrak n(k,\xi)\, d\xi} \mathcal C_n^{(a,b)}(k) 
        \qquad \text{ and } \qquad 
        \mathscr S_n^{(a,b)}(k) = \etwo{\int_a^bik \mathfrak n(k,\xi)\, d\xi}\mathcal S_n^{(a,b)}(k),
    \end{align}
    where $\mathcal C_n^{(a,b)}(k)$ and $\mathcal S_n^{(a,b)}(k)$ are defined in \eqref{eqn:CnSn}.
\end{definition}
\begin{lemma} \label{lem:Enrelations}
    With $\mathcal E_n^{(a,b)}(k)$ and $\tilde{\mathcal E}_n^{(a,b)}(k)$ defined in \eqref{eqn:En&Entilde}, and $\mathscr C_n^{(a,b)}(k)$ and $\mathscr S_n^{(a,b)}(k)$ defined in \eqref{eqn:mathscrCnSn}, we have the following relations:
    \begin{subequations} \label{eqn:relations}
        \begin{align}
            \En &= \mathcal J_n^{(a,b)}[1-(-1)^{n-p}](k), \\
            \tilde{\mathcal E}_n^{(a,b)}(k) &= \Jnp{1-(-1)^p}, \\
            \mathscr C_n^{(a,b)}(k) &= \frac12\left[ {\mathcal J}_n^{(a,b)}[1+(-1)^p](k) + {\mathcal J}_n^{(a,b)}[1-(-1)^p](k)\right]\!, \\
            \mathscr S_n^{(a,b)}(k) &= \frac1{2i}\left[ {\mathcal J}_n^{(a,b)}[1+(-1)^p](k) - {\mathcal J}_n^{(a,b)}[1-(-1)^p](k)\right]\!.
        \end{align}
    \end{subequations}
\end{lemma}
\begin{proof}
    The proofs follow immediately from the definitions in \eqref{eqn:En&Entilde} and \eqref{eqn:CnSn}.
\end{proof}
The next two lemmas give bounds and asymptotics for the function $\Jn$. 
\begin{lemma}\label{lem:Jnbounds}
    For $(x,y)\subseteq (a,b)\subseteq \D$, $k\in \Omegaext(r)$, and $r$ from Lemma~\ref{lem:Re_inu}, 
    \begin{align} \label{eqn:Jbound}
        \left| \Jn \right| \leq \frac{1}{2^nn!} \inorm{ \fracbetanu }{(a,b)}^n 
        ~~~ \text{ and } ~~~
        \left| \sum_{\ell=0}^n (-1)^{\lambda\ell} \mathcal J_{n-\ell}^{(a, x)}[\sigma_{p,n-\ell}] (k) \mathcal J_\ell^{(y,b)}[\tilde \sigma_{p,\ell}](k) \right| \leq \frac{1}{2^nn!} \inorm{ \fracbetanu }{(a,b)}^n\!\!,
    \end{align}
    \no where $\lambda = 0,1$. These inequalities hold as $a\to-\infty$ and $b\to\infty$, {provided the functions are defined}. Thus, $\mathcal J_n^{(a,b)}[\sigma_{p,n}](k)$ is well defined. The same bounds hold for $\mathcal E_n^{(a,b)}(k)$, $\tilde{\mathcal E}_n^{(a,b)}(k)$, ${\mathscr C}_n^{(a,b)}(k)$, and ${\mathscr S}_n^{(a,b)}(k)$.
\end{lemma}
\begin{proof}
    \begin{subequations}
        By Lemma~\ref{lem:Re_inu}, the exponentials in \eqref{eqn:Jndefinition} are bounded by $1$ for $k\in \Omegaext(r)$. Using Lemma~\ref{lem:beta_nu},
        \begin{align} 
            \left| \Jn \right| &\leq \frac{1}{2^{n}} \int_{a < y_1 <\cdots <y_n<b} \, \left| \prod_{p=1}^{n} \fracbetanuargs{y_p} \right| \, dy_1 \cdots dy_n = \frac{1}{2^nn!} \inorm{\fracbetanu}{(a,b)}^n \!\!,
        \end{align}
        \no so that
        \begin{align}\nonumber
            \left| \sum_{\ell=0}^n (-1)^{\lambda \ell} \mathcal J_{n-\ell}^{(a, x)}[\sigma_{p,n-\ell}](k) \mathcal J_\ell^{(y,b)}[\sigma_{p,\ell}](k) \right| 
            &\leq \sum_{\ell=0}^n \frac{1}{2^n(n-\ell)!\ell!} \inorm{\fracbetanu}{(a,x)}^{n-\ell} \inorm{\fracbetanu}{(y,b)}^{\ell}   \\ 
            &{=  \frac{1}{2^n n!} \left( \inorm{\fracbetanu}{(a,x)} +  \inorm{\fracbetanu}{(y,b)}\right)^n} \leq \frac{1}{2^nn!} \inorm{ \fracbetanu }{(a,b)}^n\!\!.
        \end{align}
        For $\mathcal E_n^{(a,b)}(k)$, $\tilde{\mathcal E}_n^{(a,b)}(k)$, ${\mathscr C}_n^{(a,b)}(k)$, and ${\mathscr S}_n^{(a,b)}(k)$, the result follows from \eqref{eqn:relations}. These bounds hold as $a\to-\infty$ or $b\to \infty$, {provided the functions are defined.}
    \end{subequations}
\end{proof}
\begin{lemma}\label{lem:Jnasymptotics} 
    \begin{subequations}
        There exists $r>\sqrt{M_\gamma}$,  such that  for any $(a,b)\subseteq \D$ and $n\geq 1$,
        \begin{align} \label{eqn:Jnlimit}
            \Jn \to 0, \quad \text{ as } |k|\to\infty,~~~k\in\Omegaext(r).
        \end{align}
        This result holds as $a\to-\infty$ and $b\to\infty$, {provided the functions are defined}. The result extends to $\mathcal E_n^{(a,b)}(k)$, $\tilde{\mathcal E}_n^{(a,b)}(k)$, ${\mathscr C}_n^{(a,b)}(k)$, and ${\mathscr S}_n^{(a,b)}(k)$.
    
        Next, we define $\lambda_{p,n} = \sigma_{p-1,n} - \sigma_{p,n}$. Using Assumption~\ref{ass:alphabetagamma2} and since $\lambda_{p,n} \neq 0$ (see Definition~\ref{def:Jndefinitions}),  we have
        \begin{align} \label{eqn:J1_asymptotics}
            \mathcal J_1^{(a,b)}[\sigma_{p,n}](k) &= \frac{1}{4\lambda_{1,1} ik} \left[ \mathfrak u(b) \exp\left( \sigma_{0,1}\int_{a}^{b} ik\mu(\xi) \, d\xi\right) - \mathfrak u(a) \exp\left(  \sigma_{1,1} \int_{a}^{b} ik\mu(\xi) \, d\xi\right) \right] +\littleoh(k^{-1}).
        \end{align}
        There exists $r>\sqrt{M_\gamma}$ and $C>1$ such that, for any $(a,b)\subseteq \D$ and $k\in\Omegaext(r)$, 
        \beq \label{eqn:Jnasymptotics}
            \big|{\mathcal J}_{n}^{(a,b)}[\sigma_{p,n}](k)\big| \leq \frac{C^{n}}{|k|^{\lfloor \frac{n+1}{2} \rfloor}},
        \eeq
        where $\lfloor \, \cdot \, \rfloor$ is the floor function.
    \end{subequations}
\end{lemma}
\begin{proof}
    From Lemma~\ref{lem:Re_inu}, for any $r>\sqrt{M_\gamma}$ and for all $k\in \Omegaext(r)$, we have
    \begin{align} \nonumber
        \left|\mathcal J_n^{(a,b)}[\sigma_{p,n}](k)\right| 
        &{{= \frac{1}{2^n} \left| \int_{a=y_0< y_1 < \cdots < y_n < y_{n+1} = b} \left(\prod_{p=1}^n \fracbetanuargs{y_p} \right) \exp\left(\sum_{p=0}^n \sigma_{p,n} \int_{y_p}^{y_{p+1}} ik\mathfrak n(k,\xi)\,d\xi\right)  dy_1\cdots dy_n \right|}} \\
        &\leq \frac{1}{2^n} \int_{a=y_0< \cdots < y_{n+1} = b} \left(\prod_{p=1}^n \left| \fracbetanuargs{y_p} \right| \right) \etwo{- m_{i\mathfrak n} |k| \sum_{p=0}^n \sigma_{p,n} (y_{p+1} - y_p)}  dy_1 \cdots dy_n,
    \end{align}
    and since $\sigma_{p,n}\geq 0$ and $\sigma_{p,n}\neq \sigma_{p+1,n}$ for any $p$, the argument of the exponential is strictly negative. Thus, by Lemma~\ref{lem:Jnbounds} and the Dominated Convergence Theorem (DCT), we have \eqref{eqn:Jnlimit}. For $\mathcal E_n^{(a,b)}(k)$, $\tilde{\mathcal E}_n^{(a,b)}(k)$, ${\mathscr C}_n^{(a,b)}(k)$, and ${\mathscr S}_n^{(a,b)}(k)$, the result follows from \eqref{eqn:relations}.

    Using the definition of $\mathfrak n(k,x)$ in Definition~\ref{def:mu}, it is straightforward to show that
    \begin{align} \label{eqn:E_approx} 
        \exp\left( ik \int_y^x \mathfrak n(k,\xi) \, d\xi \right) 
        %
        &= \etwo{ik\int_y^x \mu(\xi) \, d\xi} \big(1 + \bigoh(|x-y| k^{-1}) \big),
    \end{align}
    for which, on a finite-interval, we may omit $|x-y|$.
    Using Assumption~\enumref{ass:alphabetagamma}{enum:L1}, \eqref{eqn:mathfrak_u_def}, \eqref{eqn:betanu_expansion}, and \eqref{eqn:E_approx} in \eqref{eqn:Jndefinition} for $n = 1$, 
    we have
    \begin{align} 
        \mathcal J_1^{(a,b)}[\sigma_{p,n}](k) &= \frac{1 + \bigoh(k^{-1})}{4} \int_{a}^{b} \mathfrak u(y) \mu(y) \exp\left( \sigma_{0,1} \int_{a}^{y} + \sigma_{1,1} \int_{y}^{b} ik\mu(\xi) \, d\xi\right) d y + \bigoh (k^{-2}).
    \end{align}
    By Lemma~\ref{lem:beta_nu}, $\mathfrak u\in \mathrm{AC}(\D)$ and integration by parts gives
    \begin{align} 
        \mathcal J_1^{(a,b)}[\sigma_{p,n}](k) 
        \nonumber
        &= \frac{1}{4\lambda_{1,1}ik} \left( \mathfrak u(b) \exp\left( \sigma_{0,1} \int_{a}^{b}  ik\mu(\xi) \, d\xi\right) - \mathfrak u(a) \exp\left(  \sigma_{1,1} \int_{a}^{b} ik\mu(\xi) \, d\xi\right) \right) \\
        &~~~- \frac{1}{4\lambda_{1,1}ik} \int_{a}^{b} \mathfrak u'(y) \exp\left( \sigma_{0,1} \int_{a}^{y} + \sigma_{1,1} \int_{y}^{b} ik\mu(\xi) \, d\xi\right)  d y + \bigoh(k^{-2}).
    \end{align}
    By Lemma~\ref{lem:beta_nu} and the DCT, we obtain \eqref{eqn:J1_asymptotics}.
    
    Inequality \eqref{eqn:Jnasymptotics} for $n=0$ and $n=1$ follows from \eqref{eqn:J0} and \eqref{eqn:J1_asymptotics}, respectively. 
    Using \eqref{eqn:mathfrak_u_def}, \eqref{eqn:betanu_expansion}, and \eqref{eqn:E_approx} in \eqref{eqn:Jndefinition} for $n\geq 2$, 
    we have
    \begin{align} \label{eqn:Jn_intermediate}
        \nonumber
        \Jn 
        &= \frac{1+ \bigoh (k^{-1})}{2^{n+1}} \int_{a< \cdots <b} \left(\prod_{p=1}^{n-1} \fracbetanuargs{y_p} \right) \mathfrak u(y_n) \mu(y_n) \exp\left( \sum_{p=0}^n \sigma_{p,n} \int_{y_p}^{y_{p+1}} ik \mu(\xi) \, d\xi\right)  d \mathbf y_n \\
        &~~~+ \frac{1+ \bigoh (k^{-1})}{2^{n+1}} \int_{a<\cdots<b} \left(\prod_{p=1}^{n-1} \fracbetanuargs{y_p} \right) \left( \frac{\gamma'(y_{n})}{k^2+\gamma(y_{n})} \right) \exp\left( \sum_{p=0}^n \sigma_{p,n} \int_{y_p}^{y_{p+1}} ik\mu(\xi) \, d\xi\right) d \mathbf y_n.
    \end{align}
    Let $\mathcal I_n^{(a,b)}(k)$ denote the integral in the first line of \eqref{eqn:Jn_intermediate}:
    \begin{align} 
        \mathcal I_n^{(a,b)}(k)
        &= \int_{a< \cdots <b} \left(\prod_{p=1}^{n-1} \fracbetanuargs{y_p} \right) \mathfrak u(y_n) \mu(y_n) \exp\left( \sum_{p=0}^n \sigma_{p,n} \int_{y_p}^{y_{p+1}} ik \mu(\xi) \, d\xi\right) d \mathbf y_n.
    \end{align}
    Integration by parts with respect to $y_n\in (y_{n-1},b)$ gives
    \begin{align}
        \mathcal I_n^{(a,b)}(k)
        \nonumber
        &=  \frac{\mathfrak u(b)}{ik\lambda_{n,n}} \int_{a=y_0<\cdots <y_n=b} \left(\prod_{p=1}^{n-1} \fracbetanuargs{y_p} \right)   \exp\left( \sum_{p=0}^{n-1} \sigma_{p,n} \int_{y_p}^{y_{p+1}} ik \mu(\xi) \, d\xi \right) d \mathbf y_{n-1} \\
        \nonumber
        &~~~-  \int_{a<\cdots<b}  \left(\prod_{p=1}^{n-1} \fracbetanuargs{y_p} \right) \frac{\mathfrak u(y_{n-1})}{ik\lambda_{n,n}} \etwo{ \sigma_{n,n} \int_{y_{n-1}}^{b} ik \mu(\xi) \, d\xi +  \sum_{p=0}^{n-2} \sigma_{p,n} \int_{y_p}^{y_{p+1}} ik \mu(\xi) \, d\xi } d \mathbf y_{n-1} \\\la{eqn:Jnn}
        &~~~-  \int_{a<\cdots<b} \left(\prod_{p=1}^{n-1} \fracbetanuargs{y_p} \right) \frac{\mathfrak u'(y_n)}{ik\lambda_{n,n}} \etwo{ \sum_{p=0}^{n-1} \sigma_{p,n} \int_{y_p}^{y_{p+1}} ik \mu(\xi) \, d\xi } d \mathbf y_{n}.
    \end{align}
    In the second line of \eqref{eqn:Jnn}, we integrate over $y_{n-1}\in (a,b)$ last, leaving the remaining integral over $a=y_0<y_1<\cdots<y_{n-2} < y_{n-1}$ to be done first. Similarly, in the third line of \eqref{eqn:Jnn}, we integrate over $y_n\in (a,b)$ last and leave the remaining integral over $a=y_0<y_1<\cdots<y_{n-1}<y_n$ to be done first. Returning to \eqref{eqn:Jn_intermediate} yields
    \begin{align} 
        \Jn 
        \nonumber
        &= \frac{1+ \bigoh (k^{-1})}{2^{n+1}} \frac{\mathfrak u(b)}{\lambda_{n,n}ik} \int_{a=y_0<\cdots <y_n=b} \left(\prod_{p=1}^{n-1} \fracbetanuargs{y_p} \right)   \exp\left( \sum_{p=0}^{n-1} \sigma_{p,n} \int_{y_p}^{y_{p+1}} ik \mu(\xi) \, d\xi\right)  d \mathbf y_{n-1} \\
        \nonumber
        &~~~- \frac{1+ \bigoh(k^{-1})}{2^{n+1}} \int_a^b dy_{n-1} \, \frac{\mathfrak u(y_{n-1})}{\lambda_{n,n}ik} \fracbetanuargs{y_{n-1}} \exp\left( \sigma_{n,n} \int_{y_{n-1}}^{b} ik \mu(\xi) \, d\xi\right) \times \\
        \nonumber
        &\hspace*{1in}\times \int_{a=y_0< \cdots <y_{n-2}<y_{n-1}} \left(\prod_{p=1}^{n-2} \fracbetanuargs{y_p} \right)   \exp\left( \sum_{p=0}^{n-2} \sigma_{p,n} \int_{y_p}^{y_{p+1}} ik \mu(\xi) \, d\xi \right) d \mathbf y_{n-2} \\
        \nonumber
        &~~~- \frac{1+ \bigoh (k^{-1})}{2^{n+1}} \int_a^b dy_n \left(\frac{\mathfrak u'(y_n)}{\lambda_{n,n}ik} + \frac{\gamma'(y_{n})}{k^2+\gamma(y_{n})} \right) \exp\left( \sigma_{n,n} \int_{y_n}^{b} ik \mu(\xi) \, d\xi\right) \times \\
        &\hspace{1in}\times \int_{a=y_0< \cdots <y_{n-1}<y_n} \left(\prod_{p=1}^{n-1} \fracbetanuargs{y_p} \right)   \exp\left( \sum_{p=0}^{n-1} \sigma_{p,n} \int_{y_p}^{y_{p+1}} ik \mu(\xi) \, d\xi\right) d \mathbf y_{n-1},
    \end{align}
    which gives the asymptotic recurrence relation
    \begin{align} 
        \Jn 
        \nonumber
        &= \frac{1+ \bigoh (k^{-1})}{4} \frac{\mathfrak u(b)}{\lambda_{n,n}ik} \mathcal J_{n-1}^{(a,b)}[\sigma_{p,n}](k) \\
        \nonumber
        &~~~- \frac{1+ \bigoh (k^{-1})}{8} \int_a^b \frac{\mathfrak u(y_{n-1})}{\lambda_{n,n}ik} \fracbetanuargs{y_{n-1}} \exp\left( \sigma_{n,n} \int_{y_{n-1}}^{b} ik \mu(\xi) \, d\xi\right) \mathcal J_{n-2}^{(a,y_{n-1})}[\sigma_{p,n}](k) \, dy_{n-1} \\
        &~~~- \frac{1+ \bigoh (k^{-1})}{4} \int_a^b \left(\frac{\mathfrak u'(y_n)}{\lambda_{n,n}ik} + \frac{\gamma'(y_{n})}{k^2+\gamma(y_{n})} \right) \exp\left( \sigma_{n,n} \int_{y_n}^{b} ik \mu(\xi) \, d\xi\right) \mathcal J_{n-1}^{(a,y_n)}[\sigma_{p,n}](k) \, dy_n.
    \end{align}
    Assuming \eqref{eqn:Jnasymptotics} holds for $n=0,1,\ldots,m-1$, and using that $|\lambda_{p,n}|\geq 1$, we find
    \begin{align} 
        \left| \Jn \right| 
        &\leq \frac{1+ \bigoh (k^{-1})}{4} \left[ \frac{\|\mathfrak u\|_\infty}{|k|} \frac{C^{n-1}}{|k|^{\lfloor \frac{n}{2} \rfloor}} + \frac{\|\mathfrak u\|_\infty}{2|k|} \left\| \fracbetanu \right\|_\D \frac{C^{n-2}}{|k|^{\lfloor \frac{n-1}{2}\rfloor}} +  \left(\frac{ \|\mathfrak u'\|_\D }{|k|} + \frac{\|\gamma'\|_\D}{|k|^2-M_\gamma} \right) \frac{C^{n-1}}{|k|^{\lfloor \frac{n}{2} \rfloor}} \right]\!,
    \end{align}
    which, using Lemma~\ref{lem:beta_nu}, gives \eqref{eqn:Jnasymptotics} for $n\geq 0$.
\end{proof}
Having defined the function $\Jn$ and established some of its properties, we prove that the {function} $\Delta(k)$ is bounded and well defined, and that the {``transforms''} $\Phi_0(k,x)$, $\Phi_{f}(k,x,t)$ and $\mathcal B_m(k,x)$, of the initial condition $q_0(x)$, the inhomogeneous function $f(x,t)$, and the boundary functions $f_m(t)$, respectively, are bounded and well defined.
\begin{definition} \label{def:ftilde,Phi,q}
    We define 
    \begin{subequations} \label{eqn:Phi's}
        \begin{align} 
            \label{eqn:Phi_0}
            \Phi_{0}(k,x) &= \int_{\D} \frac{\Psi(k,x,y)q_\alpha(y)}{ \sqrt{(\beta\mathfrak n)(k,x)}\sqrt{(\beta\mathfrak n)(k,y)}} \, dy, \\ 
            \label{eqn:Phi_f}
            \Phi_{f}(k,x,t) &=\int_{\D} \frac{\Psi(k,x,y)\tilde f_\alpha(k^2,y,t)} {\sqrt{(\beta\mathfrak n)(k,x)} \sqrt{(\beta\mathfrak n)(k,y)}} \, dy,
        \end{align}
    \end{subequations}
    so that $\Phi_\psi(k,x,t) = \Phi_0(k,x) + \Phi_f(k,x,t)$.
    We define the corresponding parts of the solution as
    \begin{subequations} \label{eqn:q's}
        \begin{align}
            \label{eqn:q_0}
            q_0(x,t) &= \frac{1}{2\pi} \int_{\partial\Omega(r)} \frac{\Phi_0(k,x)}{\Delta(k)} e^{-k^2t} \, dk, \\ \label{eqn:q_f}
            q_f(x,t) &= \frac{1}{2\pi} \int_{\partial\Omega(r)} \frac{\Phi_f(k,x,t)}{\Delta(k)} e^{-k^2t} \, dk, \\ \label{eqn:q_Bm}
            q_{\mathcal B_m}(x,t) &= \frac{1}{2\pi} \int_{\partial\Omega(r)} \frac{\mathcal B_m(k,x)}{\Delta(k)} F_m(k^2,t) e^{-k^2t} \, dk, \qquad m=0,1,
        \end{align}
    \end{subequations}
    where we define $\mathcal B_m(k,x) = 0$ $(m=0,1)$ for the whole-line problem and $\mathcal B_1(k,x) = 0$ for the half-line problem. Thus $q(x,t) = q_0(x,t) + q_f(x,t) + q_{\mathcal B_0}(x,t) + q_{\mathcal B_1}(x,t)$ for the finite-interval, half-line and whole-line problems.
\end{definition}
\begin{lemma} \label{lem:Delta_asymptotics} 
    For all three problems, there exists $r>\sqrt{M_\gamma}$ and $M_\Delta>0$, so that for all $k\in \Omegaext(r)$, 
    \begin{subequations} \label{eqn:Delta_asym+b's}
        \begin{align} \label{eqn:Delta_asym}
            \Delta(k) &= \mathfrak b_0(k) (1+\varepsilon(k)) \qquad \text{ and } \qquad \frac12 |\mathfrak b_0(k)| \leq |\Delta(k)| \leq M_\Delta,
        \end{align}
        where $|\varepsilon(k)|<1/2$. For the whole-line problem,
        \begin{align} \label{eqn:b0_wl}
            \mathfrak b_0(k) = 1;
        \end{align}
        for the half-line problem, 
        \begin{align} \label{eqn:b0_hl}
            \mathfrak b_0(k) = 2\left( \frac{ia_0}{k \mathfrak n(k,x_l)} - a_1 \right)\!;
        \end{align}
        and, for the finite-interval problem,
        \begin{enumerate}
            \item if $(a:b)_{2,4}\neq 0$, then
            \begin{align} \label{eqn:b0_1} 
                \mathfrak b_0(k) = -  (a:b)_{2,4};
            \end{align}
            \item if $(a:b)_{2,4} = 0$ and $m_{\mathfrak c_0} \neq 0$, then
            \begin{align} \label{eqn:b0_2} 
                \mathfrak b_0(k) = \frac{im_{\mathfrak c_0}}{k};
            \end{align}
            \item if $(a:b)_{2,4} = 0$, $m_{\mathfrak c_0} = 0$, $m_{\mathfrak c_1} = 0$, and $(a:b)_{1,3} \neq 0$, then
            \begin{align} \label{eqn:b0_3a} 
                \mathfrak b_0(k) = -\frac{m_{\mathfrak s}}{k^2};
            \end{align}
            \item with Assumption~\ref{ass:alphabetagamma2}, if $(a:b)_{2,4} = 0$, $m_{\mathfrak c_0} = 0$, $m_{\mathfrak c_1} \neq 0$, and $m_{\mathfrak c_1} \mathfrak u_{+} - 8m_{\mathfrak s} \neq 0$,
            then
            \begin{align} \label{eqn:b0_3b} 
                \mathfrak b_0(k) = \frac{1}{8k^2} \left( m_{\mathfrak c_1} \mathfrak u_{+} - 8m_{\mathfrak s}\right)\!.
            \end{align}
        \end{enumerate}
    \end{subequations}
\end{lemma}
\begin{proof}
    For the whole-line problem, 
    \begin{align} 
        \Delta(k) = 1 + \sum_{n=1}^\infty \mathcal E_{2n}^{(-\infty,\infty)}(k)  = 1 + \varepsilon(k).
    \end{align}
    By Lemmas~\ref{lem:Jnbounds}~and~\ref{lem:Jnasymptotics} and the     DCT, 
    \begin{align} 
        \varepsilon(k) = \sum_{n=1}^\infty \mathcal E_{2n}^{(-\infty,\infty)}(k)  \to 0,
    \end{align}
    as $|k|\to \infty$. Thus, we can choose $r$ sufficiently large so that for $k\in \Omegaext(r)$, $|\varepsilon(k)|<1/2$,
    \no and
    \begin{align}
        \frac 12 \leq 1 - |\varepsilon(k)| \leq |\Delta(k)| \leq 1 + |\varepsilon(k)| < \frac 32.
    \end{align}

    For the half-line problem, we write
    \begin{align}
        \Delta(k) &= 2\sum_{n=0}^\infty \left(\frac{(-1)^nia_0}{k\mathfrak n(k,x_l)} - a_1\right) \mathcal E_n^{(x_l,\infty)}(k) 
        =2\left( \frac{ia_0}{k \mathfrak n(k,x_l)} - a_1\right)\left[1 + \sum_{n=1}^\infty \frac{\frac{(-1)^nia_0}{k \mathfrak n(k,x_l)} - a_1}{\frac{i a_0}{k \mathfrak n(k,x_l)} - a_1} \mathcal E_n^{(x_l,\infty)}(k) \right] \! \!.
    \end{align}
    Recall that we require $(a_0,a_1)\neq (0,0)$. If $a_1\neq 0$, we choose $r>\sqrt{M_\gamma}$ sufficiently large so that $|a_1|>|a_0|/(m_{\mathfrak n} r)$ for $k\in \Omegaext(r)$. Using this, for either case, $a_0\neq 0$ or $a_1\neq 0$, we have 
    \begin{align} \label{eqn:A_hl}
        \left| \frac{\frac{(-1)^nia_0}{k \mathfrak n(k,x_l)} - a_1}{\frac{ia_0}{k \mathfrak n(k,x_l)} - a_1} \right| \leq \frac{|a_1|+\frac{|a_0|}{m_{\mathfrak n} r}}{\left||a_1|- \frac{|a_0|}{m_{\mathfrak n} r}\right|} = A < \infty, 
    \end{align}
    which defines $A\geq1$. We have
    \begin{align}
        |\varepsilon(k)| = \left| \sum_{n=1}^\infty \frac{\frac{(-1)^na_0}{k \mathfrak n(k,x_l)} + ia_1}{\frac{a_0}{k \mathfrak n(k,x_l)} + ia_1} \mathcal E_n^{(x_l,\infty)}(k)\right| \leq A \sum_{n=1}^\infty \left|\mathcal E_n^{(x_l,\infty)}(k)\right| \to 0, 
    \end{align}
    by the DCT. We choose $r>\sqrt{M_\gamma}$ large enough such that $|\varepsilon(k)| < 1/2$ for $k\in\Omegaext(r)$. Then
    \begin{align}
        0 < \left| |a_1| - \left|\frac{a_0}{k \mathfrak n(k,x_l)} \right| \right| \leq |\Delta(k)| \leq 3 \left( |a_1| + \frac{|a_0|}{m_{\mathfrak n} r} \right) < \infty.
    \end{align}

    For the finite-interval problem, since
    \begin{align}
        \mathscr C_0^{(x_l,x_r)}(k) = \Xi(k) \mathcal C_0^{(x_l,x_r)}(k) = \frac{1}{2}\left(\Xi(k)^2 + 1\right) \qquad \text{ and } \qquad  \mathscr S_0^{(x_l,x_r)}(k) = \Xi(k)\mathcal S_0^{(x_l,x_r)}(k) = \frac{1}{2i}\left(\Xi(k)^2 - 1\right)\!,
    \end{align}
    where $\Xi(k)$ is defined in \eqref{eqn:Xi_def} and $\mathscr C_n^{(a,b)}(k)$ and $\mathscr S_n^{(a,b)}(k)$ are defined in \eqref{eqn:mathscrCnSn}, we factor out the $n=0$ term in \eqref{eqn:Delta_fi} and write
    \begin{align} 
        \Delta(k) &= 2i \mathfrak a(k) \Xi(k) + i \mathfrak c_0(k) \left(\Xi(k)^2 + 1\right)  +  \mathfrak s_0(k) \left(\Xi(k)^2 - 1\right) + 2i \sum_{n=1}^\infty \left( \mathfrak c_n(k)  \mathscr C_n^{(x_l,x_r)}(k)  + \mathfrak s_n(k) \mathscr S_n^{(x_l,x_r)}(k) \right)\!.
    \end{align}
    Since $\Xi(k) \to 0$ exponentially fast, we have
    \begin{align} 
        \Delta(k) &=  i \mathfrak c_0(k) - \mathfrak s_0(k)  + 2i \sum_{n=1}^\infty \left( \mathfrak c_n(k)\mathscr C_n^{(x_l,x_r)}(k)  + \mathfrak s_n(k) \mathscr S_n^{(x_l,x_r)}(k) \right) + \littleoh( k^{-2}) .
    \end{align}
    \begin{enumerate}
        \item If $(a:b)_{2,4} \neq 0$, then 
        we can write \eqref{eqn:Delta_asym} with $\mathfrak b_0(k)$ defined in \eqref{eqn:b0_1}, where
        \begin{align} 
            \varepsilon (k) 
            &= \frac{-1}{(a:b)_{2,4}}\left[ i\mathfrak c_0(k) - 
            \mathfrak s_0(k) + (a:b)_{2,4} + 2i  \sum_{n=1}^\infty \left( \mathfrak c_n(k)\mathscr C_n^{(x_l,x_r)}(k) + \mathfrak s_n(k) \mathscr S_n^{(x_l,x_r)}(k) \right) \right] + \littleoh (k^{-2} ).
        \end{align}
        Since $\mathfrak c_n(k) = \bigoh(k^{-1})$, $\mathfrak s_0(k) = (a:b)_{2,4} + \bigoh(k^{-2})$, $\mathfrak s_n(k) = \bigoh(k^0)$, and because both $\mathscr C_n^{(x_l,x_r)}(k) \to 0$  and \mbox{$\mathscr S_n^{(x_l,x_r)}(k)\to 0$}  by Lemma~\ref{lem:Jnasymptotics} and both are bounded (see Lemma~\ref{lem:Jnbounds}), we can choose $r>\sqrt{M_\gamma}$ sufficiently large so that $|\varepsilon(k)| < 1/2$ for $k\in\Omegaext(r)$, by the DCT. We have 
        \begin{align}
            \frac{1}{2}|(a:b)_{2,4}| \leq |\Delta(k)| \leq \frac{3}{2}|(a:b)_{2,4}|.
        \end{align}
        \item If $(a:b)_{2,4} = 0$ and $m_{\mathfrak c_0} \neq 0$, then we can write \eqref{eqn:Delta_asym} with $\mathfrak b_0(k)$ defined in \eqref{eqn:b0_2}, where
        \begin{align} 
            \varepsilon (k) &= \frac{k}{im_{\mathfrak c_0}}  \left[ i \mathfrak c_0(k) - \frac{im_{\mathfrak c_0}}{k} - \mathfrak s_0(k) + 2i \sum_{n=1}^\infty \left( \mathfrak c_n(k) \mathscr C_n^{(x_l,x_r)}(k)  + \mathfrak s_n(k) \mathscr S_n^{(x_l,x_r)}(k)\right) \right] + \littleoh (k^{-1} ).
        \end{align}
        Since $\mathfrak c_0(k) = m_{\mathfrak c_0}/k + \bigoh(k^{-3})$, $\mathfrak s_n(k) = \bigoh(k^{-2})$, $\mathfrak c_n(k) = \bigoh(k^{-1})$, and since $\mathscr C_n^{(x_l,x_r)}(k)\to 0$ and $\mathscr S_n^{(x_l,x_r)}(k)\to 0$ and both are bounded (see Lemma~\ref{lem:Jnbounds}), we can choose $r>\sqrt{M_\gamma}$ large enough such that $|\varepsilon(k)| < 1/2$ for $k\in\Omegaext(r)$, by the DCT. We have
        \begin{align}
            \frac{|m_{\mathfrak c_0}|}{2|k|} \leq |\Delta(k)| \leq \frac{3|m_{\mathfrak c_0}|}{2r}. 
        \end{align}
        \item If $(a:b)_{2,4} = 0$, $m_{\mathfrak c_0} = 0$, $m_{\mathfrak c_1} = 0$, and $(a:b)_{1,3} \neq 0$, then $\mathfrak c_0(k) = \mathfrak c_1(k) = 0$ and we can write \eqref{eqn:Delta_asym} with $\mathfrak b_0(k)$ defined in \eqref{eqn:b0_3a}, where
       \begin{align} 
            \varepsilon (k) 
            &= -\frac{k^2}{m_{\mathfrak s}}\left[ -\mathfrak s_0(k) + \frac{m_{\mathfrak s}}{k^2} + 2i \sum_{n=1}^\infty \mathfrak s_n(k) \mathscr S_n^{(x_l,x_r)}(k)\right] + \littleoh (k^{0}).
        \end{align}
        Since $\mathfrak s_n(k) = m_{\mathfrak s}/k^2 + \bigoh(k^{-4})$, and since $\mathscr S_n^{(x_l,x_r)}(k) \to 0$ and is bounded (see Lemma~\ref{lem:Jnbounds}), we can choose $r>\sqrt{M_\gamma}$ sufficiently large so that $|\varepsilon(k)| < 1/2$ for $k\in\Omegaext(r)$, by the DCT. We have
        \begin{align}
            \frac{|m_{\mathfrak s}|}{2|k|^2} \leq |\Delta(k)| \leq \frac{3|m_{\mathfrak s}|}{2r^2}.
        \end{align}
        \item If $(a:b)_{2,4} = 0$, $m_{\mathfrak c_0} = 0$, $m_{\mathfrak c_1} \neq 0 $, and $m_{\mathfrak c_1} \mathfrak u_{+} - 8m_{\mathfrak s} \neq 0$, we can write \eqref{eqn:Delta_asym} with $\mathfrak b_0(k)$ defined in \eqref{eqn:b0_3b}, where
        \begin{align} \nonumber
            \varepsilon (k) 
            &= \frac{8k^2}{m_{\mathfrak c_1} \mathfrak u_{+} - 8m_{\mathfrak s}}\left[ i \mathfrak c_0(k) +2i \mathfrak c_0(k) \mathscr C_2^{(x_l,x_r)}(k) - \mathfrak s_0(k) + 2i \mathfrak c_1(k) \mathscr C_1^{(x_l,x_r)}(k) - \frac{m_{\mathfrak c_1} \mathfrak u_{+} - 8m_{\mathfrak s}}{8k^2} \right. \\
            &~~~ \left. \hspace*{1.2in} + 2i \sum_{n=3}^\infty  \mathfrak c_n(k) \mathscr C_n^{(x_l,x_r)}(k)  + 2i \Xi(k) \sum_{n=1}^\infty \mathfrak s_n(k) \mathscr S_n^{(x_l,x_r)}(k)  \right] + \littleoh(k^0).
        \end{align}
        By Lemmas~\ref{lem:Enrelations}~and~\ref{lem:Jnasymptotics}, we have
        \begin{align} \nonumber
            \mathscr C_1^{(x_l,x_r)}(k) 
            &{ =\frac12\left[ \mathcal J_1^{(x_l,x_r)}[1+(-1)^p](k) + \mathcal J_1^{(x_l,x_r)}[1-(-1)^p](k) \right] } \\ 
            &=-\frac{1}{16ik} \left(\mathfrak u(x_r) + \mathfrak u(x_l)\right)\left(1- \exp\left( \int_{x_l}^{x_r} 2ik \mu(\xi) \, d\xi\right)\right)  +\littleoh(k^{-1})
            =-\frac{1}{16ik}  \mathfrak u_{+}  +\littleoh(k^{-1}).
        \end{align}
        For $n>2$, from Lemma~\ref{lem:Jnasymptotics}, we have
        \begin{align}
            \sum_{n=3}^\infty \left| k \mathscr C_n^{(x_l,x_r)}(k)\right| \leq \sum_{n=3}^\infty \frac{|k| C^{n}}{|k|^{\lfloor \frac{n+1}{2} \rfloor}} = \frac{C^4 + C^3}{k - C^2} = \bigoh( k^{-1} ). 
        \end{align}
        Since $\mathfrak c_0(k) = \bigoh(k^{-3})$, $\mathfrak c_1(k) = -m_{\mathfrak c_1}/k + \bigoh(k^{-3})$, $\mathfrak s_n(k) = m_{\mathfrak s}/k^2 + \bigoh(k^{-4})$, and since $\mathscr S _n^{(x_l,x_r)}(k) \to 0$ and $\mathscr S _n^{(x_l,x_r)}(k)$ is bounded (see Lemma~\ref{lem:Jnbounds}), we can choose $r>\sqrt{M_\gamma}$ sufficiently large so that $|\varepsilon(k)| < 1/2$ for $k\in\Omegaext(r)$, by the DCT. We have
        \begin{align}
             \frac{1}{16|k|^2}\left|m_{\mathfrak c_1} \mathfrak u_+ - 8 m_{\mathfrak s} \right| \leq |\Delta(k)| \leq \frac{3}{16r^2} \left|m_{\mathfrak c_1} \mathfrak u_+ - 8 m_{\mathfrak s} \right|\!.
        \end{align}
    \end{enumerate}
\end{proof}
\begin{rem}
    Note that for constant-coefficient IBVPs ($\alpha$, $\beta$, $\gamma$ constant), the denominator $\Delta(k)$ reduces to
    \begin{align} 
        \Delta(k) &= 2i \mathfrak a(k) \Xi(k) + i\mathfrak c_0(k) \left(\Xi(k)^2 + 1\right)  + \mathfrak s_0(k) \left(\Xi(k)^2 - 1\right)\!.
    \end{align}
    If $(a:b)_{2,4} = 0$ and $m_{\mathfrak c_0} = 0$ $($\ie $\mathfrak c_0(k) = 0$ and $\mathfrak s_0(k) = (a:b)_{1,3})$, then we require $(a:b)_{1,3} \neq 0$, so that $\Delta(k) \not\to 0$ exponentially fast (or is not identically zero). Thus, Boundary~Cases~\ref{enum:BC1}--\ref{enum:BC4} are the only allowable cases giving rise to a well-defined solution for constant-coefficient problems. If the coefficients are not constant, it may be possible to go out to higher order in the asymptotics of Lemma~\ref{lem:Delta_asymptotics}, \eg  
   $(a:b)_{2,4} = 0$, $m_{\mathfrak c_0} = 0$, $m_{\mathfrak c_1} \neq 0 $, and $m_{\mathfrak c_1} \mathfrak u_{+} - 8m_{\mathfrak s}=0$, and additional allowable boundary conditions may be identified. This requires further investigation.
\end{rem}
\begin{lemma} \label{lem:Psi_bounds}
    Consider the finite-interval, half-line, and whole-line problems. For all three, there exists an $r>\sqrt{M_\gamma}$ and $M_\Psi>0$ such that, for $k\in \Omegaext(r)$, $x\in\bar \D$, and $y\in \bar \D${\em,}
    \begin{subequations} \label{eqn:Psi_bounds}
        \begin{align} \label{eqn:Psi_bounds_direct}
            |\Psi(k,x,y)| \leq M_\Psi.
        \end{align}
        For the {\em regular problems}, 
        \begin{align} \label{eqn:Psi/Delta_bounds_reg} 
            \left| \frac{\Psi(k,x,y)}{\Delta(k)} \right| \leq M_{\Psi},
        \end{align} 
        and for the {\em irregular problems}, 
        \begin{align} \label{eqn:Psi/Delta_bounds_fi_irreg}
            \left| \frac{\Psi(k,x,y)}{\Delta(k)} \right| &\leq  M_{\Psi} \left(  1 + |k| \left( \eone{-m_{i\mathfrak n} |k|(x-x_l)} + \eone{- m_{i\mathfrak n}|k|(x_r-x)} \right) \right) \leq 3 M_\Psi|k|.
        \end{align}
    \end{subequations}
    Thus $\Psi(k,x,y)$ and $\Psi(k,x,y)/\Delta(k)$ are well-defined functions.
\end{lemma}
\begin{proof}
    For the whole-line problem, from \eqref{eqn:Psi_wl}  and Lemma~\ref{lem:Re_inu},
    \begin{align}  \label{eqn:Psi_bound_wl}
        \left|\Psi(k,x,y)\right| 
        &\leq e^{-m_{i\mathfrak n}|k||x-y|} \sum_{n=0}^\infty \frac{1}{2^nn!} \inorm{\fracbetanu }{\R}^n \leq \Etwo \!,
    \end{align}
    and \eqref{eqn:Psi_bounds_direct} follows. From Lemma~\ref{lem:Delta_asymptotics}, \eqref{eqn:Psi/Delta_bounds_reg} follows.
    For the half-line problem, with $x_l<y<x$, from \eqref{eqn:Psi_hl} 
    \begin{align} \nonumber \label{eqn:Psi_bound_hl}
        |\Psi(k,x,y)| 
        &\leq 4 \left| \etwo{ \int_y^x ik \mathfrak n(k,\xi) \, d\xi }\right|         \sum_{n=0}^\infty \sum_{\ell=0}^n \left| \left(\frac{a_0}{k \mathfrak n(k,x_l)}\mathscr S_{n-\ell}^{(x_l,y)}(k)- a_1\mathscr C_{n-\ell}^{(x_l,y)}(k) \right)\mathcal E_\ell^{(x,\infty)}(k) \right| \\
        &\leq 4e^{-m_{i\mathfrak n}|k||x-y|}\left( \frac{|a_0|}{m_{\mathfrak n} r} + |a_1| \right)\Etwo \!,
    \end{align}
    and similarly for $x_l<x<y.$ Therefore \eqref{eqn:Psi_bounds_direct} follows. From Lemma~\ref{lem:Delta_asymptotics}, we have 
    \begin{align}
        \left| \frac{\Psi(k,x,y)}{\Delta(k)}\right| &\leq \frac{4\left(|a_1| + \frac{|a_0|}{|k\mathfrak n(k,x_l)|} \right)}{\left| |a_1| - \left| \frac{a_0}{k \mathfrak n(k,x_l)}\right| \right|} \Etwo \leq 4A \Etwo \!,
    \end{align}
    where $A$ is defined in \eqref{eqn:A_hl}.
    This gives \eqref{eqn:Psi/Delta_bounds_reg}.

    For the finite-interval problem:
    \begin{enumerate}
        \item if $(a:b)_{2,4} \neq  0$, from \eqref{eqn:Psi_fi1}, we find for $x_l < y < x < r_r$,
        \begin{align} \nonumber
            |\Psi(k,x,y)| 
            &\leq 4\left( |(a:b)_{2,4}| + \frac{|(a:b)_{1,3}|}{m_{\mathfrak n}^2 r^2} + \frac{|(a:b)_{1,4}|+|(a:b)_{2,3}|}{m_{\mathfrak n} r} \right) \Etwo \\
            &~~~+ \frac{4M_\beta |(a:b)_{1,2}|}{m_\beta m_{\mathfrak n} |k|} \Etwo e^{-m_{i\mathfrak n}|k|(x_r-x_l - |x-y|)},
        \end{align}
        and similarly for $x_l < x < y < x_r$. Thus \eqref{eqn:Psi_bounds_direct} follows. From Lemma~\ref{lem:Delta_asymptotics}, \eqref{eqn:Psi/Delta_bounds_reg} follows. 
        \item If $(a:b)_{2,4} = 0$ and $m_{\mathfrak c_0} \neq 0$, from \eqref{eqn:Psi_fi1}, we find for $x_l < y < x < r_r$,
        \begin{align} \nonumber
            |\Psi(k,x,y)| 
            &\leq \frac{4}{|k|} \left( \frac{|(a:b)_{1,4}|+|(a:b)_{2,3}|}{m_{\mathfrak n}}  +\frac{|(a:b)_{1,3}|}{m_{\mathfrak n}^2 r} \right)\Etwo \\
            &~~~+ \frac{4M_\beta |(a:b)_{1,2}|}{m_\beta m_{\mathfrak n}|k|} \Etwo e^{-m_{i\mathfrak n}|k|(x_r-x_l - |x-y|)},
        \end{align}
        and similarly for $x_l < x < y < x_r$. This gives \eqref{eqn:Psi_bounds_direct}. From Lemma~\ref{lem:Delta_asymptotics}, \eqref{eqn:Psi/Delta_bounds_reg} follows.
        \item If $(a:b)_{2,4} = 0$, $m_{\mathfrak c_0} = 0$, $m_{\mathfrak c_1} = 0$, and $(a:b)_{1,3} \neq 0$, then for $x_l<y<x<x_r$,
        \def \E {\etwo{\frac12\left\|\frac{(\beta\nu)'}{(\beta\nu)}\right\|_{\D}}}
        \begin{align}
            |\Psi(k,x,y)| &\leq \frac{4}{|k|^2} \left( \frac{|(a:b)_{1,3}|}{m_{\mathfrak n}^2} 
            + \frac{4M_\beta |(a:b)_{1,2}||k|}{m_\beta m_{\mathfrak n}} e^{-m_{i\mathfrak n}|k|(x_r-x_l - |x-y|)} \right) \Etwo\!,
        \end{align}
        and similarly for $x_l<x<y<x_r$. This gives \eqref{eqn:Psi_bounds_direct}. This Boundary Case is regular if both $(a:b)_{1,2} = 0$ and $(a:b)_{3,4}=0$ and irregular if either $(a:b)_{1,2} \neq 0$ or $(a:b)_{3,4}\neq 0$, see Remark~\ref{rem:BC}. Lemma~\ref{lem:Delta_asymptotics} gives \eqref{eqn:Psi/Delta_bounds_reg} or \eqref{eqn:Psi/Delta_bounds_fi_irreg}.
        \item If $(a:b)_{2,4} = 0$, $m_{\mathfrak c_0} = 0$, $m_{\mathfrak c_1} \neq 0$, and $m_{\mathfrak c_1} \mathfrak u_+ - 8m_{\mathfrak s} \neq 0$, then, 
        \begin{align}
            \frac{(a:b)_{1,4}}{\mu(x_l)} = \frac{(a:b)_{2,3}}{\mu(x_r)} = \frac{m_{\mathfrak c_1}}2 .
        \end{align}
        From this $(a:b)_{1,4}/ \mathfrak n(k,x_l) = m_{\mathfrak c_1}/2 + \bigoh(k^{-2})$ and $(a:b)_{2,3}/\mathfrak n(k,x_r) = m_{\mathfrak c_1}/2 + \bigoh(k^{-2})$. Using Lemma~\ref{lem:Jnasymptotics}, there exists an $r>C^2$ such that
        \begin{align}
            \left|\sum_{n=1}^\infty \sum_{\ell=0}^n (-1)^{\lambda \ell} \mathcal J_{n-\ell}^{(x_l,y)}[\sigma_{p,n-\ell}](k) \mathcal J_\ell^{(x,x_r)}[\bar{\sigma}_{p,\ell}](k) \right| \leq \sum_{n=1}^\infty \sum_{\ell=0}^n \frac{C^{n-\ell}}{|k|^{\lfloor \frac{n-\ell+1}{2}\rfloor}} \frac{C^{\ell}}{|k|^{\lfloor \frac{\ell+1}{2}\rfloor}} = \frac{(k+C)^2}{(k-C^2)^2} - 1= \bigoh(k^{-1}),
        \end{align}
        for $k \in \Omegaext(r)$. For $x_l<y<x<x_r$, the $n=0$ terms involving $(a:b)_{1,4}$ and $(a:b)_{2,3}$ combine to give
        \begin{align} \nonumber
            &\frac{(a:b)_{1,4}}{k \mathfrak n(k,x_l)} \mathcal S_{0}^{(x_l,y)}(k) \mathcal C_0^{(x,x_r)}(k) - \frac{(a:b)_{2,3}}{k \mathfrak n(k,x_r)} \mathcal C_{0}^{(x_l,y)}(k) \mathcal S_0^{(x,x_r)}(k) \\ \nonumber
            &~~~= \frac{m_{\mathfrak c_1}}{2k}\left( \mathcal S_{0}^{(x_l,y)}(k) \mathcal C_0^{(x,x_r)}(k) - \mathcal C_{0}^{(x_l,y)}(k) \mathcal S_0^{(x,x_r)}(k) \right) + \bigoh (k^{-3}) \\
            &~~~= \frac{m_{\mathfrak c_1}}{2k} \sin\left( \int_{x_l}^y - \int_x^{x_r} k \mathfrak n(k,\xi) \, d\xi\right)+ \bigoh(k^{-3}),
        \end{align}
        so that, for $x_l<y<x<x_r$,
        \begin{align} 
            |\Psi(k,x,y)| &\leq 4\left\{ \frac{|m_{\mathfrak c_1}|}{4|k|} \left(\eone{-m_{i\mathfrak n} |k|(x-x_l)} + \eone{- m_{i\mathfrak n}|k|(x_r-x)} \right) + \bigoh(k^{-2}) \right\} +\frac{4M_\beta |(a:b)_{1,2}|}{m_\beta m_{\mathfrak n}|k|} e^{-m_{i\mathfrak n}|k|(x_r-x)}.
        \end{align}
        This gives \eqref{eqn:Psi_bounds_direct}. Using Lemma~\ref{lem:Delta_asymptotics}, we arrive at \eqref{eqn:Psi/Delta_bounds_fi_irreg}. The same can be shown for $x_l<x<y<x_r$.
    \end{enumerate}
\end{proof}
\begin{lemma} \label{lem:Bm_bounds}
    Consider the finite-interval and half-line problems. There exists an $r>\sqrt{M_\gamma}$ and $M_{\mathcal B}>0$ such that for $k\in \Omegaext(r)$ and $x\in \bar \D${\em,} for both the half-line $(m=0)$ and the finite-interval problem $(m=0,1)$, 
    \begin{subequations} \label{eqn:Bm_bounds}
        \begin{align} \label{eqn:Bm_bound}
            |\mathcal B_m(k,x)| \leq M_{\mathcal B}.
        \end{align}
        Further, for the half-line problem $(m=0)$, 
        \begin{align} \label{eqn:Bm/Delta_bound_hl}
            \left| \frac{\mathcal B_0(k,x)}{\Delta(k)} \right| 
            &\leq  M_{\mathcal B} |k| e^{-m_{i\mathfrak n}|k|(x-x_l)}, 
        \end{align}
        and for the finite-interval problem $(m=0,1)$, 
        \begin{align} \label{eqn:Bm/Delta_bound_fi}
            \left| \frac{\mathcal B_m(k,x)}{\Delta(k)} \right| \leq M_{\mathcal B} |k|^b \big( e^{-m_{i\mathfrak n}|k|(x_r-x)} + e^{-m_{i\mathfrak n}|k|(x-x_l)} \big).
        \end{align}
        Here, $b=1$ for {\em regular} boundary conditions, and $b=2$ for {\em irregular} boundary conditions.
    \end{subequations}   
    It follows that the functions $\mathcal B_m(k,x)$ and $\mathcal B_m(k,x)/\Delta(k)$ are well defined for the half-line and finite-interval problems.
\end{lemma}
\begin{proof}
    For the half-line problem, using Lemmas~\ref{lem:Re_inu}~and~\ref{lem:Jnbounds} in \eqref{eqn:B0_hl}, we have 
    \begin{align}
        \left| \mathcal B_0(k,x) \right| &\leq  \frac{4M_\beta}{m_\beta m_{\mathfrak n}} \Etwo \eone{-m_{i\mathfrak n} |k| (x-x_l)},
    \end{align}
    which gives \eqref{eqn:Bm_bound}. Lemma~\ref{lem:Delta_asymptotics} gives \eqref{eqn:Bm/Delta_bound_hl}. 
    Similarly, for the finite-interval problem, using Lemmas~\ref{lem:Re_inu}~and~\ref{lem:Jnbounds} in \eqref{eqn:Bm}, we have 
    \begin{align}
        \left|\mathcal B_{2-j}(k,x) \right|
        &\leq \frac{4M_\beta}{m_\beta m_{\mathfrak n}} \Etwo \left\{ \left( \frac{|a_{j1}|}{m_{\mathfrak n}|k|} + |a_{j2}| \right) e^{-m_{i\mathfrak n} |k| (x_r-x)} + \left( \frac{|b_{j1}|}{m_{\mathfrak n}|k|} + |b_{j2}| \right) e^{-m_{i\mathfrak n}|k|(x-x_l)} \right\}\!, 
    \end{align}
    for $j=1,2$,
    which gives \eqref{eqn:Bm_bound}. For the finite-interval problem with Boundary Case~\ref{enum:BC1}~or~\ref{enum:BC2} and for the irregular boundary conditions, \eqref{eqn:Bm/Delta_bound_fi} follows from the above and Lemma~\ref{lem:Delta_asymptotics}. For the regular version of Boundary Case~\ref{enum:BC3}, we have $a_{ij} = 0$ for all $i,j = 1,2,$ except for $a_{11}$ and $b_{21}$, see Remark~\ref{rem:BC}. Thus, 
    \begin{align}
        \left|\mathcal B_{2-j}(k,x) \right|
        &\leq \frac{4M_\beta}{m_\beta m_{\mathfrak n}^2 |k|} \Etwo \left( |a_{j1}| e^{-m_{i\mathfrak n} |k| (x_r-x)} + |b_{j1}| e^{-m_{i\mathfrak n}|k|(x-x_l)} \right)\!, 
    \end{align}
    from which \eqref{eqn:Bm/Delta_bound_fi} follows, using Lemma~\ref{lem:Delta_asymptotics}.
\end{proof}
\begin{lemma} \label{lem:Phi_0_bounds}
    Consider the finite-interval, half-line, and whole-line problems. For all three, there exists an $r>\sqrt{M_\gamma}$ and $M_{\Phi}>0$ such that for $k \in \Omegaext(r)$ and $x\in \bar \D,$
    \begin{subequations} \label{eqn:Phi0_bounds}
        \begin{align}
            \left|\Phi_0(k,x) \right| \leq M_\Phi \|q_0\|_\D.
        \end{align}
        For the {\em regular problems},
        \begin{align}
            \left| \frac{\Phi_{0}(k,x)}{\Delta(k)} \right| &\leq M_\Phi\|q_0\|_\D, 
        \end{align}
        and for the {\em irregular problems},
        \begin{align} \label{eqn:Phi0/Delta_bound_irreg}
            \left| \frac{\Phi_{0}(k,x)}{\Delta(k)} \right| &= M_\Phi \|q_0\|_\D  \left( 1 + |k| \left( e^{-m_{i\mathfrak n}|k|(x-x_l)} + e^{-m_{i\mathfrak n}|k|(x_r-x)} \right) \right) \leq 3 M_\Phi |k| \|q_0\|_\D.
        \end{align}
    \end{subequations}
    It follows that $\Phi_0(k,x)$ and $\Phi_0(k,x)/\Delta(k)$ are well-defined functions.
\end{lemma}
\begin{proof}
    The inequalities  \eqref{eqn:Phi_0} follow directly from Lemma~\ref{lem:Psi_bounds}.
\end{proof}
\begin{lemma} \label{lem:Phi_f_bounds}
    \begin{subequations} \label{eqn:Phi_f_bounds}
        Consider the finite-interval, half-line, and whole-line problems. For all three, there exists an $r>\sqrt{M_\gamma}$ and $M_f >0$ such that for $k\in \Omegaext(r) \backslash \Omega(r)$ (the green region of Figure~\ref{fig:Omegaext}), for $x\in \bar{\D}$, and for $t\in[0,T],$
        \begin{align} \label{eqn:Phi_f_bounds_direct}
            \big| \Phi_{f}(k,x,t) e^{-k^2t} \big| &\leq M_f \|f\|_\D. 
        \end{align}
        Further, for the {\em regular problems},
        \begin{align}\label{eqn:Phif/Delta_bounds_reg}
            \left|\frac{\Phi_{f}(k,x,t)e^{-k^2t}}{\Delta(k)}\right| &= M_f \|f\|_\D,
        \end{align}
        and for the {\em irregular problems},
        \begin{align}\label{eqn:Phif/Delta_irreg}
            \left|\frac{\Phi_{f}(k,x,t)e^{-k^2t}}{\Delta(k)}\right| &= M_{f} \|f\|_\D \left(  1 + |k| \left( \eone{-m_{i\mathfrak n} |k|(x-x_l)} + \eone{- m_{i\mathfrak n}|k|(x_r-x)} \right) \right) \leq 3M_f|k| \|f\|_\D.
        \end{align}
    \end{subequations}
    Thus, $\Phi_f(k,x,t)$ and $\Phi_f(k,x,t)/\Delta(k)$ are well-defined functions.
\end{lemma}
\begin{proof}
    For $k\in \Omegaext(r) \backslash \Omega(r)$, $|e^{-k^2(t-s)}|<1$. It follows from \eqref{eqn:q_alpha,f_alpha,ftilde_alpha} and Assumption~\enumref{ass:ffunctions}{enum:finhomogeneous} that
    \begin{align}
        \int_\D \big| \tilde f_{\alpha}(k^2,x,t) e^{-k^2t} \big| \, dx 
        \leq  \int_\D \int_0^t |f_\alpha(x,s)| \, ds dy \leq \frac{T\|f\|_\D}{m_\alpha} .
    \end{align}
    Using this and \eqref{eqn:Psi_bounds} in \eqref{eqn:Phi_f}, we obtain \eqref{eqn:Phi_f_bounds} for any $x\in \D$ and for $t\in[0,T]$.
\end{proof}
\begin{lemma} \label{lem:analytic}
    There exists an $r>\sqrt{M_\gamma}$ so that for $x\in\D$ and $t\in[0,T]$, $\mathcal J_n^{(a,b)}(k)$, $\Delta(k)$, $\Psi(k,x,y)$, and $\Phi_0(k,x)$ are analytic in $k$, for $k\in \Omegaext(r)$. The functions $\Phi_f(k,x,t)e^{-k^2t}$ and $\mathcal B_m(k,x)e^{-k^2t}$ are analytic in $k$ for $k\in \Omegaext(r) \backslash \Omega(r)$.
\end{lemma}
\begin{proof}
    Consider a closed contour $\Gamma \in \Omegaext(r)$. Then
    \begin{align}
        \oint_\Gamma \Jn \, dk &= \frac{1}{2^n} \int_{a<\cdots<b} d\mathbf y_n \oint_\Gamma dk \left(\prod_{p=1}^n \fracbetanuargs{y_p} \right)  \exp\left(\sum_{p=0}^n \sigma_{p,n}\int_{y_p}^{y_{p+1}} ik \mathfrak n(k,\xi)\, d\xi\right) = 0,
    \end{align}
    by Cauchy's theorem. We can switch the order of integration by Fubini's theorem and Lemma~\ref{lem:Jnbounds}. Therefore, by Morera's theorem, $\Jn$ is analytic for $k\in \Omegaext(r)$. For all three types of IBVPs considered, the same argument applies for the $\Delta(k)$, $\Psi(k,x,y)$, and the $\Phi_0(k,x)$ functions by Lemmas~\ref{lem:Delta_asymptotics},~\ref{lem:Psi_bounds},~and~\ref{lem:Phi_0_bounds}, and for the $\mathcal B_m(k,x)$ and $\Phi_f(k,x,t)$ functions by Lemma~\ref{lem:Bm_bounds}~and~\ref{lem:Phi_f_bounds}.
\end{proof}
The following lemmas prove that the different parts of the solution are well defined.
\begin{lemma} \label{lem:q_Bm}
    For the half-line problem $(m=0)$ and the finite-interval problem $(m=0,1)$, there exists an $r>\sqrt{M_\gamma}$ such that, for any $x\in\bar \D$ and $t\in(0,T)$, the function $q_{\mathcal B_m}(x,t)$ \eqref{eqn:q_Bm} can be written as
    \begin{subequations} \label{eqn:q_Bm_deformed2}
        \begin{align} \label{eqn:q_Bm_deformed3}
            q_{\mathcal B_m}(x,t) 
            &=  \frac{1}{2\pi} \int_{\partial\Omegaext(r)}  \frac{\mathcal B_m(k,x)} {\Delta(k)}  \mathfrak F_m(k^2,t) e^{-k^2t} \, dk,
        \end{align}
        where
        \begin{align} 
            \label{eqn:Fmfrak} 
            \mathfrak F_m(k^2,t) &= -\frac{f_m(0)}{k^2}  - \frac{1}{k^2} \int_0^t e^{k^2s} f_m'(s) \, ds,
        \end{align}
    \end{subequations}
    with the bound
    \begin{align} \label{eqn:Fmfrak_bound}
        \big| \mathfrak F_m(k^2,t)e^{-k^2t}\big| 
        &\leq \frac{\|f_m\|_\infty e^{-|k|^2\cos(2\theta)t}}{|k|^2}  + \frac{\|f_m'\|_\infty\big(1 - e^{-|k|^2\cos(2\theta)t}\big)}{|k|^4\cos(2\theta)}.
    \end{align}
    The function $q_{\mathcal B_m}(x,t)$ is well defined.
\end{lemma}
\begin{proof}
    From \eqref{eqn:F_m} and Assumption~\enumref{ass:ffunctions}{enum:fboundary}, for $k\in \Omegaext(r)\backslash \Omega(r)$,
    \begin{align}
        \big| F_m(k^2,t) e^{-k^2t} \big|
        &\leq \left| \int_0^t e^{-k^2(t-s)} f_m(s) \, ds \right| \leq T\|f_m\|_\infty .
    \end{align}
    Therefore, for $x\in \D$, we have exponential decay of the integrand of $q_{\mathcal B_m}(x,t)$ from Lemma~\ref{lem:Bm_bounds}. Using Lemma~\ref{lem:analytic}, we can deform the contour of \eqref{eqn:q_Bm} from $\Omega(r)$ to $\Omegaext(r)$.
   Assumption~\enumref{ass:ffunctions}{enum:fboundary} allows us to integrate \eqref{eqn:F_m} by parts so that
    \begin{align} 
        F_m(k^2,t)e^{-k^2t} = \frac{f_m(t)}{k^2} - \frac{f_m(0)e^{-k^2t}}{k^2} - \frac{1}{k^2} \int_0^t e^{-k^2(t-s)} f_m'(s) \, ds, 
    \end{align}
    which gives \eqref{eqn:q_Bm_deformed2}, after using Cauchy's theorem on the $f_m(t)$ term. 
    Equation \eqref{eqn:Fmfrak_bound} follows from \eqref{eqn:Fmfrak} and Assumption~\enumref{ass:ffunctions}{enum:fboundary}.
    \begin{subequations} \label{eqn:q_Bm_bounds}
        From Lemma~\ref{lem:Bm_bounds}, for the half-line problem,
        \begin{align} 
            |q_{\mathcal B_m}(x,t) |
            &\leq  \frac{M_{\mathcal B}}{2\pi} \int_{\partial\Omegaext(r)}  |k| e^{-m_{i\mathfrak n}|k|(x-x_l)} \big| \mathfrak F_m(k^2,t)e^{-k^2t}\big|  \, dk,
        \end{align}
        and for the finite-interval problem,
        \begin{align} 
            |q_{\mathcal B_m}(x,t) |
            &\leq  \frac{M_{\mathcal B}}{2\pi} \int_{\partial\Omegaext(r)}  |k|^b \left( e^{-m_{i\mathfrak n}|k|(x_r-x)} + e^{-m_{i\mathfrak n}|k|(x-x_l)} \right) \big| \mathfrak F_m(k^2,t)e^{-k^2t}\big|  \, dk.
        \end{align}
        From \eqref{eqn:q_Bm_bounds}, we see that $q_{\mathcal B_m}(x,t)$ is well defined for $x\in \bar \D$ and for $t\in (0,T)$. 
    \end{subequations}
\end{proof}
\begin{lemma} \label{lem:q_0}
    Consider the finite-interval, half-line, and whole-line problems. There exists an $r>\sqrt{M_\gamma}$ so that for $x\in\bar \D$ and $t\in (0,T)$, $q_0(x,t)$ \eqref{eqn:q_0} can be written as
    \begin{align} \label{eqn:q_0_deformed}
        q_0(x,t) &= \frac{1}{2\pi} \int_{\partial \Omegaext(r)} \frac{\Phi_0(k,x)}{\Delta(k)} e^{-k^2t} \, dk,
    \end{align}
    which is well defined.
\end{lemma}
\begin{proof}
    By Lemmas~\ref{lem:Phi_0_bounds}~and~\ref{lem:analytic}, $\Phi_0(k,x)/\Delta(k)$ is bounded, well defined, and analytic for $k\in\Omegaext(r)$. Let \break \mbox{$C_R = \{k\in \C: |k|=R$ and $\theta_0 < \theta < \pi/4$ or $3\pi/4 < \theta < \pi - \theta_0\}$}, see Figure~\ref{fig:Omegaext}. For the {\em regular problems}, using symmetry, 
    \begin{align}
        \left|\int_{C_R} \frac{\Phi_0(k,x)}{\Delta(k)} e^{-k^2t} \, dk \right| &\leq 2 M_{\Phi} \|q_0\|_\infty \int_{\theta_0}^{\frac{\pi}{4}} e^{-R^2\cos(2\theta)t} R \, d\theta 
        \leq \frac{\pi M_{\Phi} \|q_0\|_\D \big(1-e^{-R^2t}\big)}{2Rt} \to 0,
    \end{align}
    as $R\to \infty$.
    Thus we can deform the contour by Cauchy's theorem to conclude \eqref{eqn:q_0_deformed}. For the {\em irregular problems}, the above holds for the integral over the first term of \eqref{eqn:Phi0/Delta_bound_irreg} and for $x\in \D$, the second term is exponentially decaying, and we again conclude \eqref{eqn:q_0_deformed}. It follows that for all three problems
    \begin{align}
        \left|q_0(x,t)\right| = \frac{M_{\Phi} \|q_0\|_\D}{2\pi} \int_{\partial \Omegaext} \left| k e^{-k^2t} \right| \, |dk| < \infty.
    \end{align}
\end{proof}
\begin{lemma} \label{lem:q_f}
    Consider the finite-interval, half-line, and whole-line problems.  There exists an $r>\sqrt{M_\gamma}$ so that for $x\in \overline \D$ and $t\in (0,T)$, $q_f(x,t)$ \eqref{eqn:q_f} can be written as
    \begin{subequations}\label{eqn:q_f_deformed2}
        \begin{align} \label{eqn:q_f_deformed2_direct}
            q_f(x,t) 
            &= \frac{1}{2\pi} \int_{\partial\Omegaext(r)} \frac{\Phi_{\mathfrak f}(k,x,t)e^{-k^2t}}{\Delta(k)}  dk,
        \end{align}
        where
        \begin{align} 
            \label{eqn:Phif_frak}
            \Phi_{\mathfrak f}(k,x,t) 
            = \int_{\D} \frac{\Psi(k,x,y) \mathfrak f_\alpha(k^2,y,t)} { \sqrt{(\beta \mathfrak n)(k,x)} \sqrt{(\beta \mathfrak n)(k,y)}}  \, dy, 
        \end{align}
        and
        \begin{align}
            \label{eqn:f_frak} 
            \mathfrak f_\alpha(k^2,y,t) = -\frac{f_\alpha(y,0)}{k^2} - \frac{1}{k^2} \int_0^t f_{\alpha,s}(y,s) e^{k^2s} \, ds.
        \end{align}
    \end{subequations}
    Further, we have the bound
    \begin{align} \label{eqn:f_frak_bounds}
        \int_\D \big| \mathfrak f_\alpha(k^2,y,t) e^{-k^2t} \big| \, dy \leq \frac{\|f\|_\D e^{-|k|^2\cos(2\theta) t}}{m_\alpha |k|^2} + \frac{\|f_{t}\|_\D\big(1-e^{-|k|^2\cos(2\theta)t}\big)}{m_\alpha |k|^4\cos(2\theta)}.
    \end{align}
    \begin{subequations} \label{eqn:Phif_frak_bounds}
        For all three problems, there exists an $M_f>0$ such that
        \begin{align} \label{eqn:Phif_frak_bounds_direct}
            \big| \Phi_{\mathfrak f}(k,x,t) e^{-k^2t} \big|
            \leq M_f \int_{\D} \big| \mathfrak f_\alpha(k^2,y,t) e^{-k^2t} \big| \, dy.
        \end{align}
        For the {\em regular problems}
        \begin{align} \label{eqn:Phif_frak_bounds_reg}
            \left| \frac{\Phi_{\mathfrak f}(k,x,t)e^{-k^2t}}{\Delta(k)} \right|
            \leq M_f \int_{\D} \big| \mathfrak f_\alpha(k^2,y,t) e^{-k^2t} \big| \, dy, 
        \end{align}
        and for the {\em irregular problems}, 
        \begin{align} \label{eqn:Phif_frak_bounds_irreg}
            \left| \frac{\Phi_{\mathfrak f}(k,x,t)e^{-k^2t}}{\Delta(k)} \right|
            \leq M_f \left( 1 + |k| \left( e^{-m_{i\mathfrak n}|k|(x-x_l)} + e^{-m_{i\mathfrak n}|k|(x_r-x)} \right) \right) \int_{\D} \big| \mathfrak f_\alpha(k^2,y,t) e^{-k^2t} \big| \, dy.
        \end{align}
    \end{subequations}
    It follows that $q_f(x,t)$ is well defined for all three problems.
\end{lemma}
\begin{proof}
    By Lemmas~\ref{lem:Phi_f_bounds}~and~\ref{lem:analytic}, $\Phi_f(k,x,t)e^{-k^2 t}/\Delta(k)$ is bounded, well defined, and analytic for $k\in\Omegaext(r)/\Omega(r)$. Let $C_R$ be defined as in the proof of Lemma~\ref{lem:q_0}, see Figure~\ref{fig:Omegaext}. Then, for the {\em regular problems}, using symmetry, 
    \begin{align}
        \left|\int_{C_R} \frac{\Phi_f(k,x,t) e^{-k^2t}}{\Delta(k)} \, dk \right| 
        &\leq 2 M_{f} \|f\|_\D  \int_0^{\frac{\pi}{4}} R e^{-R^2\cos(2\theta)t} \, d\theta \to 0,
    \end{align} 
    as $R\to\infty$. Thus, we can deform the integral in \eqref{eqn:q_f} from $\Omega(r)$ to $\Omegaext(r)$. For the {\em irregular problems}, the above holds for the integral over the first term of \eqref{eqn:Phif/Delta_irreg} and for $x\in \D$ the second term is exponentially decaying. Thus, we can still deform from $\Omega(r)$ to $\Omegaext(r)$. Using Assumption~\enumref{ass:ffunctions}{enum:finhomogeneous}, we can integrate \eqref{eqn:q_alpha,f_alpha,ftilde_alpha} by parts, to obtain
    \begin{align}
        \tilde f_\alpha(k^2,x,t) &= \frac{f_\alpha(x,t) e^{k^2t} - f_\alpha(x,0)}{k^2} - \frac{1}{k^2} \int_0^t f_{\alpha,s}(x,s) e^{k^2s} \, ds,
    \end{align}
    which gives \eqref{eqn:q_f_deformed2}, after using Cauchy's theorem on the $f_\alpha(x,t)$ term. Equation \eqref{eqn:f_frak_bounds} follows directly from \eqref{eqn:f_frak} and Assumption~\enumref{ass:ffunctions}{enum:finhomogeneous}, and equation \eqref{eqn:Phif_frak_bounds} follows from Lemma~\ref{lem:Psi_bounds}. From \eqref{eqn:Phif_frak_bounds}, we see that the integrand in $q_f(x,t)$ is absolutely integrable and is therefore well defined for all $x\in \bar \D$ and any $t\in(0,T)$ (or for any $x\in \D$ and for all $t\in[0,T]$).
\end{proof}
Finally, we combine all the results obtained. 
\begin{theorem} \label{thm:welldefined}
    There exists an $r>\sqrt{M_\gamma}$ such that the functions \eqref{eqn:q_wl}, \eqref{eqn:q_hl}, and \eqref{eqn:q_fi} are well defined for all $x\in \bar \D$ and for any $t\in(0,T)$.
\end{theorem}
\begin{proof}
    Combining Lemmas~\ref{lem:q_Bm},~\ref{lem:q_0},~and~\ref{lem:q_f}, we obtain our result.
\end{proof}
\section{Proofs: the solution expressions solve the evolution equation} \label{sec:proofs_evolution}
\newcommand{\Jn}{{\mathcal{J}}_n^{(a,b)}[\sigma_{p,n}](k)}
In this appendix, we prove that the solution expressions \eqref{eqn:q_wl}, \eqref{eqn:q_hl}, and \eqref{eqn:q_fi} for the whole-line, half-line, and finite-interval problems, respectively, solve the evolution equation \eqref{eqn:IVP} in their respective domains. Naturally, we are in need of lemmas on the derivatives of various quantities defining the solution expressions. The following lemma deals with derivatives with respect to the spatial variable. 
\begin{lemma} \label{lem:Jnderivatives}
    For $n\geq 0$, 
   the derivatives of $\mathcal E_n^{(x,\infty)}(k)$ and $\tilde{\mathcal E}_n^{(-\infty,x)}(k)$ are given by
    \begin{subequations} \label{eqn:Enderivatives}
        \begin{align} 
            \partial_x \mathcal E_n^{(x,\infty)}(k) 
            &= -\frac12 \fracbetanuargs{x} \mathcal E_{n-1}^{(x,\infty)}(k) - (1-(-1)^{n})ik \mathfrak n(k,x)\mathcal E_n^{(x,\infty)}(k), \\
            \label{eqn:Entilde_derivative}
            \partial_x \tilde{\mathcal E}_{n}^{(-\infty,x)}(k) 
            &= \frac12 \fracbetanuargs{x} \tilde{\mathcal E}_{n-1}^{(-\infty,x)}(k) + (1-(-1)^{n}) ik\mathfrak n(k,x) \tilde{\mathcal E}_{n}^{(-\infty,x)}(k),
            \intertext{and those of $\mathcal C_n^{(a,b)}(k)$ and $\mathcal S_n^{(a,b)}(k)$ are}
            \partial_x \mathcal C_{n}^{(x_l,x)}(k) &= \frac12 \fracbetanuargs{x} \mathcal C_{n-1}^{(x_l,x)}(k) - (-1)^{n} k \mathfrak n(k,x) \mathcal S_{n}^{(x_l,x)}(k), \\
            \partial_x \mathcal C_{n}^{(x,x_r)}(k) &= -\frac12 \fracbetanuargs{x} \mathcal C_{n-1}^{(x,x_r)}(k) + k \mathfrak n(k,x) \mathcal S_{n}^{(x,x_r)}(k), \\
            \partial_x \mathcal S_{n}^{(x_l,x)}(k) &= \frac12 \fracbetanuargs{x} \mathcal S_{n-1}^{(x_l,x)}(k) + (-1)^{n}k \mathfrak n(k,x)\mathcal C_{n}^{(x_l,x)}(k), \\
            \partial_x \mathcal S_{n}^{(x,x_r)}(k) &= \frac12 \fracbetanuargs{x} \mathcal S_{n-1}^{(x,x_r)}(k) - k \mathfrak n(k,x) \mathcal C_{n}^{(x,x_r)}(k).
        \end{align}
    \end{subequations}
\end{lemma}
\begin{proof}
    Since $(\beta\mathfrak n)'/(\beta \mathfrak n) \in L^1(\D)$ by Lemma~\ref{lem:beta_nu}, the proof is by direct calculation of the derivatives of \eqref{eqn:En&Entilde} and \eqref{eqn:CnSn} \cite{analysis_mcdonald}. We show one such calculation. From \eqref{eqn:Entilde}, 
    \begin{align}
        \partial_x \tilde{\mathcal E}_n^{(a,x)}(k) 
        &= \frac{\partial_x}{2^n} \int_a^x dy_1 \int_{y_1}^x dy_2 \cdots \int_{y_{n-2}}^x dy_{n-1} \int_{y_{n-1}}^x dy_n \left(\prod_{p=1}^n \frac{(\beta\mathfrak n)'}{(\beta\mathfrak n)} \right) \exp\left(\sum_{p=0}^n (1-(-1)^{p}) \int_{y_p}^{y_{p+1}} ik \mathfrak n(k,\xi)\,d\xi\right) \!,
    \end{align}
    so that
    \begin{align} \nonumber
        \partial_x \tilde{\mathcal E}_n^{(a,x)}(k)
        &= \frac{(\beta\mathfrak n)'(k,x)}{(\beta\mathfrak n)(k,x)} \frac{1}{2^n} \int_{a=y_0<\cdots<y_{n}=x} \left(\prod_{p=1}^{n-1} \frac{(\beta\mathfrak n)'}{(\beta\mathfrak n)} \right) \exp\left(\sum_{p=0}^{n-1} (1-(-1)^{p}) \int_{y_p}^{y_{p+1}} ik \mathfrak n (k,\xi)\,d\xi\right) d\mathbf{y}_{n-1} \\ 
        &~~~+ (1-(-1)^n)\frac{ik \mathfrak n(k,x)}{2^n} \int_{a=y_0<\cdots<y_{n+1}=x} \left(\prod_{p=1}^n \frac{(\beta\mathfrak n)'}{(\beta\mathfrak n)} \right) \exp\left(\sum_{p=0}^n (1-(-1)^{p}) \int_{y_p}^{y_{p+1}} i\nu(k,\xi)\,d\xi\right) d\mathbf{y}_n,
    \end{align}
    which gives \eqref{eqn:Entilde_derivative}.
\end{proof}
In Lemma~\ref{lem:Jn_identity}, we prove a general summation identity for the generalized accumulation functions $\Jn$. This identity is used to prove the problem-specific identities in Lemma~\ref{lem:Delta_identity}. In turn, these are used to prove the relation between $\chi(k,x)$ \eqref{eqn:delta_def} and $\Delta(k)$ in Lemma~\ref{lem:delta(k,x)}.
\begin{lemma} \label{lem:Jn_identity}
    Let $\bar \sigma_{p,n-\ell}$ and $\tilde{\sigma}_{p,\ell}$ be two non-negative integer-valued functions as described in Definition~\ref{def:Jndefinitions}. Denote
    \begin{subequations}\label{eqn:Jn_identity}
        \begin{align} \label{eqn:Jn_identity_ass}
            \sigma_{p,n} = \case{1}{ \bar \sigma_{p,n-\ell}, & \text{ if } 0\leq p \leq n-\ell, \\ \tilde \sigma_{p-(n-\ell),\ell}, & \text{ if } n-\ell < p \leq n.} 
        \end{align}
        If $\bar \sigma_{n-\ell,n-\ell} = \tilde \sigma_{0,\ell}$ and $\sigma_{p,n}$ is independent of $\ell$, then
        \begin{align} \label{eqn:Jn_identity_eq}
            \mathcal J_n^{(a,b)}[\sigma_{p,n}](k) &= \sum_{\ell=0}^n \mathcal J_{n-\ell}^{(a,x)}[\bar \sigma_{p,n-\ell}](k) \mathcal J_{\ell}^{(x,b)}[\tilde \sigma_{p,\ell}](k).
        \end{align}
    \end{subequations}
\end{lemma}
\begin{proof}
    Define
    \begin{align}
        \mathfrak j_n^{(a,b)}[\sigma_{p,n}](k) &= \sum_{\ell=0}^n \mathcal J_{n-\ell}^{(a,x)}[\bar \sigma_{p,n-\ell}](k) \mathcal J_{\ell}^{(x,b)}[\tilde \sigma_{p,\ell}](k).
    \end{align}
    By the definition of $\Jn$ \eqref{eqn:Jndefinition},
    \begin{align} \nonumber
        \mathfrak j_n^{(a,b)}[\sigma_{p,n}](k)
        &= \frac{1}{2^n} \sum_{\ell=0}^n \int_{a<\cdots<x} \left(\prod_{p=1}^{n-\ell} \fracbetanuargs{y_p} \right) \exp\left( \sum_{p=0}^{n-\ell} \bar \sigma_{p,n-\ell} \int_{y_p}^{y_{p+1}} i k \mathfrak n(k,\xi) \, d\xi \right) dy_1 \cdots dy_{n-\ell} \times \\
        &~~~ \times \int_{x < \cdots < b} \Bigg(\prod_{p=n-\ell+1}^{n} \fracbetanuargs{y_p} \Bigg) \exp\left(\sum_{p=n-\ell}^{n} \tilde \sigma_{p-(n-\ell),\ell} \int_{y_p}^{y_{p+1}} ik \mathfrak n(k,\xi) \, d\xi \right) dy_{n-\ell+1} \cdots dy_n. 
    \end{align}
    In the exponential of the first integral, for the $p=n-\ell$ term, $y_{n-\ell+1}$ is defined as $x$. In the exponential in the second integral, for the $p=n-\ell$ term, $y_{n-\ell} = x$. Since $\bar \sigma_{n-\ell,n-\ell} = \tilde \sigma_{0,\ell} = \sigma_{n-\ell,n}$, multiplying the exponentials and adding these terms together, we have
    \begin{align}
        \bar \sigma_{n-\ell,n-\ell} \int_{y_{n-\ell}}^x ik\mathfrak n(k,\xi) \, d\xi + \tilde \sigma_{0,\ell} \int_{x}^{y_{n-\ell+1}} ik \mathfrak n(k,\xi) \, d\xi = \sigma_{n-\ell,n} \int_{y_{n-\ell}}^{y_{n-\ell+1}} ik \mathfrak n(k,\xi) \, d\xi,
    \end{align}
    and the two integrals are combined as
    \begin{align}
        \mathfrak j_n^{(a,b)}[\sigma_{p,n}](k)
        &= \frac{1}{2^n} \sum_{\ell=0}^n \int_{a<\cdots<y_{n-\ell}<x< y_{n-\ell+1} < \cdots  < b} \left(\prod_{p=1}^{n} \fracbetanuargs{y_p} \right) \exp\left(\sum_{p=0}^{n} \sigma_{p,n} \int_{y_p}^{y_{p+1}} ik \mathfrak n(k,\xi) \, d\xi \right) dy_1 \cdots dy_{n}.
    \end{align}
    Summing over $\ell$ is equivalent to adding up all possibilities of $x$ lying between one of the $y_1,\ldots, y_n$. Since $\sigma_{p,n}$ is independent of $\ell$ by assumption, the integrand is independent of $\ell$, and
    \begin{align}
        \mathfrak j_n^{(a,b)}[\sigma_{p,n}](k)
        &= \frac{1}{2^n} \int_{a= y_0<\cdots< y_{n+1} < b} \left(\prod_{p=1}^{n} \fracbetanuargs{y_p} \right) \exp\left(\sum_{p=0}^{n} \sigma_{p,n} \int_{y_p}^{y_{p+1}} ik \mathfrak n(k,\xi) \, d\xi \right) dy_1 \cdots dy_{n},
    \end{align}
    which is \eqref{eqn:Jn_identity}.
\end{proof}
From the identity in Lemma~\ref{lem:Jn_identity}, we can prove the following more specific forms of  \eqref{eqn:Jn_identity_eq}.
\begin{lemma} \label{lem:Delta_identity}
    For the whole-line problem, if $n$ is even,
    \begin{subequations} \label{eqn:Delta_id}
        \begin{align} 
            \label{eqn:Delta_id_wl}
             \mathcal E_n^{(-\infty,\infty)}(k) &= \sum_{\ell=0}^n \tilde{\mathcal E}_{n-\ell}^{(-\infty, x)}(k) {\mathcal E}_\ell^{(x,\infty)}(k).
        \end{align}
        For the half-line problem, for any $n$,
        \begin{align}
             \label{eqn:Delta_id_hl}
             \mathcal E_n^{(x_l,\infty)}(k) &=  \sum_{\ell=0}^n \left( \mathscr C_{n-\ell}^{(x_l,x)}(k)  -  (-1)^n i\mathscr S_{n-\ell}^{(x_l,x)}(k) \right) {\mathcal E}_\ell^{(x,\infty)}(k).
        \end{align}
        Finally, for the finite-interval problem, for any $n$,
        \begin{align}
             \label{eqn:Delta_id_fiC}
             \mathscr C_n^{(x_l,x_r)} (k) &= \sum_{\ell=0}^n \left( \mathscr C_{n-\ell}^{(x_l,x)}(k) \mathscr C_\ell^{(x,x_r)}(k) -  (-1)^{n-\ell} \mathscr S_{n-\ell}^{(x_l,x)}(k) \mathscr S_\ell^{(x,x_r)}(k) \right)\!, \\
             \label{eqn:Delta_id_fiS}
             \mathscr S_n^{(x_l,x_r)}(k) &= \sum_{\ell=0}^n \left( \mathscr S_{n-\ell}^{(x_l,x)}(k) \mathscr C_\ell^{(x,x_r)}(k) + (-1)^{n-\ell} \mathscr C_{n-\ell}^{(x_l,x)}(k)  \mathscr S_\ell^{(x,x_r)}(k) \right) \!.
        \end{align}
    \end{subequations}
\end{lemma}
\begin{proof}
    For the whole-line problem, we define $\mathfrak e_\text{wl}$ as the right-hand side of \eqref{eqn:Delta_id_wl}.
    Using \eqref{eqn:relations}, we write $\mathfrak e_\text{wl}$ as
    \begin{align}
        \mathfrak e_\text{wl} &= \sum_{\ell=0}^n \mathcal J_{n-\ell}^{(-\infty, x)}[1-(-1)^p](k) \mathcal J_\ell^{(x,\infty)}[1-(-1)^{\ell-p}](k).
    \end{align}
    From Lemma~\ref{lem:Jn_identity}, if $n$ is even, $\bar \sigma_{p,n-\ell} = 1-(-1)^p$ and $\tilde \sigma_{p,\ell} = 1-(-1)^{\ell-p}$ so that $\bar \sigma_{n-\ell,n-\ell} =\tilde \sigma_{0,\ell}$ and 
    \def \spone {\hspace*{-0.2in}}
    \def \sptwo {\hspace*{-0.1in}}
    \begin{align}
        \sigma_{p,n} = \case{1}{ \bar \sigma_{p,n-\ell}, & \spone \text{if } 0\leq p \leq n-\ell, \\ \tilde \sigma_{p-(n-\ell),\ell}, & \spone \text{if } n-\ell < p \leq n,} \sptwo = \case{1}{ 1- (-1)^p, & \spone \text{if } 0\leq p \leq n-\ell, \\ 1-(-1)^{\ell-(p-(n-\ell))}, & \spone \text{if } n-\ell < p \leq n,} \sptwo = 1-(-1)^{n-p},
    \end{align} 
    is independent of $\ell$ so that \eqref{eqn:Delta_id_wl} follows.
    
    For the half-line problem, we define $\mathfrak e_\text{hl}$ as the right-hand side of \eqref{eqn:Delta_id_hl}.
    Using \eqref{eqn:relations},
    \begin{align}
        \mathfrak e_\text{hl} &= \frac12 \sum_{\ell=0}^n \left( (1-(-1)^n) \mathcal J_{n-\ell}^{(x_l,x)}[1+(-1)^p](k) + (1+(-1)^n) \mathcal J_{n-\ell}^{(x_l,x)}[1-(-1)^p](k) \right) \mathcal J_\ell^{(x,\infty)}[1-(-1)^{\ell-p}](k),
    \end{align}
    which is simplified to
    \begin{align}
        \mathfrak e_\text{hl}
        &= \sum_{\ell=0}^n  \mathcal J_{n-\ell}^{(x_l,x)}[1-(-1)^{n-p}](k)  \mathcal J_\ell^{(x,\infty)}[1-(-1)^{\ell-p}](k).
    \end{align}
    From Lemma~\ref{lem:Jn_identity}, $\bar \sigma_{p,n-\ell} = 1-(-1)^{n-p}$ and $\tilde \sigma_{p,\ell} = 1-(-1)^{\ell-p}$ so that $\bar \sigma_{n-\ell,n-\ell} = \tilde \sigma_{0,\ell}$ and \mbox{$\sigma_{p,n} = 1-(-1)^{n-p}$}. Equation \eqref{eqn:Delta_id_hl} follows.
    
    For the finite-interval problem, we define $\mathfrak e_c$ and $\mathfrak e_s$ as the right-hand side of \eqref{eqn:Delta_id_fiC} and \eqref{eqn:Delta_id_fiS}, respectively.
    Using \eqref{eqn:relations}, we write these in terms of $\Jn$, obtaining
    \begin{subequations}
        \begin{align} \nonumber
             \mathfrak e_c 
             &= \frac{1}{4} \sum_{\ell=0}^n (1+(-1)^{n-\ell}) \left( \mathcal J_{n-\ell}^{(x_l,x)}[1+(-1)^p](k) \mathcal J_{\ell}^{(x,x_r)}[1+(-1)^p](k) \right. \\
             \nonumber
             &\hspace{1.5in}\left.+ \mathcal J_{n-\ell}^{(x_l,x)}[1-(-1)^p](k) \mathcal J_{\ell}^{(x,x_r)}[1-(-1)^p](k) \right)\\
             \nonumber
             &~~~+ \frac{1}{4} \sum_{\ell=0}^n (1-(-1)^{n-\ell}) \left( \mathcal J_{n-\ell}^{(x_l,x)}[1+(-1)^p](k) \mathcal J_{\ell}^{(x,x_r)}[1-(-1)^p](k) \right. \\
             &\hspace{1.5in}\left. + \mathcal J_{n-\ell}^{(x_l,x)}[1-(-1)^p](k) \mathcal J_{\ell}^{(x,x_r)}[1+(-1)^p](k) \right)\!, \\
             \nonumber
             \mathfrak e_s 
             &= \frac{1}{4i} \sum_{\ell=0}^n (1+(-1)^{n-\ell})\left( 
             \mathcal J_{n-\ell}^{(x_l,x)}[1+(-1)^p](k) 
             \mathcal J_{\ell}^{(x,x_r)}[1+(-1)^p](k) \right. \\
             \nonumber
             &\hspace{1.5in}\left.- \mathcal J_{n-\ell}^{(x_l,x)}[1-(-1)^p](k) \mathcal J_{\ell}^{(x,x_r)}[1-(-1)^p](k) \right) \\
             \nonumber
             &~~~+ \frac{1}{4i} \sum_{\ell=0}^n 
             (1-(-1)^{n-\ell}) \left( \mathcal J_{n-\ell}^{(x_l,x)}[1+(-1)^p](k)
             \mathcal J_{\ell}^{(x,x_r)}[1-(-1)^p](k) \right. \\
             &\hspace{1.5in}\left.- \mathcal J_{n-\ell}^{(x_l,x)}[1-(-1)^p](k) \mathcal J_{\ell}^{(x,x_r)}[1+(-1)^p](k) \right)\!,
        \end{align}
    \end{subequations}
    which simplify to
    \begin{subequations}
        \begin{align}
             \mathfrak e_c 
             &= \frac{1}{2} \sum_{\ell=0}^n \left( \mathcal J_{n-\ell}^{(x_l,x)}[1+(-1)^p](k) \mathcal J_{\ell}^{(x,x_r)}[1+(-1)^{n-\ell+p}](k) + \mathcal J_{n-\ell}^{(x_l,x)}[1-(-1)^p](k) \mathcal J_{\ell}^{(x,x_r)}[1-(-1)^{n-\ell+p}](k) \right)\!, \\
             \mathfrak e_s 
             &= \frac{1}{2i} \sum_{\ell=0}^n \left( 
             \mathcal J_{n-\ell}^{(x_l,x)}[1+(-1)^p](k) 
             \mathcal J_{\ell}^{(x,x_r)}[1+(-1)^{n-\ell+p}](k) - \mathcal J_{n-\ell}^{(x_l,x)}[1-(-1)^p](k) \mathcal J_{\ell}^{(x,x_r)}[1-(-1)^{n-\ell+p}](k) \right) \!.
        \end{align}
    \end{subequations}
    For the first terms of $\mathfrak e_c$ and $\mathfrak e_s$, $\bar \sigma_{p,n-\ell} = 1+(-1)^p$, and $\tilde \sigma_{p,\ell} = 1+(-1)^{n-\ell+p}$ so that $\bar \sigma_{n-\ell,n-\ell} = \tilde \sigma_{0,\ell}$ and $\sigma_{p,n} = 1+(-1)^p$. 
    For the second terms of $\mathfrak e_c$ and $\mathfrak e_s$, $\bar \sigma_{p,n-\ell} = 1-(-1)^p$, and $\tilde \sigma_{p,\ell} = 1-(-1)^{n-\ell+p}$ so that $\bar \sigma_{n-\ell,n-\ell} = \tilde \sigma_{0,\ell}$ and $\sigma_{p,n} = 1-(-1)^p$.
    Equations \eqref{eqn:Delta_id_fiC} and \eqref{eqn:Delta_id_fiS} follow.
\end{proof}
Now, we begin taking derivatives of the solution expressions. In Definition~\ref{def:Psibar_tilde}, we introduce some functions that appear in the derivatives of the solution expressions. In Lemmas~\ref{lem:Psibar}--\ref{lem:delta(k,x)}, we prove some properties of these functions.
\begin{definition} \label{def:Psibar_tilde}
    We define
    \begin{subequations} \label{eqn:Psibar_tilde}
        \def \sp {-0.55in}
        \begin{align}
            \label{eqn:Psibar_def}
            \bar \Psi(k,x,y) &= \sqrt{(\beta\mathfrak n)(k,x)} \frac{\partial}{\partial x}\left(\frac{\Psi(k,x,y)}{\sqrt{(\beta \mathfrak n)(k,x)}} \right) &&\hspace*{\sp}= \Psi_x(k,x,y) - \frac12 \fracbetanuargs{x} \Psi(k,x,y), ~~~\mbox{and}\\
            \label{eqn:Psitilde_def}
            \tilde \Psi(k,x,y) &=\sqrt{(\beta\mathfrak n)(k,x)} \frac{\partial}{\partial x}\left(\frac{(\beta \bar \Psi)(k,x,y)}{\sqrt{(\beta \mathfrak n)(k,x)}} \right) &&\hspace*{\sp}= (\beta\bar\Psi)_x(k,x,y) - \frac12 \fracbetanuargs{x} \left(\beta\bar\Psi \right)(k,x,y),
        \end{align}
    \end{subequations}
    where we use the notation $(\beta \bar \Psi)(k,x,y) = \beta(x) \bar \Psi(k,x,y)$. We also define
    \begin{subequations} \label{eqn:Psilimits}
        \begin{align} \label{eqn:Psilimit_def}
            \Psi(k,x,x^\pm) = \lim_{y\to x^\pm} \Psi(k,x,y),~~~~~~
            \bar \Psi(k,x,x^\pm) = \lim_{y\to x^\pm} \bar \Psi(k,x,y),
        \end{align}
        and
        \begin{align} \label{eqn:delta_def}
            \chi(k,x) = (\beta \bar \Psi)(k,x,x^-) - (\beta \bar \Psi)(k,x,x^+).
        \end{align}
    \end{subequations}
\end{definition}

\begin{lemma} \label{lem:Psibar}
    For the whole-line problem, for $y<x$, 
    \begin{subequations} \label{eqn:Psibar_wl}
        \begin{align} \label{eqn:Psibar_wl_lower}
            \bar \Psi(k,x,y) &= ik \mathfrak n(k,x)\exp\left( \int_y^x ik \mathfrak n(k,\xi) \, d\xi\right) \sum_{n=0}^\infty \sum_{\ell=0}^n  \tilde{\mathcal E}_{n-\ell}^{(-\infty,y)} (k) \mathcal E_\ell^{(x,\infty)}(k), 
        \end{align}
        and for $x < y$,
        \begin{align} \label{eqn:Psibar_wl_upper}
            \bar \Psi(k,x,y) &= - ik \mathfrak n(k,x) \exp\left( \int_x^y ik \mathfrak n(k,\xi) \, d\xi\right) \sum_{n=0}^\infty \sum_{\ell=0}^{n} (-1)^{n}  \tilde{\mathcal E}_{n-\ell}^{(-\infty,x)}(k) \mathcal E_\ell^{(y,\infty)}(k).
        \end{align}
    \end{subequations}
    For the half-line problem, for $x_l<y<x$,
    \begin{subequations} \label{eqn:Psibar_hl}
        \begin{align} \label{eqn:Psibar_hl_lower}
            \bar \Psi(k,x,y) &=  4ik \mathfrak n(k,x) \exp\left(\int_{x_l}^{x} ik  \mathfrak n(k,\xi) \, d\xi\right)\sum_{n=0}^\infty \sum_{\ell=0}^n \left(\frac{a_0}{k \mathfrak n(k,x_l)} \mathcal S_{n-\ell}^{(x_l,y)}(k)- a_1\mathcal C_{n-\ell}^{(x_l,y)}(k) \right) \mathcal E_\ell^{(x,\infty)}(k), 
        \end{align}
        and for $x_l<x<y$,
        \begin{align} \label{eqn:Psibar_hl_upper}
            \bar \Psi(k,x,y) 
            &= 4 k \mathfrak n(k,x) \exp\left(\int_{x_l}^{y} ik \mathfrak n(k,\xi) \, d\xi\right)\sum_{n=0}^\infty \sum_{\ell=0}^n (-1)^n \left(\frac{a_0}{k \mathfrak n(k,x_l)} \mathcal C_{n-\ell}^{(x_l,x)}(k) + a_1 \mathcal S_{n-\ell}^{(x_l,x)}(k) \right) {\mathcal E}_\ell^{(y,\infty)}(k) .
        \end{align}
    \end{subequations}
    For the finite-interval problem, for $x_l<y<x<x_r$,
    \begin{subequations} \label{eqn:Psibar_fi}
        \begin{align} \label{eqn:Psibar_fi_lower} \nonumber
            \bar \Psi(k,x,y)
            &= 4k \mathfrak n(k,x) {{ \Eint}} \left\{ -\frac{\beta(x_r)(a:b)_{1,2}}{k \sqrt{(\beta\mathfrak n)(k,x_l)}\sqrt{(\beta\mathfrak n)(k,x_r)}} \sum_{n=0}^\infty   (-1)^n\mathcal C_n^{(y,x)}(k) \right. \\  \nonumber
            &\hspace*{0.75in}-(a:b)_{2,4} \sum_{n=0}^\infty \sum_{\ell=0}^n (-1)^\ell \mathcal C_{n-\ell}^{(x_l,y)}(k)  \mathcal S_\ell^{(x,x_r)}(k) - \frac{(a:b)_{1,3}}{k^2 \mathfrak n(k,x_l) \mathfrak n(k,x_r)}\sum_{n=0}^\infty \sum_{\ell=0}^n \mathcal S_{n-\ell}^{(x_l,y)}(k) \mathcal C_\ell^{(x,x_r)}(k) \\ 
            &\hspace*{0.75in}\left. + \frac{(a:b)_{1,4}}{k \mathfrak n(k,x_l)} \sum_{n=0}^\infty \sum_{\ell=0}^n (-1)^\ell \mathcal S_{n-\ell}^{(x_l,y)}(k) \mathcal S_\ell^{(x,x_r)}(k) + \frac{(a:b)_{2,3}}{k \mathfrak n(k,x_r)} \sum_{n=0}^\infty \sum_{\ell=0}^n \mathcal C_{n-\ell}^{(x_l,y)}(k) \mathcal C_\ell^{(x,x_r)}(k)    \right\}\!,
        \end{align}
        and for $x_l<x<y<x_r$,
        \begin{align} \label{eqn:Psibar_fi_upper} \nonumber
            \bar \Psi(k,x,y) 
            &= 4 k \mathfrak n(k,x) {{\Eint}} \left\{ \frac{\beta(x_l) (a:b)_{3,4}}{k \sqrt{(\beta\mathfrak n)(k,x_l)}\sqrt{(\beta\mathfrak n)(k,x_r)}} \sum_{n=0}^\infty   \mathcal C_n^{(x,y)}(k) \right. \\ \nonumber
            &\hspace*{0.25in}+(a:b)_{2,4} \sum_{n=0}^\infty \sum_{\ell=0}^n (-1)^n \mathcal S_{n-\ell}^{(x_l,x)}(k) \mathcal C_\ell^{(y,x_r)}(k) + \frac{(a:b)_{1,3}}{k^2 \mathfrak n(k,x_l)\mathfrak n(k,x_r)} \sum_{n=0}^\infty \sum_{\ell=0}^n (-1)^{n-\ell} \mathcal C_{n-\ell}^{(x_l,x)}(k)  \mathcal S_\ell^{(y,x_r)}(k) \\ 
            &\hspace*{0.25in}\left.+ \frac{(a:b)_{1,4}}{k \mathfrak n(k,x_l)} \sum_{n=0}^\infty \sum_{\ell=0}^n (-1)^n  \mathcal C_{n-\ell}^{(x_l,x)}(k) \mathcal C_\ell^{(y,x_r)}(k)+ \frac{(a:b)_{2,3}}{k \mathfrak n(k,x_r)} \sum_{n=0}^\infty \sum_{\ell=0}^n (-1)^{n-\ell} \mathcal S_{n-\ell}^{(x_l,x)}(k) \mathcal S_\ell^{(y,x_r)}(k)  \right\}\!.
        \end{align}
    \end{subequations}
\end{lemma}
\begin{proof}
    Using \eqref{eqn:Enderivatives} in \eqref{eqn:Psi_wl}, \eqref{eqn:Psi_hl}, and \eqref{eqn:Psi_fi}, we find \eqref{eqn:Psibar_wl}, \eqref{eqn:Psibar_hl}, and \eqref{eqn:Psibar_fi} for the whole-line, half-line, and finite-interval problems, respectively.
\end{proof}
\begin{lemma} \label{lem:Psibar_bounds}
    Consider the finite-interval, half-line, and whole-line problems. There exists an $r>\sqrt{M_\gamma}$ and $M_\Psi>0$ so that for $k\in \Omegaext(r)$, for $x\in\bar \D$, and for $y\in \bar \D$ \\
    \begin{subequations}
        \begin{align} \label{eqn:Psibar_bounds}
            |\bar \Psi(k,x,y)| \leq M_{\Psi} |k|. 
        \end{align}
        For the {\em regular problems} 
        \begin{align} \label{eqn:Psibar/Delta_bounds_reg}        
            \left| \frac{\bar\Psi(k,x,y)}{\Delta(k)} \right| \leq M_{\Psi} |k|,
        \end{align} 
        and for the {\em irregular problems} 
        \begin{align} \label{eqn:Psibar/Delta_bounds_fi_irreg}
            \left| \frac{\bar \Psi(k,x,y)}{\Delta(k)} \right| &\leq  M_{\Psi} |k| \left(  1 + |k| \left( \eone{-m_{i\mathfrak n} |k|(x-x_l)} + \eone{- m_{i\mathfrak n}|k|(x_r-x)}  \right) \right)\!.
        \end{align}
    \end{subequations}
    Therefore, $\bar \Psi(k,x,y)$ and $\bar \Psi(k,x,y)/\Delta(k)$ are well-defined functions.
\end{lemma}
\begin{proof}
    The proof is identical to that of Lemma~\ref{lem:Psi_bounds} in Appendix~\ref{sec:proofs_welldefined}. Note that $M_\Psi$ here and from Lemma~\ref{lem:nubounds} are identical up to a factor of $M_{\mathfrak n}$. Without loss of generality, we take them to be the same.
\end{proof}
\begin{lemma} \label{lem:delta(k,x)}
    For the finite-interval, half-line, and whole-line problems, 
    \begin{align} \label{eqn:chi_identity}
        \chi(k,x) &= 2ik (\beta\mathfrak n)(k,x) \Delta(k),
    \end{align}
    where $\chi(k,x)$ is defined in \eqref{eqn:delta_def}.
\end{lemma}
\begin{proof}
    For the whole-line problem, using Lemma~\ref{lem:Psibar} in \eqref{eqn:delta_def},
    \begin{align}
        \chi(k,x) 
        &= ik (\beta\mathfrak n)(k,x) \sum_{n=0}^\infty (1+ (-1)^{n}) \sum_{\ell=0}^n  \tilde{\mathcal E}_{n-\ell}^{(-\infty,x)}(k) \mathcal E_\ell^{(x,\infty)}(k),
    \end{align}
    which gives \eqref{eqn:chi_identity}, using Lemma~\ref{lem:Delta_identity}.
    Similarly, for the half-line problem, 
    \begin{align} 
        \chi(k,x) 
        &= - 4k (\beta\mathfrak n)(k,x) \exp\left(\int_{x_l}^{x}ik \mathfrak n(k,\xi) \, d\xi\right)\sum_{n=0}^\infty \left( \frac{(-1)^n a_0}{k \mathfrak n(k,x_l)} + ia_1 \right) \sum_{\ell=0}^n \left( \mathcal C_{n-\ell}^{(x_l,x)}(k)  -  (-1)^n i\mathcal S_{n-\ell}^{(x_l,x)}(k) \right) {\mathcal E}_\ell^{(x,\infty)}(k).
    \end{align}
    Using Lemma~\ref{lem:Delta_identity}, 
    \begin{align} 
        \chi(k,x) 
        &= - 4k (\beta\mathfrak n)(k,x) \sum_{n=0}^\infty \left( \frac{(-1)^n a_0}{k \mathfrak n(k,x_l)} + ia_1 \right) \mathcal E_{n}^{(x_l,\infty)}(k) = 2ik (\beta\mathfrak n)(k,x) \Delta(k).
    \end{align}
    Finally, for the finite-interval problem, since $\mathcal C_n^{(x,x)}(k) = \delta_{0n}$ and $\mathcal S_n^{(x,x)}(k) = 0$,
    \begin{align} \nonumber
        \chi(k,x) 
        &= -4k (\beta\mathfrak n)(k,x) {{ \Eint}} \left\{ \mathfrak a(k) +  \sum_{n=0}^\infty \mathfrak c_n(k) \sum_{\ell=0}^n \left( \mathcal C_{n-\ell}^{(x_l,x)}(k) \mathcal C_\ell^{(x,x_r)}(k) -  (-1)^{n-\ell} \mathcal S_{n-\ell}^{(x_l,x)}(k) \mathcal S_\ell^{(x,x_r)}(k) \right) \right. \\ 
        &\hspace*{1.4in} \left. + \sum_{n=0}^\infty  \mathfrak s_n(k) \sum_{\ell=0}^n \left( \mathcal S_{n-\ell}^{(x_l,x)}(k) \mathcal C_\ell^{(x,x_r)}(k) + (-1)^{n-\ell} \mathcal C_{n-\ell}^{(x_l,x)}(k)  \mathcal S_\ell^{(x,x_r)}(k) \right) \right\}\!, 
    \end{align}
    which gives \eqref{eqn:chi_identity}, using Lemma~\ref{lem:Delta_identity}.
\end{proof}
\begin{lemma} \label{lem:Bmx}
    For the half-line problem, 
    \begin{subequations} \label{eqn:Bmx}
        \begin{align}
            \label{eqn:B0x_hl}
            \mathcal B_{0,x}(k,x) 
            &= \frac{4\beta(x_l) ik \mathfrak n(k,x) \exp\left( \int_{x_l}^{x} ik \mathfrak n(k,\xi) \, d\xi\right)}{\sqrt{(\beta\mathfrak n)(k,x_l)}\sqrt{(\beta \mathfrak n)(k,x)}} \sum_{n=0}^\infty \mathcal E_n^{(x,\infty)}(k),
        \end{align}
        and there exists an $r>\sqrt{M_\gamma}$ and $M_{\mathcal B} >0$ so that for $k\in \Omegaext(r)$ and $x\in \bar \D$,
        \begin{align} \label{eqn:Bmx/Delta_bound_hl}
            \left| \mathcal B_{0,x}(k,x) \right| \leq M_{\mathcal B}|k|
            \qquad \text{ and } \qquad
            \left| \frac{\mathcal B_{0,x}(k,x)}{\Delta(k)} \right| 
            &\leq  M_{\mathcal B} |k|^2 e^{-m_{i\mathfrak n}|k|(x-x_l)}. 
        \end{align}
    \end{subequations}
    \begin{subequations}
        For the finite-interval problem, we have for $j=1,2$,
        \begin{align}
            \label{eqn:Bmx_fi} \nonumber
            \mathcal B_{2-j,x}(k,x) 
            &= -(-1)^{j}\frac{4k \mathfrak n(k,x) \Eint}{\sqrt{(\beta\mathfrak n)(k,x)}} \left \{ \frac{\beta(x_r)}{\sqrt{(\beta\mathfrak n)(k,x_r)}}  \left[ \frac{a_{j1}}{k \mathfrak n(k,x_l)} \sum_{n=0}^\infty (-1)^{n}\mathcal C_{n}^{(x_l,x)}(k) + a_{j2} \sum_{n=0}^\infty (-1)^{n} \mathcal S_{n}^{(x_l,x)}(k) \right] \right.  \\ 
            & \hspace*{1.5in} \left. + \frac{\beta(x_l)}{\sqrt{(\beta\mathfrak n)(k,x_l)}}  \left[  \frac{b_{j1}}{k \mathfrak n(k,x_r)} \sum_{n=0}^\infty \mathcal C_{n}^{(x,x_r)}(k) - b_{j2} \sum_{n=0}^\infty (-1)^n \mathcal S_{n}^{(x,x_r)}(k) \right]\right\}\!,
        \end{align}
        and there exists an $r>\sqrt{M_\gamma}$ and $M_{\mathcal B} >0$ so that for $k\in \Omegaext(r)$ and $x\in \bar \D$,
        \begin{align} \label{eqn:Bmx/Delta_bound_fi}
            \left| \mathcal B_{m,x}(k,x) \right| \leq M_{\mathcal B}|k|
            \qquad \text{ and } \qquad
            \left| \frac{\mathcal B_{m,x}(k,x)}{\Delta(k)} \right| \leq M_{\mathcal B} |k|^{b+1} \big( e^{-m_{i\mathfrak n}|k|(x_r-x)} + e^{-m_{i\mathfrak n}|k|(x-x_l)} \big). 
        \end{align}
        For {\em regular} boundary conditions $b=1$, and for {\em irregular} boundary conditions $b=2$. Therefore, the functions $\mathcal B_{m,x}(x,t)$ and $\mathcal B_{m,x}(x,t)/\Delta(k)$ are well defined for the half-line and finite-interval problems.
    \end{subequations}
\end{lemma}
\begin{proof}
    Lemma~\ref{lem:Jnderivatives} and a direct calculation gives \eqref{eqn:B0x_hl} and \eqref{eqn:Bmx_fi}. The proofs for \eqref{eqn:Bmx/Delta_bound_hl} and \eqref{eqn:Bmx/Delta_bound_fi} are identical to the proof of Lemma~\ref{lem:Bm_bounds}. Note that, as in Lemma~\ref{lem:Psibar_bounds}, the $M_{\mathcal B}$'s differ only by a factor of $M_{\mathfrak n}$ (see Lemma~\ref{lem:nubounds}). Without loss of generality, we may take them to be identical.
\end{proof}
\begin{lemma} \label{lem:Phi0x}
    Consider the finite-interval, half-line, and whole-line problems. We have
    \begin{align} 
        \label{eqn:Phi0x}
        \Phi_{0,x}(k,x) &= \int_{\D} \frac{\bar \Psi(k,x,y)q_\alpha (y)}{ \sqrt{(\beta \mathfrak n)(k,x)} \sqrt{(\beta \mathfrak n)(k,y)}} \, dy,
    \end{align}
    where $\Phi_0(k,x)$ is defined in \eqref{eqn:Phi_0}. There exists an $M_\Phi>0$ so that
    \begin{align}\label{eqn:Phi0x_bounds}
        \left| \Phi_{0,x}(k,x) \right| \leq M_{\Phi} |k|  \|q_0\|_\D
        \qquad \text{ and } \qquad 
        \left| \frac{\Phi_{0,x}(k,x)}{\Delta(k)} \right| \leq M_{\Phi} |k|^2 \|q_0\|_\D . 
    \end{align}
    Thus $\Phi_{0,x}(k,x)$ and $\Phi_{0,x}(k,x)/\Delta(k)$ are well defined for all three problems.
\end{lemma}
\begin{proof}
    Breaking up the integral over $\D$ in \eqref{eqn:Phi_0} into two integrals over the regions $y<x$ and $y>x$ and using the Leibniz integral rule, we obtain
    \begin{align}
        \Phi_{0,x}(k,x) &= \frac{(\Psi(k,x,x^-) - \Psi(k,x,x^+))q_\alpha(x)}{(\beta\mathfrak n)(k,x)} + \int_{\D} \frac{\bar \Psi(k,x,y)q_\alpha(y)}{ \sqrt{(\beta \mathfrak n)(k,x)}\sqrt{(\beta \mathfrak n)(k,y)}} \, dy. 
    \end{align}
    Since $\Psi(k,x,x^-) = \Psi(k,x,x^+)$, we find \eqref{eqn:Phi0x}.
    We obtain \eqref{eqn:Phi0x_bounds} from Lemma~\ref{lem:Psibar_bounds}. 
    Since the integrand in \eqref{eqn:Phi0x} is absolutely integrable, differentiation under the integral is allowed.
\end{proof}
\begin{lemma} \label{lem:Phifx}
    Consider the finite-interval, half-line, and whole-line problems. For $k\in \Omegaext$, $x\in \bar \D$, and $t\in(0,T)$,
    \begin{align} 
        \label{eqn:Phifx}
        \Phi_{\mathfrak f,x}(k,x,t) 
        &= \int_{\D} \frac{\bar \Psi(k,x,y)\mathfrak f_\alpha (k^2,y,t)} { \sqrt{(\beta\mathfrak n)(k,x)} \sqrt{(\beta\mathfrak n)(k,y)}} \,  dy.
    \end{align}
    Further, there exists an $M_f>0$ so that
    \begin{subequations} \label{eqn:Phifx_bounds}
        \begin{align} 
            \label{eqn:Phifx_bound_direct}
            \big|\Phi_{\mathfrak f,x}(k,x,t) e^{-k^2t} \big|
            &\leq M_f |k| \int_{\D} \big| \mathfrak f_\alpha (k^2,y,t) e^{-k^2t} \big| \,  dy. 
        \end{align}
        For the {\em regular problems},
        \begin{align} 
            \label{eqn:Phifx/Delta_bound_reg}
            \left| \frac{\Phi_{\mathfrak f,x}(k,x,t)e^{-k^2t}}{\Delta(k)} \right|
            &\leq M_f |k| \int_{\D} \big| \mathfrak f_\alpha (k^2,y,t) e^{-k^2t} \big| \,  dy,
        \end{align}
        and for the {\em irregular problems},
        \begin{align} 
            \label{eqn:Phifx/Delta_bound_irreg}
            \left| \frac{\Phi_{\mathfrak f}(k,x,t)e^{-k^2t}}{\Delta(k)} \right|
            \leq M_f |k| \Big( 1 + |k| \big( e^{-m_{i\mathfrak n}|k|(x-x_l)} + e^{-m_{i\mathfrak n}|k|(x_r-x)} \big) \Big) \int_{\D} \big| \mathfrak f_\alpha(k^2,y,t) e^{-k^2t} \big| \, dy.
        \end{align}
    \end{subequations}
    where $\Phi_{\mathfrak f}(k,x,t)$ and $\mathfrak f_\alpha(k^2,y,t)$ are defined in \eqref{eqn:Phif_frak} and \eqref{eqn:f_frak}, respectively.
\end{lemma}
\begin{proof}
    Breaking up the integral over $\D$ in \eqref{eqn:Phi_0} into two integrals over the regions $y<x$ and $y>x$
    and using the Leibniz integral rule, we obtain
    \begin{align}
        \Phi_{\mathfrak f,x}(k,x,t) 
        &=  \frac{\left(\Psi(k,x,x^-) - \Psi(k,x,x^+)\right) \mathfrak f_\alpha(k^2,x,t)} {\sqrt{(\beta \mathfrak n)(k,x)} \sqrt{(\beta\mathfrak n)(k,x)}} +  \int_{\D} \frac{\bar \Psi(k,x,y) \mathfrak f_\alpha(k^2,y,t)} {\sqrt{(\beta\mathfrak n)(k,x)} \sqrt{(\beta\mathfrak n)(k,y)}}  \, dy. 
    \end{align}
    Since $\Psi(k,x,x^-) = \Psi(k,x,x^+)$ for all three problems, we obtain \eqref{eqn:Phifx}. Equation \eqref{eqn:Phifx_bounds} follows from \eqref{eqn:Psibar_bounds}. 
    Since the integrand \eqref{eqn:Phifx} is absolutely integrable, differentiation under the integral is allowed.
\end{proof}
\begin{lemma} \label{lem:qx}
    Consider the finite-interval, half-line, and whole-line problems. For $x\in  \D$ and $t\in(0,T)$,
    \begin{subequations}
        \begin{align} 
            \label{eqn:q0x}
            q_{0,x}(x,t) &= \frac{1}{2\pi} \int_{\partial\Omegaext} \frac{\Phi_{0,x}(k,x)}{\Delta(k)} e^{-k^2t} \, dk, \\
            \label{eqn:qfx}
            q_{f,x}(x,t) 
            &= \frac{1}{2\pi} \int_{\partial\Omegaext} \frac{\Phi_{\mathfrak f,x}(k,x,t)e^{-k^2t}}{\Delta(k)}  dy, \\
            \label{eqn:q_Bmx}
            q_{\mathcal B_m,x}(x,t) 
            &= \frac{1}{2\pi} \int_{\partial\Omegaext}  \frac{\mathcal B_{m,x} (k,x)} {\Delta(k)} \mathfrak F_m(k^2,t) e^{-k^2t}  \, dk,
        \end{align}
    \end{subequations}
    are well defined, \ie we can differentiate under the integral sign. Furthermore, $q_{0,x}(x,t)$ and $q_{f,x}(x,t)$ are well defined for $x\in \bar\D$. For the {\em regular problems}, $q_{\mathcal B_m,x}(x,t)$ is well defined for $x\in\bar \D$.
\end{lemma}
\begin{proof}
    The integrand in $q_{0,x}(x,t)$ is exponentially decaying for $t\in(0,T)$, and therefore is well defined for $x\in \bar \D$. From \eqref{eqn:Phifx_bounds} and \eqref{eqn:f_frak_bounds}, we see that, for any $t\in (0,T)$, $q_{f,x}(x,t)$ is also well defined for $x\in\bar\D$. For $t\in(0,T)$, from \eqref{eqn:Bmx/Delta_bound_hl}, \eqref{eqn:Bmx/Delta_bound_fi}, and \eqref{eqn:Fmfrak_bound}, we see that for $x\in \D$, $q_{\mathcal B_m,x}(x,t)$ has exponential decay and is well defined. For the {\em regular problems}, $q_{\mathcal B_m,x}(x,t)$ is absolutely integrable for $x\in \bar \D$ and is well defined. 
\end{proof}
\begin{remark}
    For the {\em irregular problems}, $q_{\mathcal B_m,x}(x,t)$ may be ill defined at the boundaries, but the boundary conditions \eqref{eqn:BC1_fi} and \eqref{eqn:BC2_fi} are well defined and satisfied, see Section~\ref{sec:proofs_BC}.
\end{remark}
\begin{lemma} \label{lem:Psitilde}
    Consider the finite-interval, half-line, and whole-line problems. For $x,y\in \D$ and $k\in \Omegaext$,
    \begin{subequations} \label{eqn:Psixx}
        \begin{align} \label{eqn:Psitilde}
            \tilde \Psi(k,x,y) &= -\frac{k^2+\gamma(x)}{\alpha(x)} \Psi(k,x,y).
        \end{align}
        For the half-line $(m=0)$ and the finite-interval problems $(m=0,1)$,
        \begin{align} \label{eqn:Bmxx}
            (\beta\mathcal B_{m,x})_x(k,x) = - \frac{k^2+\gamma(x)}{\alpha(x)} \mathcal B_m(k,x),
        \end{align}
        for $x\in \D$, $t\in(0,T)$, and $k\in\Omega$.
    \end{subequations}
\end{lemma}
\begin{proof}
    For the whole-line problem, a direct calculation using Lemma~\ref{lem:Jnderivatives} gives \eqref{eqn:Psitilde}  from \eqref{eqn:Psibar_wl_lower} for $y<x$ and from \eqref{eqn:Psibar_wl_upper} for $y>x$. 
    Similarly, for the half-line problem, we obtain \eqref{eqn:Psitilde} from \eqref{eqn:Psibar_hl_lower} for $x_l<y<x$ and from \eqref{eqn:Psibar_hl_upper} for $x_l<x<y$. Equation \eqref{eqn:Bmxx} follows from \eqref{eqn:B0x_hl}.
    For the finite-interval problem, we obtain \eqref{eqn:Psitilde} from \eqref{eqn:Psibar_fi_lower} for $x_l<y<x<x_r$ and from \eqref{eqn:Psibar_fi_upper} for $x_l<x<y<x_r$. Finally, \eqref{eqn:Bmxx} follows from \eqref{eqn:Bmx_fi}.
\end{proof}
\begin{lemma} \label{lem:Phixx}
    Consider the finite-interval, half-line, and whole-line problems. With $\mathfrak f(k^2,x,t) = \alpha(x) \mathfrak f_\alpha(k^2,x,t)$,
    \begin{subequations} \label{eqn:Phixx_wl}
        \begin{align} 
            \label{eqn:Phi0xx_wl}
            \alpha(x) (\beta\Phi_{0,x})_x(k,x) &= 2ik\Delta(k) q_0(x) - (k^2+\gamma(x)) \Phi_0(k,x), \\
            \label{eqn:Phifxx_wl}
            \alpha(x) (\beta\Phi_{\mathfrak f,x})_x(k,x,t) 
            &= 2ik \Delta(k) \mathfrak f(k^2,x,t) - (k^2+\gamma(x)) \Phi_{\mathfrak f}(k,x,t),
        \end{align}
    \end{subequations}
\end{lemma}
\begin{proof}
    \begin{subequations}
        Using Lemmas~\ref{lem:Phi0x}~and~\ref{lem:Phifx}, we split $\D$ into the two parts $y<x$ and $y>x$, and the Leibniz integral rule gives
        \begin{align}
            (\beta\Phi_{0,x})_x(k,x) 
            &= \frac{\chi(k,x) q_0(x)}{\alpha(x) (\beta\mathfrak n)(k,x)} + \int_{\D} \frac{\tilde \Psi(k,x,y)q_0(y)}{\alpha(y)\sqrt{(\beta \mathfrak n)(k,x)} \sqrt{(\beta\mathfrak n)(k,y)}} \, dy, \\
            (\beta\Phi_{\mathfrak f,x})_x(k,x,t) 
            &= \frac{\chi(k,x) \mathfrak f(k^2,x,t)} { \alpha(x) (\beta \mathfrak n)(k,x)} + \int_{\D} \frac{\tilde \Psi(k,x,y)\mathfrak f(k^2,y,t)} { \alpha(y) \sqrt{(\beta \mathfrak n)(k,x)} \sqrt{(\beta\mathfrak n)(k,y)}} \,  dy,
        \end{align}
        where $\chi(k,x)$ and $\tilde \Psi(k,x,y)$ are defined in Definition~\ref{def:Psibar_tilde}. Using Lemmas~\ref{lem:delta(k,x)}~and~\ref{lem:Psitilde} gives \eqref{eqn:Phixx_wl}.
    \end{subequations}
\end{proof}
\begin{lemma} \label{lem:qt}
    Consider the finite-interval, half-line, and whole-line problem. For $x\in \D$ and $t\in(0,T)$, the $t$-derivatives of $q_0(x,t)$, $q_f(x,t)$ and $q_{\mathcal B_m}(x,t)$ are
    \begin{subequations} \label{eqn:qt}
        \begin{align}
            \label{eqn:q0t}
            q_{0,t}(x,t) &= -\frac{1}{2\pi} \int_{\partial\Omegaext} \frac{k^2\Phi_{0}(k,x)}{\Delta(k)} e^{-k^2t} \, dk, \\
            \label{eqn:qft}
            q_{f,t}(x,t) &= - \frac{1}{2\pi} \int_{\partial\Omegaext(r)} \frac{ k^2 \Phi_{\mathfrak f}(k,x,t)e^{-k^2t}}{\Delta(k)}  dk, \\
            \label{eqn:qBmt}
            q_{\mathcal B_m,t}(x,t) &= - \frac{1}{2\pi} \int_{\partial\Omegaext(r)}  \frac{k^2\mathcal B_m(k,x)} {\Delta(k)}  \mathfrak F_m(k^2,t) e^{-k^2t} \, dk, \qquad m=0,1.
        \end{align}
    \end{subequations}
    These functions are well defined.
\end{lemma}
\begin{proof}
    Differentiating \eqref{eqn:q_0_deformed} with respect to $t$ gives \eqref{eqn:q0t}, since the integrand is absolutely integrable. 
    From \eqref{eqn:f_frak}, $\mathfrak f_{\alpha,t}(k^2,x,t)e^{-k^2t} =  - f_{\alpha,t}(x,t)/k^2$, and differentiating \eqref{eqn:Phif_frak} with respect to $t$ yields
    \begin{align} 
        \label{eqn:Phift}
        \Phi_{\mathfrak f,t}(k,x,t) e^{-k^2t} 
        = - \int_{\D} \frac{\Psi(k,x,y)f_{\alpha,t}(y,t)} {k^2 \sqrt{(\beta \mathfrak n)(k,x)} \sqrt{(\beta \mathfrak n)(k,y)}}  \, dy, 
    \end{align}
    so that, using Lemma~\ref{lem:Psi_bounds},
    \begin{align} \label{eqn:Phi_ft_bounds}
        \left| \frac{\Phi_{\mathfrak f,t}(k,x,t)e^{-k^2t}}{\Delta(k)} \right|
        \leq \frac{M_{f}}{|k|^2} \Big(  1 + |k| \big( \eone{-m_{i\mathfrak n} |k|(x-x_l)} + \eone{- m_{i\mathfrak n}|k|(x_r-x)} \big) \Big)  \|f_{\alpha,t}\|_\D.
    \end{align}
    Differentiating \eqref{eqn:q_f_deformed2_direct} with respect to $t$, we obtain
    \begin{align} 
        q_{f,t}(x,t) 
        &= \frac{1}{2\pi} \int_{\partial\Omegaext(r)} \frac{\Phi_{\mathfrak f,t}(k,x,t)e^{-k^2t} }{\Delta(k)}  dk - \frac{1}{2\pi} \int_{\partial\Omegaext(r)} \frac{ k^2 \Phi_{\mathfrak f}(k,x,t)e^{-k^2t}}{\Delta(k)}  dk.
    \end{align}
    From \eqref{eqn:Phi_ft_bounds}, it follows that the first contour integral can be closed in the upper half plane, implying it is zero by Cauchy's theorem, resulting in \eqref{eqn:qft}. From \eqref{eqn:Phif_frak_bounds},  $q_{f,t}(x,t)$ is well defined, for $x\in \D$.
    Since $\mathfrak F_{m,t}(k^2,t)e^{-k^2t} = - f_m'(t)/k^2$, differentiating \eqref{eqn:q_Bm_deformed3} with respect to $t$, 
    \begin{align} 
        q_{\mathcal B_m,t}(x,t) 
        &= - \frac{f_m'(t)}{2\pi} \int_{\partial\Omegaext(r)}  \frac{\mathcal B_m(k,x)} {k^2\Delta(k)} \, dk - \frac{1}{2\pi} \int_{\partial\Omegaext(r)}  \frac{k^2\mathcal B_m(k,x)} {\Delta(k)}  \mathfrak F_m(k^2,t) e^{-k^2t} \, dk.
    \end{align}
    As above, \eqref{eqn:Bm_bounds} allows us to close the contour of the first integral in the upper half plane, showing the first term is zero by Cauchy's theorem, obtaining \eqref{eqn:qBmt}. From \eqref{eqn:Bm_bounds} and \eqref{eqn:Fmfrak_bound},  $q_{\mathcal B_m,t}(x,t)$ is well defined for $x\in \D$.
\end{proof}
\begin{lemma} \label{lem:qxx}
    For $x\in \D$ and $t\in(0,T)$, the derivatives
    \begin{subequations} \label{eqn:qxx}
        \begin{align} 
            \label{eqn:q0xx}
            \alpha(x)(\beta q_{0,x})_x(x,t) + \gamma(x) q_0(x,t) &=  q_{0,t}(x,t), \\
            \label{eqn:qfxx}
            \alpha(x)(\beta q_{f,x})_x(x,t) + \gamma(x) q_f(x,t) + f(x,t) &=q_{f,t}(x,t), \\
            \label{eqn:qBmxx}
            \alpha(x) (\beta q_{\mathcal B_{m},x})_x(x,t) + \gamma(x) q_{\mathcal B_m}(x,t) &= q_{\mathcal B_m,t}(x,t), \qquad m=0,1.
        \end{align}
    \end{subequations}
    are well defined, \ie differentiation under the integral sign is allowed.
\end{lemma}
\begin{proof}
    Direct differentiation of the results in Lemma~\ref{lem:qx} yields 
    \begin{subequations} 
        \begin{align} 
            (\beta q_{0,x})_x(x,t) &= \frac{1}{2\pi} \int_{\partial\Omegaext} \frac{(\beta\Phi_{0,x})_x(k,x)}{\Delta(k)} e^{-k^2t} \, dk, \\
            (\beta q_{f,x})_x(x,t) &= \frac{1}{2\pi} \int_{\partial\Omegaext} \frac{(\beta \Phi_{\mathfrak f,x})_x(k,x,t)e^{-k^2t}}{\Delta(k)}  dk, \\
            (\beta q_{\mathcal B_{m},x})_x(x,t) &= \frac{1}{2\pi} \int_{\partial\Omega} \frac{(\beta\mathcal B_{m,x})_x(k,x)}{\Delta(k)} \mathfrak F_m(k^2,t) e^{-k^2t} \, dk, \qquad m=0,1.
        \end{align}
    \end{subequations}
    Using Lemmas~\ref{lem:Psitilde},~\ref{lem:Phixx},~and~\ref{lem:qt}, 
    \begin{subequations} \label{eqn:qxx_intermediate}
        \begin{align} 
            \label{eqn:q0xx_intermediate}
            \alpha(x)(\beta q_{0,x})_x(x,t) + \gamma(x) q_0(x,t)
            &= q_{0,t}(x,t) -\frac{q_0(x)}{i\pi} \int_{\partial\Omegaext} k e^{-k^2t} \, dk , \\
            \label{eqn:qfxx_intermediate}
            \alpha(x)(\beta q_{f,x})_x(x,t) + \gamma(x) q_f(x,t) &=q_{f,t}(x,t) - \frac{1}{i\pi} \int_{\partial\Omegaext} k \mathfrak f(k^2,x,t)  e^{-k^2t} dk, \\
            \label{eqn:qBmxx_intermediate}
            \alpha(x) (\beta q_{\mathcal B_{m},x})_x(x,t) + \gamma(x) q_{\mathcal B_m}(x,t) &= q_{\mathcal B_m,t}(x,t), \qquad m=0,1.
        \end{align}
    \end{subequations}
    Since the integrands in \eqref{eqn:qxx_intermediate} are absolutely integrable, the differentiation inside the integral is justified. The path for the remaining integral in \eqref{eqn:q0xx_intermediate} can be deformed down to the real line showing it is zero. Using \eqref{eqn:f_frak}, the remaining integral in \eqref{eqn:qfxx_intermediate} is evaluated as
    \begin{align}
        -\int_{\partial\Omegaext} k \mathfrak f(k^2,x,t)  e^{-k^2t} dk 
        &= \int_{\partial\Omegaext} \left( \frac{f(y,0)}{k} + \frac{1}{k} \int_0^t f_{s}(y,s) e^{k^2s} \, ds \right)  e^{-k^2t} dk,
    \end{align}
    which may also be deformed to an indented contour on the real line. The principal-value part integral is zero, while the indentation integral evaluates to
    \begin{align}
        \frac{1}{i\pi} \int_{\partial\Omegaext} k \mathfrak f(k^2,x,t)  e^{-k^2t} dk 
        &= \Res\left( \left( \frac{f(y,0)}{k} + \frac{1}{k} \int_0^t f_{s}(y,s) e^{k^2s} \, ds \right)  e^{-k^2t}; \, k=0 \right) 
        = f(y,t).
    \end{align}
    Equation \eqref{eqn:qxx_intermediate} yields \eqref{eqn:qxx}.
\end{proof}
\begin{theorem} \label{thm:evolution_equation}
    The solution expressions \eqref{eqn:q_wl}, \eqref{eqn:q_hl}, and \eqref{eqn:q_fi} each solve the evolution equation \eqref{eqn:evolution_equation}.
\end{theorem}
\begin{proof}
    Since $q(x,t) = q_0(x,t) + q_{f}(x,t) + q_{\mathcal B_0}(x,t) + q_{\mathcal B_1}(x,t)$, \eqref{eqn:qxx} gives the result.
\end{proof}
\section{Proofs: the solution expressions satisfy the boundary values} 
%

%
\newcommand\StepSubequations{
  \stepcounter{parentequation}
  \gdef\theparentequation{D.\arabic{parentequation}}
  \setcounter{equation}{0}
}
\newcommand{\newlink}[2]{$\mathrm{(}$\hyperlink{#1}{$\mathrm{#2}$}$\mathrm{)}$}
\renewcommand{\sp}{\hspace*{0.5in}}
\begin{definition} \label{def:boundary_functions}
    In this appendix, $\ell=0$ corresponds to the half-line problem, while $\ell=1,2$ correspond to the finite-interval problem. We define, 
    for $k\in\Omegaext$ and $y \in \D$,
    \setlength{\belowdisplayskip}{10pt}%
    \setlength{\abovedisplayskip}{10pt}%
    \begin{subequations} \label{eqn:mathfrakPl}
        \begin{align}
            \hypertarget{eqnmathfrakPl}{}
            \xdef \eqnmathfrakPl {\theequation}
            \label{eqn:mathfrakPl_hl}
            \mathfrak P^{(0)}(k,y)
            &= \frac{a_{0}\Psi(k,x_l,y) + a_{1} \bar \Psi(k,x_l,y)}{\sqrt{(\beta \mathfrak n)(k,x_l)}}, \\
            \label{eqn:mathfrakPl_fi}
            \mathfrak P^{(\ell)}(k,y) 
            &= \frac{a_{\ell1}\Psi(k,x_l,y) + a_{\ell2} \bar \Psi(k,x_l,y)}{\sqrt{(\beta \mathfrak n)(k,x_l)}}+ \frac{b_{\ell1}\Psi(k,x_r,y)+ b_{\ell2} \bar \Psi(k,x_r,y)}{\sqrt{(\beta \mathfrak n)(k,x_r)}}, 
            && \ell=1,2.
            \intertext{For $k\in\Omegaext$,}
            \StepSubequations
            \hypertarget{eqnmathfrakBml}{}
            \xdef \eqnmathfrakBml {\theequation}
            \label{eqn:mathfrakBml_hl}
            \mathfrak B_0^{(0)}(k) 
            &= a_{0} \mathcal B_0(k,x_l) + a_{1} \mathcal B_{0,x}(k,x_l), \\
            \label{eqn:mathfrakBml_fi}
            \mathfrak B_m^{(\ell)}(k) 
            &= a_{\ell1} \mathcal B_m(k,x_l) + a_{\ell2} \mathcal B_{m,x}(k,x_l) + b_{\ell1} \mathcal B_m(k,x_r) + b_{\ell2} \mathcal B_{m,x}(k,x_r) ,
            && \ell=1,2,~~m=0,1, \\
            \intertext{and}
            \StepSubequations
            \hypertarget{eqnmathcalP0}{}
            \xdef \eqnmathcalPnaught {\theequation}
            \label{eqn:mathcalP0_hl}
            \mathcal P_0^{(0)}(k) &= a_{0} \Phi_{0}(k,x_l)+ a_{1} \Phi_{0,x}(k,x_l),\\
            \label{eqn:mathcalP0_fi}
            \mathcal P_0^{(\ell)}(k) &= a_{\ell1} \Phi_{0}(k,x_l)+ a_{\ell2} \Phi_{0,x}(k,x_l)+b_{\ell1} \Phi_{0}(k,x_r)+b_{\ell2} \Phi_{0,x}(k,x_r), && \ell=1,2. \\
            \intertext{For $k\in\Omegaext$ and $t\in(0,T)$,}
            \StepSubequations 
            \hypertarget{eqnmathcalPfl}{}
            \xdef \eqnmathcalPfl {\theequation}
            \label{eqn:mathcalPfl_hl}
            \mathcal P_f^{(0)}(k,t) &= a_{0} \Phi_{\mathfrak f}(k,x_l,t)+ a_{1} \Phi_{\mathfrak f,x}(k,x_l,t), \\
            \label{eqn:mathcalPfl_fi}
            \mathcal P_f^{(\ell)}(k,t) &= a_{\ell1} \Phi_{\mathfrak f}(k,x_l,t)+ a_{\ell2} \Phi_{\mathfrak f,x}(k,x_l,t)+b_{\ell1} \Phi_{\mathfrak f}(k,x_r,t)+b_{\ell2} \Phi_{\mathfrak f,x}(k,x_r,t), && \ell=1,2. \\
            \intertext{Finally, for $t\in (0,T)$,}
            \StepSubequations 
            \hypertarget{eqnmathcalQBm}{}
            \xdef \eqnmathcalQBm {\theequation}
            \label{eqn:mathcalQBm_hl}
            \mathcal Q_{\mathcal B_m}^{(0)} (t) &= a_{0} q_{\mathcal B_m}(x_l,t) +a_{1} q_{\mathcal B_m,x}(x_l,t) , \\
            \label{eqn:mathcalQBm_fi}
            \mathcal Q_{\mathcal B_m}^{(\ell)} (t) &= a_{\ell1} q_{\mathcal B_m}(x_l,t) +a_{\ell2} q_{\mathcal B_m,x}(x_l,t) + b_{\ell1} q_{\mathcal B_m}(x_r,t) + b_{\ell2} q_{\mathcal B_m,x}(x_r,t), && \ell=1,2, \\
            \intertext{and}
            \StepSubequations
            \hypertarget{eqnmathcalQ0}{ }
            \xdef \eqnmathcalQnaught {\theequation}
            \label{eqn:mathcalQ0_hl}
            \mathcal Q_{0}^{(0)} (t) &= a_{0} q_{0}(x_l,t) + a_{1} q_{0,x}(x_l,t),\\
            \label{eqn:mathcalQ0_fi}
            \mathcal Q_{0}^{(\ell)} (t) &= a_{\ell1} q_{0}(x_l,t) +a_{\ell2} q_{0,x}(x_l,t) + b_{\ell1} q_{0}(x_r,t) + b_{\ell2} q_{0,x}(x_r,t), && \ell=1,2, \\
            \intertext{and}
            \StepSubequations
            \hypertarget{eqnmathcalQf}{}
            \xdef \eqnmathcalQf {\theequation}
            \label{eqn:mathcalQf_hl}
            \mathcal Q_{f}^{(0)} (t) &= a_{0} q_{f}(x_l,t) +a_{1} q_{f,x}(x_l,t), \\
            \label{eqn:mathcalQf_fi}
            \mathcal Q_{f}^{(\ell)} (t) &= a_{\ell1} q_{f}(x_l,t) +a_{\ell2} q_{f,x}(x_l,t) + b_{\ell1} q_{f}(x_r,t) + b_{\ell2} q_{f,x}(x_r,t), && \ell=1,2.
        \end{align}
    \end{subequations}
\end{definition}
\begin{lemma} \label{lem:mathfrakPl=0}
    For both the half-line problem and the finite-interval problem, for $k\in \Omegaext$ and $y\in \bar\D$, 
    \begin{align} \label{eqn:mathfrakPl=0}
        \mathfrak P^{(\ell)}(k,y) &= 0, ~~\ell=0,1,2.
    \end{align}
\end{lemma}
\def \sp {0.3in}
\begin{proof}
    For the half-line, using \eqref{eqn:Psi_hl} (with $x_l=x<y<x_r$) and \eqref{eqn:Psibar_hl_upper} in \eqref{eqn:mathfrakPl_hl}, gives \eqref{eqn:mathfrakPl=0}.
    \begin{align} \nonumber
        \Psi(k,x_l,y) 
        &= -4\exp\left(\int_{x_l}^{y}ik \mathfrak n(k,\xi)\,d\xi \right) \sum_{n=0}^\infty(-1)^n a_1 {\mathcal E}_n^{(y,\infty)}(k), \\ \nonumber
        \bar \Psi(k,x_l,y) 
        &= 4 \exp\left(\int_{x_l}^{y}ik\mathfrak n(k,\xi)\, d\xi \right) \sum_{n=0}^\infty(-1)^n a_0 {\mathcal E}_n^{(y,\infty)}(k), \\
        \Rightarrow 
        \hspace*{\sp} \mathfrak P^{(0)}(k,y) &= \frac{a_{0}\Psi(k,x_l,y) + a_{1} \bar \Psi(k,x_l,y)}{\sqrt{(\beta\mathfrak n)(k,x_l)}} = 0.
    \end{align}
    Using \eqref{eqn:Psi_fi} and \eqref{eqn:Psibar_fi} in \eqref{eqn:mathfrakPl_fi}, the calculations for the finite-interval case are equally straightforward albeit more tedious. 
\end{proof}
\begin{lemma} \label{lem:mathfrakBml_eqn}
    For the half-line problem $(m=0)$,  and for the finite-interval problem $(m=0,1)$, for $k\in \Omegaext$, 
    \begin{subequations} \label{eqn:mathfrakBml_eqns}
        \begin{align} \label{eqn:mathfrakBml_eqn}
            \mathfrak B_m^{(\ell)}(k) &= -2ik \Delta(k) \tilde \delta_{\ell-1,m},~~~\ell=0,1,2.
        \end{align}
        Here
        \begin{align} \label{eqn:deltatilde}
            \tilde \delta_{\ell-1,m} = \case{1}{1, & \ell = 0, \, m=0, \\ 1, & \ell\neq 0, \, m=\ell-1, \\ 0, & \ell\neq 0,\, m\neq \ell-1.}
        \end{align}
    \end{subequations}
\end{lemma}
\begin{proof}
    For the half-line problem, using \eqref{eqn:B0_hl} and \eqref{eqn:B0x_hl} in \eqref{eqn:mathfrakBml_hl}, we find \eqref{eqn:mathfrakBml_eqns}:
    \begin{align} \nonumber
        \mathcal B_0(k,x_l) 
        &=  \frac{4}{\mathfrak n(k,x_l)} \sum_{n=0}^\infty (-1)^n {\mathcal E}_{n}^{(x_l,\infty)}(k), \\ \nonumber
        \mathcal B_{0,x}(k,x_l) 
        &= \frac{4 ik \mathfrak n(k,x_l) }{\mathfrak n(k,x_l)} \sum_{n=0}^\infty \mathcal E_n^{(x_l,\infty)}(k), \\
        \Rightarrow 
        \hspace*{\sp} 
        \mathfrak B_0^{(0)}(k) &= a_{0} \mathcal B_0(k,x_l) + a_{1} \mathcal B_{0,x}(k,x_l) = 4 \sum_{n=0}^\infty \left(\frac{(-1)^n a_0}{\mathfrak n(k,x_l)}  +   a_1 ik \right) \mathcal E_n^{(x_l,\infty)}(k) = -2i k \Delta(k).
    \end{align}
    The finite-interval case (using \eqref{eqn:Bm} and \eqref{eqn:Bmx_fi} in \eqref{eqn:mathfrakBml_fi}) is similar but more tedious. Its details are omitted.
\end{proof}
\begin{lemma} \label{lem:Phi_BCs}
    For both the half-line problem and the finite-interval problem, for $k\in\Omegaext$ and $t\in [0,T]$,
    \begin{subequations} \label{eqn:Phi_BCs}
        \begin{alignat}{2} 
            \label{eqn:Phi0_BCs}
            \mathcal P_0^{(\ell)}(k) &= 0, \\
            \label{eqn:Phif_BCs}
            \mathcal P_f^{(\ell)}(k,t) &= 0.
        \end{alignat}
    \end{subequations}
\end{lemma}
\begin{proof}
    Using \eqref{eqn:Phi_0} and \eqref{eqn:Phi0x} in 
    \eqref{eqn:mathcalP0_hl} and \eqref{eqn:mathcalP0_fi}, 
    we find
    \begin{subequations}
        \begin{align}
            \mathcal P_0^{(\ell)}(k) 
            &= \int_{\D} \frac{\mathfrak P^{(\ell)}(k,y) q_\alpha(y)}{ \sqrt{(\beta\mathfrak n)(k,y)}} \, dy,
        \end{align}
        which gives \eqref{eqn:Phi0_BCs}, using Lemma~\ref{lem:mathfrakPl=0}. Similarly, using \eqref{eqn:Phif_frak} and \eqref{eqn:Phifx} in 
        \eqref{eqn:mathcalPfl_hl} and \eqref{eqn:mathcalPfl_fi}, 
        we find
        \begin{align}
            \mathcal P_f^{(\ell)}(k,t) 
            &= \int_{\D} \frac{\mathfrak P^{(\ell)}(k,y) \mathfrak f_\alpha(k^2,y,t)} { \sqrt{(\beta\mathfrak n)(k,y)}} \,  dy,
        \end{align}
    \end{subequations}
    which gives \eqref{eqn:Phif_BCs}, using Lemma~\ref{lem:mathfrakPl=0}. 
\end{proof}
\begin{lemma} \label{lem:q_BCs}
    For the half-line $(m=0)$ and the finite-interval problem $(m=0,1)$, for $k\in\Omegaext$ and $t\in [0,T]$,
    \begin{subequations} \label{eqn:q_BCs}
        \begin{alignat}{2}
            \label{eqn:q0_BCs}
            \mathcal Q_0^{(\ell)}(t) &= 0,\\
            \label{eqn:qf_BCs}
            \mathcal Q_f^{(\ell)}(t) &= 0,\\
            \label{eqn:Bm_BCs}
            \mathcal Q_{\mathcal B_m}^{(\ell)}(t) &= f_m(t) \tilde \delta_{\ell-1,m},
        \end{alignat}
    \end{subequations}
    where $\tilde \delta_{\ell-1,m}$ is defined in \eqref{eqn:deltatilde}.
\end{lemma}
\begin{proof}
    From Lemmas~\ref{lem:q_0},~\ref{lem:q_f},~and~\ref{lem:qx}, $\mathcal Q_0^{(\ell)}(t)$ \newlink{eqnmathcalQ0}{\eqnmathcalQnaught} and $\mathcal Q_f^{(\ell)}(t)$ \newlink{eqnmathcalQf}{\eqnmathcalQf} are well-defined functions. Similarly, for the {\em regular problems}, $Q_{\mathcal B_m}^{(\ell)}(t)$ \newlink{eqnmathcalQBm}{\eqnmathcalQBm} is a well-defined function from Lemmas~\ref{lem:q_Bm}~and~\ref{lem:qx}. For the {\em irregular problems}, for Boundary Case~\ref{enum:BC3}, $q_{\mathcal B_m,x}(x,t)$ may be undefined at the boundary, but the linear combination of boundary terms $\mathcal Q_{\mathcal B_m}^{(\ell)}(t)$ \eqref{eqn:mathcalQBm_fi} is well defined. For Boundary Case~\ref{enum:BC4}, using Assumption~\ref{ass:alphabetagamma2}, $q_{\mathcal B_m,x}(x,t)$ is well defined at the boundary and therefore $Q_{\mathcal B_m}^{(\ell)}(t)$ is well defined.
    
    For the irregular Boundary Case~\ref{enum:BC3}, see Remark~\enumref{rem:BC}{enum:rem:BC3}, from \eqref{eqn:Bmx_fi}
    for $x\approx x_l$, 
    \begin{align}
        \mathcal B_{2-j,x}(k,x) 
        &= (-1)^j\frac{4k \mathfrak n(k,x) \Eint}{\sqrt{(\beta\mathfrak n)(k,x)}} \left \{  \frac{\beta(x_l) b_{j2} }{\sqrt{(\beta\mathfrak n)(k,x_l)}}  \sum_{n=0}^\infty (-1)^n \mathcal S_{n}^{(x,x_r)}(k) \right\} + \bigoh(k^0).
    \end{align}
    We can prove that either (i) $b_{12} = 0 = b_{22}$, in which case $\mathcal B_{m,x}(k,x_l) = \bigoh(k^0)$, $\mathcal B_{m,x}(k,x_l)/\Delta(k) = \bigoh(k^{-2})$, and $q_{\mathcal B_m,x}(x_l,t)$ is well defined, see Lemma~\ref{lem:qx}; or (ii) if $(b_{12},b_{22}) \neq (0,0)$, then $a_{12} = 0 = a_{22}$, in which case $q_{\mathcal B_m,x}(x,t)$ does not appear in $\mathcal Q_{\mathcal B_m}^{(\ell)}(t)$. The same holds for $x\approx x_r$. It follows that  $\mathcal Q_{\mathcal B_m}^{(\ell)}(t)$ is well defined.
    
    For Boundary Case~\ref{enum:BC4}, with Assumption~\ref{ass:alphabetagamma2}, we integrate $\mathfrak F_m(k^2,t)$ \eqref{eqn:Fmfrak} by parts to obtain
    \begin{align} 
        \mathfrak F_m(k^2,t) &= -\frac{f_m(0)}{k^2}  - \frac{e^{k^2t} f_m'(t) - f_m'(0)}{k^4} + \frac{1}{k^4} \int_0^t e^{k^2s} f_m''(s) \, ds ,
    \end{align}
    so that we may write $q_{\mathcal B_m,x}(x,t)$ \eqref{eqn:q_Bmx} as
    \begin{align}
        q_{\mathcal B_m,x}(x,t) 
        &= \frac{1}{2\pi} \int_{\partial\Omegaext}  \frac{\mathcal B_{m,x} (k,x)} {\Delta(k)} \tilde{\mathfrak F}_m(k^2,t) e^{-k^2t}  \, dk,
    \end{align}
    where 
    \begin{align} 
        \tilde{\mathfrak F}_m(k^2,t) &= -\frac{f_m(0)}{k^2}  + \frac{f_m'(0)}{k^4} + \frac{1}{k^4} \int_0^t e^{k^2s} f_m''(s) \, ds ,
    \end{align}
    and where the integral of the $f_m'(t)$ term is zero by Cauchy's theorem (before the $x$-differentiation). The first two terms of $\tilde{\mathfrak F}_m(k^2,t)e^{-k^2t}$ are exponentially decaying for $t\in(0,T)$ and the last term is $\bigoh(k^{-4})$, by Assumption~\enumref{ass:alphabetagamma2}{enum:fm_extra}. Therefore $q_{\mathcal B_m,x}(x,t)$ is well defined for $x\in \bar \D$ and $t\in(0,T)$. Consequentially, $\mathcal Q_{\mathcal B_m}^{(\ell)}(t)$ is well defined.
    
    Using \eqref{eqn:q_0_deformed} and \eqref{eqn:q0x} in 
    \newlink{eqnmathcalQ0}{\eqnmathcalQnaught}, we find
    \begin{subequations}
        \begin{align}
            \mathcal Q_0^{(\ell)} (t)  &=\frac{1}{2\pi} \int_{\partial\Omegaext} \frac{\mathcal P_0^{(\ell)}(k)}{\Delta(k)} e^{-k^2t} \, dk,
        \end{align}
        which gives \eqref{eqn:q0_BCs}, using Lemma~\ref{lem:Phi_BCs}.
        Similarly, using \eqref{eqn:q_f_deformed2} and \eqref{eqn:qfx} in 
        \newlink{eqnmathcalQf}{\eqnmathcalQf}, we find
        \begin{align}
            \mathcal Q_{f}^{(\ell)} (t) 
            &= \frac{1}{2\pi} \int_{\partial\Omegaext} \frac{\mathcal P_f^{(\ell)}(k,t) e^{-k^2t}}{\Delta(k)}   dk.
        \end{align}
        Using Lemma~\ref{lem:Phi_BCs}, this gives \eqref{eqn:qf_BCs}.
        Using \eqref{eqn:q_Bm_deformed2} and \eqref{eqn:q_Bmx} in \newlink{eqnmathcalQBm}{\eqnmathcalQBm}, 
        \begin{align}
            \mathcal Q_{\mathcal B_m}^{(\ell)} (t) 
            = \frac{1}{2\pi} \int_{\partial\Omegaext} \frac{\mathfrak B_m^{(\ell)}(k)}{\Delta(k)} \mathfrak F_m(k^2,t) e^{-k^2t} \, dk.
        \end{align}
        Finally, using Lemma~\ref{lem:mathfrakBml_eqn} and \eqref{eqn:Fmfrak}, we obtain
        \begin{align}
            \mathcal Q_{\mathcal B_m}^{(\ell)} (t) 
            &= -\frac{\tilde \delta_{\ell-1,m}}{i\pi} \int_{\partial\Omegaext} \left( \frac{f_m(0)e^{-k^2t}}{k} + \frac{1}{k} \int_0^t e^{-k^2(t-s)} f_m'(s) \, ds \right) dk.
        \end{align}
        Since the integrand is $\bigoh(k^{-3})$, we can deform the path of integration to the real axis. Using the oddness of the integrand, the principal value integral vanishes and only the residue contribution at the origin needs to be calculated:
        \begin{align}
            \mathcal Q_{\mathcal B_m}^{(\ell)} (t) 
            &= \tilde \delta_{\ell-1,m} \Res \left( \frac{f_m(0)e^{-k^2t}}{k} + \frac{1}{k} \int_0^t e^{-k^2(t-s)} f_m'(s) \, ds; \, k=0 \right) = f_m(t) \tilde \delta_{\ell-1,m}.
        \end{align}    
    \end{subequations}
\end{proof}
\begin{lemma} \label{lem:BC_limit}
    Consider any $t\in(0,T)$, fixed. Then 
    \begin{align} \label{eqn:BC_limit}
        \lim_{|x|\to\infty} q(x,t) = 0
        \qquad \text{ and } \qquad
        \lim_{x\to\infty} q(x,t) = 0,
    \end{align}
    for the whole-line and half-line problems, respectively.
\end{lemma}
\begin{proof}
    For any fixed $t\in(0,T)$, we have absolute integrability in \eqref{eqn:q_Bm_deformed2}, \eqref{eqn:q_0_deformed}, and \eqref{eqn:q_f_deformed2_direct}. Therefore, we may switch the limit and integrals. Since, from \eqref{eqn:Psi_bound_wl} and \eqref{eqn:Psi_bound_hl}, 
    \begin{align}
        \lim_{|x|\to \infty} \left| \frac{\Psi(k,x,y)}{\Delta(k)} \right| \leq M_\Psi e^{-m_{i\nu}|k||x-y|} = 0,
        &&
        \lim_{x\to \infty} \left| \frac{\Psi(k,x,y)}{\Delta(k)} \right| \leq M_\Psi e^{-m_{i\nu}|k||x-y|} = 0,
    \end{align}
    for the whole-line problem and the half-line problem, respectively,  \eqref{eqn:BC_limit} follows.
\end{proof}
\begin{rem}
    Since we have absolute integrability in \eqref{eqn:q0x}, \eqref{eqn:qfx}, and in \eqref{eqn:q_Bmx}, we conclude that also
    \begin{align} \label{eqn:BCx_limit}
        \lim_{|x|\to\infty} q_x(x,t) = 0
        \qquad \text{ and } \qquad
        \lim_{x\to\infty} q_x(x,t) = 0,
    \end{align}
    for the whole-line and half-line problems, respectively.
\end{rem}
\begin{theorem}
    Consider the finite-interval, the half-line, and the whole-line problems. For all three problems, the solution expression \eqref{eqn:q_all} satisfies the appropriate boundary conditions.
\end{theorem}
\begin{proof}
    Lemma~\ref{lem:BC_limit} shows the boundary conditions for the whole-line problem and the right boundary condition for the half-line problem are satisfied. From Lemma~\ref{lem:q_BCs}, 
    \begin{subequations}
        \begin{align}
            a_{0} q(x_l,t) +a_{1} q(x_l,t) = \mathcal Q_0^{(0)}(t) + \mathcal Q_f^{(0)}(t) + \mathcal Q_{\mathcal B_m}^{(0)}(t) = f_0(t).
        \end{align}
        Similarly, for the finite-interval problem, 
        \begin{align} \nonumber
            a_{\ell1} q(x_l,t) +a_{\ell2} q_{x}(x_l,t) + b_{\ell1} q(x_r,t) + b_{\ell2} q_{x}(x_r,t) 
            &= \mathcal Q_0^{(\ell)}(t) + \mathcal Q_f^{(\ell)}(t) + \mathcal Q_{\mathcal B_0}^{(\ell)}(t) + \mathcal Q_{\mathcal B_1}^{(\ell)}(t) \\
            &= f_0(t) \tilde \delta_{\ell-1,0} + f_1(t) \tilde \delta_{\ell-1,1}.
        \end{align}
    \end{subequations}
\end{proof}
 \label{sec:proofs_BC}

\section{Proofs: the solution expressions satisfy the initial condition} \label{sec:proofs_IC}
\begin{theorem} \label{thm:qf,qBmlimits}
    Consider the finite-interval, half-line, and whole-line problms. For $x\in \D$, fixed, 
    \begin{subequations}
        \begin{align}
            \label{eqn:qf_limit}
            \lim_{t\to 0^+} q_f(x,t) &= 0, \\
            \label{eqn:qBm_limit}
            \lim_{t\to 0^+} q_{\mathcal B_m} (x,t) &= 0.
        \end{align}
    \end{subequations}
\end{theorem}
\begin{proof}
    Since the integral in \eqref{eqn:q_f_deformed2} is absolutely convergent, we can pass the limit $t\to0^+$ inside the integral to obtain \eqref{eqn:qf_limit} by using Cauchy's theorem. Similarly, we move the limit in the integral in \eqref{eqn:q_Bm_deformed2} to obtain \eqref{eqn:qBm_limit}.
\end{proof}
\begin{lemma} \label{lem:Psi/Delta_asymptotics}
    For fixed $x\in \D$, for $y\in \bar \D$, for the finite-interval, half-line, and whole-line problems, 
    \begin{align} \label{eqn:Psi/Delta_asymptotics} 
        \frac{\Psi(k,x,y)}{\Delta(k)} 
        &=  \exp\left( \sgn(x-y) ik \int_y^x \mathfrak n(k,\xi) \, d\xi \right) \big(1 + \littleoh(k^0) \big) 
        + \littleoh(k^{-1}),
    \end{align}
    as $|k|\to \infty$ for $k\in\Omegaext$.
\end{lemma}
\begin{proof}
    For the whole-line problem, from \eqref{eqn:Psi_wl}, for $y<x$,
    \begin{align} \label{eqn:Psi_asym_wl_intermediate}
        \Psi(k,x,y) &= \etwo{\sgn(x-y)ik\int_y^x \mathfrak n(k,\xi) \, d\xi}         \left(1 + \sum_{n=1}^\infty \sum_{\ell=0}^n (-1)^\ell \tilde{\mathcal E}_{n-\ell}^{(-\infty,y)}(k) \mathcal E_\ell^{(x,\infty)}(k) \right)\!. 
    \end{align}
    By Lemma~\ref{lem:Jnasymptotics} and the DCT,
    \begin{align}
        \sum_{n=1}^\infty \sum_{\ell=0}^n (-1)^\ell \tilde{\mathcal E}_{n-\ell}^{(-\infty,y)}(k) \mathcal E_\ell^{(x,\infty)}(k) = \littleoh(k^0).
    \end{align}
    Dividing \eqref{eqn:Psi_asym_wl_intermediate} by $\Delta(k)$ and using Lemma~\ref{lem:Delta_asymptotics}, we obtain \eqref{eqn:Psi/Delta_asymptotics}. The proof for $x<y$ is identical.
    
    For the half-line problem, for $x_l<y<x$, we write \eqref{eqn:Psi_hl} as
    \begin{align} \label{eqn:Psi_intermediate_hl} 
        \Psi(k,x,y) 
        &= 4 \etwo{ik \int_y^x \mathfrak n(k,\xi) \, d\xi}  \left[ \left( \frac{a_0}{k \mathfrak n(k,x_l)} \mathscr S_0^{(x_l,y)}(k) - a_1 \mathscr C_0^{(x_l,y)}(k) \right)+ \left( \frac{|a_0|}{m_{\mathfrak n}|k|} + |a_1| \right)\littleoh(k^0) \right] \!.
    \end{align}
    Using
    \begin{subequations} \label{eqn:C0S0_identity}
        \begin{align} 
            \label{eqn:C0_identity}
            \mathscr C_0^{(a,b)}(k)
            &= \frac{1}{2} \left( \exp\left( 2ik \int_{a}^b \mathfrak n(k,\xi) \, d\xi \right) + 1\right)\!, \\
            \label{eqn:S0_identity}
            \mathscr S_0^{(a,b)}(k) 
            &= \frac{1}{2i} \left( \exp\left( 2ik \int_{a}^b \mathfrak n(k,\xi) \, d\xi \right) - 1 \right)\!,
        \end{align}
    \end{subequations}
    in \eqref{eqn:Psi_intermediate_hl}, we find
    \begin{align} 
        \Psi(k,x,y) &= 2 \etwo{ik\int_y^x \mathfrak n(k,\xi) \, d\xi} 
        \left[ \frac{ia_0}{k\mathfrak n(k,x_l)} - a_1 + \left( \frac{|a_0|}{m_{\mathfrak n}|k|} + |a_1| \right)\littleoh(k^0) \right] 
        +\bigoh\big(\eone{-m_{i\mathfrak n}|k|(x-x_l)} \big),
    \end{align}
    which, from \eqref{eqn:Delta_asym} with \eqref{eqn:b0_hl}, gives \eqref{eqn:Psi/Delta_asymptotics}. The proof is identical for $x_l<x<y$.
    
    For the finite-interval problem we consider the 4 different cases. 
    \begin{enumerate}
        \item If $(a:b)_{2,4} \neq 0$, then
        for $x_l<y<x<x_r$, using \eqref{eqn:Delta_asym} with \eqref{eqn:b0_1}, we write \eqref{eqn:Psi_fi1} as
        \begin{align} \nonumber
            \frac{\Psi(k,x,y)}{\mathfrak b_0(k)} &= \frac{-4}{(a:b)_{2,4}} \etwo{ ik \int_y^x \mathfrak n(k,\xi) \, d\xi} \left\{  
            -(a:b)_{2,4} 
            \mathscr C_{0}^{(x_l,y)}(k) \mathscr C_0^{(x,x_r)}(k) 
            + \littleoh(k^0) \right\} \\
            &~~~ + \frac{4\beta(x_r)(a:b)_{1,2}}{(a:b)_{2,4} k \sqrt{(\beta \mathfrak n)(k,x_l)} \sqrt{(\beta \mathfrak n)(k,x_r)}} \sum_{n=0}^\infty \Xi(k) \mathcal S_n^{(y,x)}(k).
        \end{align}
        Using \eqref{eqn:C0_identity} and dividing by $1+\varepsilon(k)$, we arrive at \eqref{eqn:Psi/Delta_asymptotics}.
        The proof for $x_l<x<y<x_r$ is identical.
        \item If $(a:b)_{2,4} = 0$ and $m_{\mathfrak c_0} \neq 0$, then 
        for $x_l<y<x<x_r$, using \eqref{eqn:C0S0_identity} and \eqref{eqn:Delta_asym} with \eqref{eqn:b0_2}, we write \eqref{eqn:Psi_fi1} as
        \begin{align}
            \frac{\Psi(k,x,y)}{\mathfrak b_0(k)} 
            &= \frac{4k}{im_{\mathfrak c_0}} \etwo{ ik \int_y^x \mathfrak n(k,\xi) \, d\xi} \left\{  - \frac1{4i} \left(\frac{(a:b)_{1,4}}{k \mathfrak n(k,x_l)} + \frac{(a:b)_{2,3}}{k \mathfrak n(k,x_r)}  \right) + \littleoh (k^{-1}) \right\}  + \littleoh( k^{-1} ).
        \end{align}
        Using $1/\mathfrak n(k,x) = 1/\mu(x) + \bigoh(k^{-2})$ and dividing by $1+\varepsilon(k)$, we obtain \eqref{eqn:Psi/Delta_asymptotics}. The proof for $x_l<x<y<x_r$ is identical.
        \item If $(a:b)_{2,4} = 0$, $m_{\mathfrak c_0} = 0$, $m_{\mathfrak c_1}=0$, and $(a:b)_{1,3} \neq 0$, then 
        for $x_l<y<x<x_r$, using \eqref{eqn:Delta_asym} with \eqref{eqn:b0_3a}, we write \eqref{eqn:Psi_fi1} as
        \begin{align}
            \frac{\Psi(k,x,y)}{\mathfrak b_0(k)} &= -\frac{4k^2}{m_{\mathfrak s}} \etwo{ ik \int_y^x \mathfrak n(k,\xi) \, d\xi} \left\{ -\frac{1}{4}  \frac{(a:b)_{1,3}}{k^2 \mathfrak n(k,x_l) \mathfrak n(k,x_r)}  + \littleoh (k^{-2}) \right\} + \littleoh(k^{-1}).
        \end{align}
        Using the asymptotics for $1/\mathfrak n(k,x)$ and dividing by $1+\varepsilon(k)$ gives \eqref{eqn:Psi/Delta_asymptotics}. The proof is identical for $x_l<x<y<x_r$.
        \item If $(a:b)_{2,4} = 0$, $m_{\mathfrak c_0} = 0$, $m_{\mathfrak c_1}\neq 0$, and $m_{\mathfrak c_1} \mathfrak u_+ - 8m_{\mathfrak s} \neq 0$, then 
        for $x_l<y<x<x_r$, using that
        \begin{align}
            \sum_{n=3}^\infty \left| k \mathscr C_n^{(a,b)}(k)\right| = \bigoh(k^{-1}),
        \end{align}
        using the asymptotics of $1/\mathfrak n(k,x)$, the fact that $(a:b)_{1,4}/\mu(x_r) =(a:b)_{2,3}/\mu(x_l) =  m_{\mathfrak c_1}/2$, and  \eqref{eqn:Delta_asym} with \eqref{eqn:b0_3b}, we write \eqref{eqn:Psi_fi1} as
        \begin{align} \nonumber
            \frac{\Psi(k,x,y)}{\mathfrak b_0(k)} 
            &= \frac{32k^2}{m_{\mathfrak c_1} \mathfrak u_+ - 8m_{\mathfrak s}} \etwo{ ik \int_y^x \mathfrak n(k,\xi) \, d\xi} \left\{  \frac{m_{\mathfrak c_1}}{2k} \sum_{n=1}^2 \sum_{\ell=0}^n  (-1)^\ell \mathscr S_{n-\ell}^{(x_l,y)}(k) \mathscr C_\ell^{(x,x_r)}(k)\right. \\
            &~~~\qquad \left. - \frac{m_{\mathfrak c_1}}{2k} \sum_{n=1}^2 \sum_{\ell=0}^n \mathscr C_{n-\ell}^{(x_l,y)}(k) \mathscr S_\ell^{(x,x_r)}(k) + \littleoh(k^{-2}) - \frac{m_{\mathfrak s}}{4k^2} \right\} + \littleoh (k^{-1}).
        \end{align}
        Using integration by parts as in Lemma~\ref{lem:Jnasymptotics}, we derive
        \begin{subequations}
            \begin{align}
                \mathscr C_1^{(a,b)}(k) &= \frac{1}{16ik} \mathfrak u_+(a,b) \left( \etwo{ 2ik \int_a^b \mathfrak n(k,\xi) \, d\xi} - 1\right) + \littleoh(k^{-1}), \\
                \mathscr S_1^{(a,b)}(k) &= -\frac{1}{16k} \mathfrak u_-(a,b) \left( \etwo{ 2ik \int_a^b \mathfrak n(k,\xi) \, d\xi} + 1\right) + \littleoh(k^{-1}), \\
                \mathscr C_2^{(a,b)}(k) &= \frac{1}{64ik} m_\mathrm{int}(a,b) \left( \etwo{ 2ik \int_a^b \mathfrak n(k,\xi) \, d\xi} - 1\right) + \littleoh(k^{-1}), \\
                \mathscr S_2^{(a,b)}(k) &= -\frac{1}{64k} m_\mathrm{int}(a,b) \left( \etwo{ 2ik \int_a^b \mathfrak n(k,\xi) \, d\xi} + 1\right) + \littleoh(k^{-1}), 
            \end{align}
        \end{subequations}
        where $\mathfrak u_\pm (a,b) = \mathfrak u(b) \pm \mathfrak u(a)$, and
        \begin{align}
            m_\mathrm{int}(a,b) 
            = \int_a^b \frac{1}{\mu(y)}\left( \frac{(\beta\mu)'(y)}{(\beta\mu)(y)} \right)^2 dy,
        \end{align}
        with $\mathfrak u(x)$ defined in \eqref{eqn:mathfrak_u_def}.
        We find
        \begin{align} \nonumber
            \frac{\Psi(k,x,y)}{\mathfrak b_0(k)} &= \frac{32k^2}{m_{\mathfrak c_1} \mathfrak u_+ - 8m_{\mathfrak s}} \etwo{ ik \int_y^x \mathfrak n(k,\xi) \, d\xi} \left\{ - \frac{m_{\mathfrak s}}{4k^2} + \frac{m_{\mathfrak c_1}}{64k^2} \left( \mathfrak u_+(x_l,y) - \mathfrak u_-(x_l,y) \right) \right. \\
            &\hspace{2in} \left. + \frac{m_{\mathfrak c_1}}{64k^2} \left( \mathfrak u_+(x,x_r) +  \mathfrak u_-(x,x_r) \right)  + \littleoh(k^{-2})\right\}+ \littleoh (k^{-1}).
        \end{align}
        Combining terms, 
        \begin{align}
            \frac{\Psi(k,x,y)}{\mathfrak b_0(k)} &= \frac{32k^2}{m_{\mathfrak c_1} \mathfrak u_+ - 8m_{\mathfrak s}} \etwo{ ik \int_y^x \mathfrak n(k,\xi) \, d\xi} \left\{ - \frac{m_{\mathfrak s}}{4k^2} + \frac{m_{\mathfrak c_1}}{32k^2} \mathfrak u_+  + \littleoh (k^{-2})\right\}+ \littleoh (k^{-1}),
        \end{align}
        which, after dividing by $1+\varepsilon(k)$, gives \eqref{eqn:Psi/Delta_asymptotics}. The proof is identical for $x_l<x<y<x_r$.
    \end{enumerate}
\end{proof}
\begin{theorem} \label{thm:q0limit}
    Consider the finite-interval, half-line, and whole-line problems.  If $q_0\in L^1(\D)$, then for almost every $x\in \D$,
    \begin{align} \label{eqn:q0_limit}
        \lim_{t\to 0^+} q_0(x,t) &= q_0(x).
    \end{align}
\end{theorem}
\begin{proof}
    Using the change of variables $k = \lambda z$ with $\lambda = 1/\sqrt{t}$ in \eqref{eqn:q_0_deformed}, 
    \begin{align} 
        q_0(x,t) &= \frac{\lambda}{2\pi} \int_{\partial\Omegaext} \frac{\Phi_0(\lambda z,x)}{{\Delta(\lambda z)}} e^{-z^2} \, dz 
        = \frac{\lambda}{2\pi} \int_{\partial\Omegaext} \frac{ e^{-z^2} }{{\Delta(\lambda z)}} \int_\D \frac{\Psi(\lambda z,x,y) q_\alpha(y)}{\sqrt{(\beta\mathfrak n)(\lambda z,x)}\sqrt{(\beta\mathfrak n)(\lambda z,y)}} \, dy \, dz.
    \end{align}
    By Lemma~\ref{lem:Psi_bounds}, we can use the Fubini-Tonelli theorem to write this as
    \begin{align} 
        q_0(x,t) &= \frac{\lambda}{2\pi} \int_\D q_\alpha(y) \int_{\partial\Omegaext} \frac{\Psi(\lambda z,x,y)}{{\Delta(\lambda z)}}   \frac{ e^{-z^2} }{\sqrt{(\beta\mathfrak n)(\lambda z,x)}\sqrt{(\beta\mathfrak n)(\lambda z,y)}} \, dz \, dy.
    \end{align}
    Using \eqref{eqn:Psi/Delta_asymptotics},
    \begin{align} \label{eqn:q_0_limit_intermediate}
        q_0(x,t) &= \frac{\lambda\big(1 + \littleoh(\lambda^0) \big)}{2\pi} \int_\D \frac{q_\alpha(y)}{\sqrt{(\beta\mu)(x)}\sqrt{(\beta\mu)(y)}} \int_{\partial \Omegaext} \exp\left( \sgn(x-y) i\lambda z \int_y^x \mathfrak n(\lambda z,\xi) \, d\xi \right) e^{-z^2}  \, dz \, dy+ \littleoh(\lambda^{0})  .
    \end{align}
    Since
    \begin{align} 
        \left| \lambda \exp\left( \sgn(x-y) i\lambda z \int_y^x \mu(\xi) \, d\xi \right) \bigoh(|x-y| \lambda^{-1}) e^{-z^2} \right| \leq \bigoh(|x-y|\lambda^0) \big|e^{-z^2}\big|,
    \end{align}
    is absolutely integrable, we may use the DCT on the remainder term from \eqref{eqn:E_approx}. Substituting this result in \eqref{eqn:q_0_limit_intermediate}, we obtain  
    \begin{align} 
        q_0(x,t) &= \frac{\lambda}{2\pi} \int_\D \frac{q_\alpha(y)}{\sqrt{(\beta\mu)(x)}\sqrt{(\beta\mu)(y)}} \int_{\partial\Omegaext} \exp\left( \sgn(x-y) i\lambda z \int_y^x \mu(\xi) \, d\xi \right) e^{-z^2}  \, dz \, dy + \littleoh(\lambda^{0}),
    \end{align}
    as $\lambda \to \infty$.
    Define $M_x(y) = \int_x^y \mu(\xi) \, d\xi$. Deforming $\partial \Omegaext$ down to the real axis and  integrating the $z$-integral gives
    \begin{align} 
        q_0(x,t) &= \frac{\lambda}{2\sqrt{\pi}} \int_\D \frac{q_\alpha(y)\eone{-\frac{1}{4}\lambda^2 M_x^2(y)}}{\sqrt{(\beta\mu)(x)}\sqrt{(\beta\mu)(y)}}  \, dy+ \littleoh(\lambda^{0}).
    \end{align}
    For a fixed $x \in \D$, if $q_0(x)$ is finite, using that $q_\alpha(y)/(\mu(y) \sqrt{(\beta\mu)(y)}) \in L^1(\D)$, it follows that for any $\epsilon>0$ and for each $\lambda$, there exists $\varphi \in \mathrm{AC}(\D) \cap C_0(\D)$ \cite{folland}, so that
    \begin{align}
        \int_\D \bigg| \frac{q_\alpha(y)}{\mu(y) \sqrt{(\beta\mu)(y)}} - \varphi(y) \bigg| \, dy \leq \lambda^{-2}, 
        \qquad \text{ and } \qquad
        \left|\varphi(x) - \frac{q_\alpha(x)}{\mu(x) \sqrt{(\beta\mu)(x)}} \right| \leq \frac{\epsilon}{2}.
    \end{align}
    Using this,
    \begin{align} \label{eqn:q_0_limit_intermediate2}
        q_0(x,t) &= \frac{\lambda}{2\sqrt{\pi}\sqrt{(\beta\mu)(x)}} \left[ \int_\D \left( \frac{q_\alpha(y)}{\mu(y) \sqrt{(\beta\mu)(y)}} -\varphi(y) \right) \mu(y) \eone{-\frac{1}{4}\lambda^2 M_x^2(y)} \, dy + \int_\D \varphi(y) \mu(y) \eone{-\frac{1}{4}\lambda^2 M_x^2(y)} \, dy \right]+ \littleoh(\lambda^{0}).
    \end{align}
    For the first integral and any $\epsilon>0$, we can find $\lambda$ sufficiently large, so that
    \begin{align}
        \left| \frac{\lambda}{2\sqrt{\pi}\sqrt{(\beta\mu)(x)}}\int_\D \left( \frac{q_\alpha(y)}{\mu(y) \sqrt{(\beta\mu)(y)}} -\varphi(y) \right) \mu(y) \eone{-\frac{1}{4}\lambda^2 M_x^2(y)} \, dy \right| &\leq \frac{M_{\mathfrak n}}{2\sqrt{\pi m_\beta m_{\mathfrak n}}\lambda} \leq \frac{\epsilon}{2}.
    \end{align}
    Since $\varphi \in \mathrm{AC}(\D)$, we may integrate the second integral of \eqref{eqn:q_0_limit_intermediate2} by parts to obtain
    \begin{align}
       \frac{\lambda}{2\sqrt{\pi}} \int_\D \varphi(y) \mu(y) \eone{-\frac{1}{4}\lambda^2 M_x^2(y)} \, dy 
       &= - \frac12  \int_\D \varphi'(y)\erf\left(\frac{\lambda M_x(y)}{2} \right) dy .
    \end{align}
    \no At this point, we may take the limit as $\lambda \to \infty$ using the DCT. Since $\arg(\mu) \in (-\pi/4, \pi/4)$, if $y>x$, then $\arg(M_x(y)) \in (-\pi/4,\pi/4)$ and the error function limits to $1$  as $\lambda \to \infty$ \cite{dlmf}. If $y<x$, then \mbox{$\arg(M_x(y)) \in (3\pi/4,5\pi/4)$} and the error function limits to $-1$. It follows that
    \begin{align}
      \lim_{\lambda\to \infty} \frac{\lambda}{2\sqrt{\pi}} \int_\D \varphi(y) \mu(y) \eone{-\frac{1}{4}\lambda^2 M_x^2(y)} \, dy
       &=\varphi(x),
    \end{align}
    and we have
    \begin{align}
        q_0(x,t) \to \frac{\varphi(x)}{\sqrt{\beta\mu(x)}} \to q_0(x),
    \end{align}
    as $t\to 0^+$ and $\epsilon\to 0^+$. Since $q_0 \in L^1(\D)$, $q_0(x)$ is finite for almost every $x\in \D$, concluding the proof.
\end{proof}

\end{appendices}
\section*{Acknowledgements}
The authors thank David Smith and Vishal Vasan for useful conversations. 
%
%
{\small
\bibliographystyle{abbrv}
\bibliography{references}
}
\end{document}